\documentclass{birkart}
\usepackage{amsmath}
\usepackage{amsxtra}
\usepackage{amscd}
\usepackage{amsthm}
\usepackage{amsfonts}
\usepackage{amssymb}
\usepackage{eucal}

\newtheorem{thm}{Theorem}[section]
\newtheorem{conj}{Conjecture}
\newtheorem{cor}[thm]{Corollary}
\newtheorem{lem}[thm]{Lemma}
\newtheorem{prop}[thm]{Proposition}

\theoremstyle{remark}
\newtheorem{remark}[thm]{Remark}

\theoremstyle{definition}

\numberwithin{equation}{section}

\newcommand{\bean}{\begin{eqnarray}}
\newcommand{\eean}{\end{eqnarray}}
\newcommand{\be}{\begin{displaymath}}
\newcommand{\ee}{\end{displaymath}}
\newcommand{\bea}{\begin{eqnarray*}}
\newcommand{\eea}{\end{eqnarray*}}

\newcommand{\thmref}[1]{Theorem~\ref{#1}}
\newcommand{\secref}[1]{Section~\ref{#1}}
\newcommand{\lemref}[1]{Lemma~\ref{#1}}
\newcommand{\propref}[1]{Proposition~\ref{#1}}
\newcommand{\corref}[1]{Corollary~\ref{#1}}

\newcommand{\conjref}[1]{Conjecture~\ref{#1}}
\newcommand{\nc}{\newcommand}
\nc{\on}{\operatorname}
\nc{\ch}{\mbox{ch}}
\nc{\Z}{{\mathbb Z}}
\nc{\C}{{\mathbb C}}
\nc{\pone}{{\mathbb P}^1}
\nc{\pa}{\partial}
\nc{\F}{{\mathcal F}}
\nc{\arr}{\rightarrow}
\nc{\larr}{\longrightarrow}
\nc{\al}{\alpha}
\nc{\ri}{\rangle}
\nc{\lef}{\langle}
\nc{\W}{{\mathcal W}}
\nc{\la}{\lambda}
\nc{\ep}{\epsilon}
\nc{\su}{\widehat{{\mathfrak s}{\mathfrak l}}_2}
\nc{\sw}{{\mathfrak s}{\mathfrak l}}
\nc{\g}{{\mathfrak g}}
\nc{\h}{{\mathfrak h}}
\nc{\n}{{\mathfrak n}}
\nc{\N}{\widehat{\n}}
\nc{\G}{\widehat{\g}}
\nc{\De}{\Delta}
\nc{\gt}{\widetilde{\g}}
\nc{\Ga}{\Gamma}
\nc{\one}{{\mathbf 1}}
\nc{\z}{{\mathfrak Z}}
\nc{\La}{\Lambda}
\nc{\wt}{\widetilde}
\nc{\wh}{\widehat}
\nc{\cri}{_{\kappa_c}}
\nc{\kk}{h^\vee}
\nc{\sun}{\widehat{\sw}_N}
\nc{\si}{\sigma}
\nc{\el}{\ell}
\nc{\bi}{\bibitem}
\nc{\om}{\omega}
\nc{\ol}{\overline}
\nc{\ds}{\displaystyle}
\nc{\dzz}{\frac{dz}{z}}
\nc{\Res}{\on{Res}}
\nc{\mc}{\mathcal}
\nc{\Cal}{\mathcal}
\nc{\bb}{{\mathfrak b}}
\nc{\ot}{\otimes}
\nc{\R}{{\mc R}}
\nc{\yy}{{\mc Y}}
\nc{\ga}{\gamma}

\nc{\us}{\underset}
\nc{\opl}{\oplus}
\nc{\beq}{\begin{equation}}
\nc{\Fq}{{\mathcal F}}
\nc{\Mq}{{\mathcal M}}
\nc{\Rep}{\on{Rep}}
\nc{\sssec}{\subsubsection}
\nc{\ssec}{\subsection}
\nc{\lan}{\langle}
\nc{\ran}{\rangle}

\nc{\D}{\mathcal D}
\nc{\Vect}{\on{Vect}}
\nc{\ghat}{\G}
\nc{\T}{\mc T}
\nc{\Tloc}{\T^\g_{\on{loc}}}
\nc{\vac}{|0\ran}
\nc{\Wick}{{\mb :}}
\nc{\mb}{\mathbf}
\nc{\delz}{\partial_z}
\nc{\K}{{\cali K}}
\nc{\cali}{\mathcal}
\nc{\li}{\mathfrak l}
\nc{\lt}{\widetilde{\li}}
\nc{\astar}{a^*}
\nc{\cA}{{\mc A}}
\nc{\ka}{\kappa}

\nc{\OO}{{\mc O}}
\nc{\AutO}{\on{Aut}\OO}
\nc{\DerO}{\on{Der}\OO}
\nc{\DerpO}{\on{Der}_+\OO}
\nc{\Au}{{\mc A}ut}
\nc{\mf}{\mathfrak}
\nc{\V}{{\mathbb V}}
\nc{\hh}{\wh{\h}}

\nc{\pp}{{\mathfrak p}}
\nc{\mm}{{\mathfrak m}}
\nc{\rr}{{\mathfrak r}}
\nc{\ket}{\rangle}
\nc{\zz}{{\mathfrak z}}
\nc{\gr}{\on{gr}}
\nc{\Spe}{\on{Spec}}
\nc{\rv}{\crho}
\nc{\can}{\on{can}}
\nc{\CC}{\on{Op}_G(D))}
\nc{\Op}{\on{Op}_G(D)}
\nc{\MOp}{\on{MOp}_G(D)}
\nc{\Db}{{\mathbb D}}
\nc{\ww}{w}

\nc{\af}{{\mathbb A}^1}
\nc{\bs}{\backslash}
\nc{\laa}{(\la_i)}
\nc{\zn}{(z_i)}

\nc{\cla}{\check{\la}}
\nc{\cmu}{\check{\mu}}
\nc{\crho}{\check{\rho}}
\nc{\chal}{\check{\al}}
\nc{\cc}{{\mathfrak c}}

\nc{\M}{{\mathbb M}}

\nc{\ZZ}{{\mc Z}_{(z_i)}}

\nc{\UU}{{\mathbb U}}

\nc{\Conn}{\on{Conn}(\Omega^{\crho})}
\nc{\Con}{\on{Conn}(\Omega^{-\rho})}
\nc{\Co}{\on{Conn}(\Omega^{\rho})}

\begin{document}

\title{Gaudin model and opers}

\author{Edward Frenkel}\thanks{Partially supported by grants from the
NSF and DARPA}

\address{Department of Mathematics, University of California,
  Berkeley, CA 94720, USA}

\begin{abstract}

This is a review of our previous works \cite{FFR,F:icmp,F:opers}
(some of them joint with B. Feigin and N. Reshetikhin) on the
Gaudin model and opers. We define a commutative subalgebra in the
tensor power of the universal enveloping algebra of a simple Lie
algebra $\g$. This algebra includes the hamiltonians of the Gaudin
model, hence we call it the Gaudin algebra. It is constructed as a
quotient of the center of the completed enveloping algebra of the
affine Kac-Moody algebra $\ghat$ at the critical level. We
identify the spectrum of the Gaudin algebra with the space of
opers associated to the Langlands dual Lie algebra $^L \g$ on the
projective line with regular singularities at the marked points.
Next, we recall the construction of the eigenvectors of the Gaudin
algebra using the Wakimoto modules over $\ghat$ of critical level.
The Wakimoto modules are naturally parameterized by Miura opers
(or, equivalently, Cartan connections), and the action of the
center on them is given by the Miura transformation. This allows
us to relate solutions of the Bethe Ansatz equations to Miura
opers and ultimately to the flag varieties associated to the
Langlands dual Lie algebra $^L \g$.
\end{abstract}

\maketitle

\section*{Introduction}

Let $\g$ be a finite-dimensional simple Lie algebra over $\C$ and
$U(\g)$ its universal enveloping algebra. Choose a basis $\{ J_a \},
a=1,\ldots,d$, of $\g$, and let $\{ J^a \}$ the dual basis with respect
to a non-degenerate invariant bilinear form on $\g$. Let
$z_1,\ldots,z_N$ be a collection of distinct complex numbers.

The {\em Gaudin Hamiltonians} are the following elements of the
algebra $U(\g)^{\otimes N}$:
\begin{equation} \label{the gaudin ham}
\Xi_i = \sum_{j\neq i} \sum_{a=1}^d \frac{J_a^{(i)}
J^{a(j)}}{z_i-z_j}, \qquad i=1,\ldots,N,
\end{equation}
where for $A \in \g$ we denote by $A^{(i)}$ the element of
$U(\g)^{\otimes N}$ which is the tensor product of $A$ in the $i$th
factor and $1$ in all other factors. One checks easily that these
elements commute with each other and are invariant with respect to the
diagonal action of $\g$ on $U(\g)^{\otimes N}$.

For any collection $M_1,\ldots,M_N$ of $\g$-modules, the Gaudin
Hamiltonians give rise to commuting linear operators on
$\bigotimes_{i=1}^N M_i$. We are interested in the diagonalization of
these operators. More specifically, we will consider the following two
cases: when all of the $M_i$'s are Verma modules and when they are
finite-dimensional irreducible modules.

It is natural to ask: are there any other elements in
$U(\g)^{\otimes N}$ which commute with the Gaudin Hamiltonians?
Clearly, the $N$-fold tensor product $Z(\g)^{\otimes N}$ of the
center $Z(\g)$ of $U(\g)$ is the center of $U(\g)^{\otimes N}$,
and its elements obviously commute with the $\Xi_i$'s. As shown in
\cite{FFR}, if $\g$ has rank grater than one, then in addition to
the Gaudin Hamiltonians and the central elements there are other
elements in $U(\g)^{\otimes N}$ of orders higher than two which
commute with the Gaudin operators and with each other (but
explicit formulas for them are much more complicated and unknown
in general). Adjoining these ``higher Gaudin Hamiltonians'' to the
$\Xi_i$'s together with the center $Z(\g)^{\otimes N}$, we obtain
a large commutative subalgebra of $U(\g)^{\otimes N}$. We will
call it the {\em Gaudin algebra} and denote it by $\ZZ(\g)$.

The construction of $\ZZ(\g)$ will be recalled in \secref{higher
gaudin}. The key point is the realization of $U(\g)^{\otimes N}$
as the space of {\em coinvariants} of induced modules over the
affine Kac-Moody algebra $\ghat$ on the projective line. Using
this realization, we obtain a surjective map from the center of
the completed universal enveloping algebra of $\ghat$ at the {\em
critical level} onto $\ZZ(\g)$.

The next natural question is what is the spectrum of $\ZZ(\g)$, i.e.,
the set of all maximal ideals of $\ZZ(\g)$, or equivalently, algebra
homomorphisms $\ZZ(\g) \to \C$. Knowing the answer is important,
because then we will know how to think about the common eigenvalues of
the higher Gaudin operators on the tensor products $\bigotimes_{i=1}^N
M_i$ of $\g$-modules. These common eigenvalues correspond to points
of the spectrum of $\ZZ(\g)$.

The answer comes from the description of the center of the
completed universal enveloping algebra of $\ghat$ at the critical
level. In \cite{FF:gd,F:wak} it is shown that the spectrum of this
center (more precisely, the center of the corresponding vertex
algebra) is canonically identified with the space of the so-called
$^L G$-{\em opers}, where $^L G$ is the {\em Langlands dual Lie
group} to $\g$ (of adjoint type), on the formal disc. This result
leads us to the following description of the spectrum of the
algebra $\ZZ(\g)$ of higher Gaudin Hamiltonians: it is the space
of $^L G$-opers on $\pone$ with regular singularities at the
points $z_1,\ldots,z_N$ and $\infty$.

We obtain this description from some basic facts about the spaces
of coinvariants from \cite{FB}. Recall that the space of
coinvariants $H_V(X;(x_i);(M_i))$ is defined in \cite{FB} for any
(quasi-conformal) vertex algebra $V$, a smooth projective curve
$X$, a collection $x_1,\ldots,x_n$ of distinct points of $X$ and a
collection of $V$-modules $M_1,\ldots,M_n$ attached to those
points. Suppose that $V$ is a commutative vertex algebra, and so
in particular it is a commutative algebra. Suppose that the
spectrum of $V$ is the space $S(D)$ of certain geometric objects,
such as $^L G$-opers, on the disc $D = \on{Spec} \C[[t]]$. Then a
$V$-module is the same as a smooth module over the complete
topological algebra $U(V)$ of functions on $S(D^\times)$, which is
the space of our objects (such as $^L G$-opers) on the punctured
disc $D^\times = \on{Spec} \C((t))$. Suppose in addition that each
$V$-module $M_i$ is the space of functions on a subspace $S_i$ of
$S(D^\times)$ (with its natural $\on{Fun} S(D^\times)$-module
structure). Then the space of coinvariants $H_V(X;(x_i);(M_i))$ is
naturally a commutative algebra, and its spectrum is the space of
our objects (such as $^L G$-opers) on $X \bs \{ x_1,\ldots,x_n
\}$ whose restriction to the punctured disc $D_{x_i}^\times$
around $x_i$ belongs to $S_i \subset S(D_{x_i}^\times),
i=1,\ldots,n$.

For example, if $\g=\sw_2$, then $^L G=PGL_2$, and $PGL_2$-opers are
the same as second order differential operators $\pa_t^2 - q(t)$
acting from sections of the line bundle $\Omega^{-1/2}$ to sections of
$\Omega^{3/2}$. A $PGL_2$-oper on $\pone$ with regular singularities
at $z_1,\ldots,z_N$ and $\infty$ may be written as the Fuchsian
differential operator of second order with regular singularities at
$z_1,\ldots,z_N$,
$$
\pa_t^2 - \sum_{i=1}^N \frac{c_i}{(t-z_i)^2} - \sum_{i=1}^N
\frac{\mu_i}{t-z_i},
$$
satisfying the condition $\sum_{i=1}^N \mu_i = 0$ that insures that it
also has regular singularity at $\infty$. Defining such an operator is
the same as giving a collection of numbers $c_i,\mu_i, i=1,\ldots,N$,
such that $\sum_{i=1}^N \mu_i = 0$. The set
$$
\left\{ c_i,\mu_i, i=1,\ldots,N \, \left| \sum_{i=1}^N \mu_i = 0
\right. \right\},
$$
is then the spectrum of the Gaudin algebra $\ZZ(\g)$, which in this
case is the polynomial algebra generated by the Casimir operators $C_i
= \frac{1}{2} \sum_a J_a^{(i)} J^{a(i)}, i=1,\ldots,N$, and the Gaudin
Hamiltonians $\Xi_i$, subject to the relation $\sum_{i=1}^N \Xi_i =
0$. In other words, the numbers $c_i$ record the eigenvalues of the
$C_i$'s, while the numbers $\mu_i$ record the eigenvalues of the
$\Xi_i$'s.

For a general simple Lie algebra $\g$, the Gaudin algebra ${\mc
Z}_N(\g)$ has many more generators, and its spectrum does not have
such a nice system of coordinates as the $c_i$'s and the $\mu_i$'s in
the above example. Therefore the description of the spectrum as a
space of $^L G$-opers is very useful. In particular, we obtain that
common eigenvalues of the higher Gaudin Hamiltonians are encoded by
$^L G$-opers on $\pone$, with regular singularities at prescribed
points. These $^L G$-opers appear as generalizations of the above
second order Fuchsian operators.

Next, we ask which points in the spectrum of $\ZZ(\g)$ might occur as
the common eigenvalues on particular tensor products
$\bigotimes_{i=1}^N M_i$. We answer this question first in the case
when each $M_i$ admits a central character: namely, it turns out that
the central character of $M_i$ fixes the residue of the $^L G$-oper
at the point $z_i$. We then show that if all $\g$-modules $M_i$ are
finite-dimensional and irreducible, then the $^L G$-opers encoding
possible eigenvalues of the higher Gaudin Hamiltonians in
$\bigotimes_{i=1}^N M_i$ necessarily have {\em trivial monodromy
representation}.

We conjecture that there is a bijection between the eigenvalues of the
Gaudin Hamiltonians on $\bigotimes_{i=1}^N M_i$, where the $M_i$'s are
irreducible finite-dimensional $\g$-modules, and $^L G$-opers on
$\pone$ with prescribed singularities at $z_1,\ldots,z_N,\infty$ and
trivial monodromy.

Thus, we obtain a correspondence between two seemingly unrelated
objects: the eigenvalues of the generalized Gaudin Hamiltonians and
the $^L G$-opers on $\pone$. The connection between the eigenvalues of
the Gaudin operators and differential operators of some sort has been
observed previously, but it was not until \cite{FFR,F:icmp} that this
phenomenon was explained conceptually.

We present a more geometric description the $^L G$-opers without
monodromy (which occur as the eigenvalues of the Gaudin
Hamiltonians) as isomorphism classes of holomorphic maps from
$\pone$ to $^L G/{}^L B$, the flag manifold of $^L G$, satisfying
a certain transversality condition. For example, if $^L G =
PGL_2$, they may be described as holomorphic maps $\pone \to
\pone$ whose derivative vanishes to prescribed orders at the
marked points $z_1,\ldots,z_N$ and $\infty$, and does not vanish
anywhere else (these orders correspond to the highest weights of
the finite-dimensional representations inserted at those points).

If the $^L G$--oper is non-degenerate (in the sense explained in
\secref{from}), then we can associate to it an eigenvector of the
Gaudin Hamiltonians called a {\em Bethe vector}. The procedure to
construct eigenvectors of the Gaudin Hamiltonians that produces these
vectors is known as the {\em Bethe Ansatz}. In \cite{FFR} we explained
that this procedure can also be understood in the framework of
coinvariants of $\ghat$-modules of critical level. We need to use a
particular class of $\ghat$-modules, called the {\em Wakimoto
modules}.

Let us recall that the Wakimoto modules of critical level are
naturally parameterized by objects closely related to opers, which
we call {\em Miura opers}. They may also be described more
explicitly as certain connections on a particular $^L H$-bundle
$\Omega^{-\rho}$ on the punctured disc, where $^L H$ is the Cartan
subgroup of $^L G$. The center acts on the Wakimoto module
corresponding to a Cartan connection by the Miura transformation
of this connection (see \cite{F:wak}). The idea of \cite{FFR} was
to use the spaces of coinvariants of the tensor product of the
Wakimoto modules to construct eigenvectors of the generalized
Gaudin Hamiltonians. We found in \cite{FFR} that the eigenvalues
of the Gaudin Hamiltonians on these vectors are encoded by the $^L
G$-opers which are obtained by applying the Miura transformation
to certain Cartan connections on $\pone$.

More precisely, the Bethe vector depends on an $m$-tuple of
complex numbers $w_{j}$, where $j=1,\ldots,m$, with an extra datum
attached to each of them, $i_j \in I$, where $I$ is the set of
nodes of the Dynkin diagram of $\g$ (or equivalently, the set of
simple roots of $\g$). These numbers have to be distinct from the
$z_i$'s and satisfy the following system of {\em Bethe Ansatz
equations}:
\begin{equation} \label{BAE}
\sum_{i=1}^N \frac{\langle \la_i,\chal_{i_j} \rangle}{w_j-z_i} -
\sum_{s \neq j} \frac{\langle \al_{i_s},\chal_{i_j}
\rangle}{w_j-w_s} = 0, \quad j=1,\ldots,m,
\end{equation}
where $\la_i$ denotes the highest weight of the finite-dimensional
$\g$-module $M_i = V_{\la_i}, i=1,\ldots,N$.

We can compute explicitly the $^L G$-oper encoding the
eigenvalues of the generalized Gaudin Hamiltonians on this vector.
As shown in \cite{FFR}, this $^L G$-oper is obtained by applying
the Miura transformation of the connection
\begin{equation} \label{cc}
\pa_t + \sum_{i=1}^N \frac{\la_i}{t-z_i} - \sum_{j=1}^m
\frac{\al_{i_j}}{t-w_j}
\end{equation}
on the $^L H$-bundle $\Omega^{-\rho}$ on $\pone$. This $^L
G$-oper automatically has trivial monodromy.

The Bethe vector corresponding to a solution of the system
\eqref{BAE} is a highest weight vector in $\bigotimes_{i=1}^N
V_{\la_i}$ of weight $$\mu = \sum_{i=1}^N \la_i - \sum_{j=1}^m
\al_{i_j},$$ so it can only be non-zero if $\mu$ is a dominant
integral weight of $\g$. But it is still interesting to describe
the set of {\em all} solutions of the Bethe Ansatz equations
\eqref{BAE}, even for non-dominant weights $\mu$.

While the eigenvalues of the Gaudin Hamiltonians are parameterized
by $^L G$-opers, it turns out that the solutions of the Bethe
Ansatz equations are parameterized by the (non-degenerate) Miura
$^L G$-opers. As mentioned above, those may in turn be related to
very simple objects, namely, connections on an $^L H$-bundle
$\Omega^{-\rho}$ of the kind given above in formula \eqref{cc}.

A $^L G$-oper on a curve $X$ is by definition a triple
$(\F,\nabla,\F_{^L B})$, where $\F$ is a $^L G$-bundle on $X$,
$\nabla$ is a connection on $\F$ and $\F_{^L B}$ is a reduction of
$\F$ to a Borel subgroup $^L B$ of $^L G$, which satisfies a certain
transversality condition with $\nabla$. A Miura $^L G$-oper is by
definition a quadruple $(\F,\nabla,\F_{^L B},\F'_{^L B})$ where
$\F'_{^L B}$ is another $^L B$-reduction of $\F$, which is preserved
by $\nabla$. The space of Miura opers on a curve $X$ whose underlying
oper has regular singularities and trivial monodromy representation
(so that $\F$ is isomorphic to the trivial bundle) is isomorphic to
the {\em flag manifold} $^L G/{}^L B$ of $^L G$. Indeed, in order to
define the $^L B$-reduction $\F'_{^L B}$ of such $\F$ everywhere, it
is sufficient to define it at one point $x \in X$ and then use the
connection to ``spread'' it around. But choosing a $^L B$-reduction
at one point means choosing an element of the twist of $^L G/{}^L B$
by $\F_x$, and so we see that the space of all reductions is
isomorphic to the flag manifold of $^L G$.

The relative position of the two reductions $\F_{^L B}$ and $\F'_{^L
B}$ at each point of $X$ is measured by an element $w$ of the Weyl
group $W$ of $G$. The two reductions are in generic relative position
(corresponding to $w=1$) almost everywhere on $X$. The mildest
possible non-generic relative positions correspond to the simple
reflections $s_i$ from $W$. We call a Miura oper on $\pone$ with
marked points $z_1,\ldots,z_N$ {\em non-degenerate} if the two
reductions $\F_B$ and $\F'_B$ are in generic position at
$z_1,\ldots,z_N$, and elsewhere on $\pone$ their relative position is
either generic or corresponds to a simple reflection. We denote the
points where the relative position is not generic by $w_j,
j=1,\ldots,m$; each point $w_j$ comes together with a simple
reflection $s_{i_j}$, or equivalently a simple root $\al_{i_j}$
attached to it.

It is then easy to see that this collection satisfies the
equations \eqref{BAE}, and conversely any solution of \eqref{BAE}
corresponds to a non-degenerate Miura oper (or to an $^L
H$-connection \eqref{cc}). Thus, we obtain that there is a
bijection between the set of solutions of \eqref{BAE} (for all
possible collections $\{ i_1,\ldots,i_m \}$) and the set of
non-degenerate Miura $^L G$-opers such that the underlying $^L
G$-opers have prescribed residues at the points
$z_1,\ldots,z_N,\infty$ and trivial monodromy.

Now let us fix $\la_1,\ldots,\la_N$ and $\mu$. Then every $^L G$-oper
$\tau$ on $\pone$ with regular singularities at $z_1,\ldots,z_N$ and
$\infty$ and with prescribed residues corresponding to
$\la_1,\ldots,\la_N$ and $\mu$ and trivial monodromy admits a
horizontal $^L B$-reduction $\F'_{^L B}$ satisfying the conditions of
a non-degenerate Miura oper. Since these are open conditions, we find
that for such $z_1,\ldots,z_N$ the non-degenerate Miura oper
structures on a particular $^L G$-oper $\tau$ on $\pone$ form an open
dense subset in the set of all Miura oper structures on $\tau$. But
the set of all Miura structures on a given $^L G$-oper $\tau$ is
isomorphic to the flag manifold $^L G/{}^L B$. Therefore we find that
the set of non-degenerate Miura oper structures on $\tau$ is an open
dense subset of $^L G/{}^L B$!

Recall that the set of all solutions of the Bethe Ansatz equations
\eqref{BAE} is the union of the sets of non-degenerate Miura oper
structures on all $^L G$-opers with trivial monodromy. Hence it is
naturally a disjoint union of subsets, parameterized by these $^L
G$-opers. We have now identified each of these sets with an open and
dense subset of the flag manifold $^L G/{}^L B$.

\medskip
\noindent
Let us summarize our results:

\begin{itemize}

\item the eigenvalues of the Hamiltonians of the Gaudin model
associated to a simple Lie algebra $\g$ on the tensor product of
finite-dimensional representations are encoded by $^L G$-opers on
$\pone$, where $^L G$ is the Langlands dual group of $G$, which have
regular singularities at the marked points $z_1,\ldots,z_N,\infty$ and
trivial monodromy;

\item if such an oper $\tau$ is non-degenerate, then we can associate
  to it a solution of the Bethe Ansatz equations \eqref{BAE},
which gives rise to the Bethe eigenvector of dominant integral weight
whose eigenvalues are encoded by $\tau$;

\item there is a one-to-one correspondence between the set of all
solutions of the Bethe Ansatz equations \eqref{BAE} and the set of
non-degenerate Miura opers corresponding to a fixed $^L G$-oper;

\item the set of non-degenerate Miura opers corresponding to the same
underlying $^L G$-oper is an open dense subset of the flag manifold
$^L G/{}^L B$ of the Langlands dual group, and therefore the set of
all solutions of the Bethe Ansatz equations \eqref{BAE} is the union
of certain open dense subsets of the flag manifold of the Langlands
dual group, one for each $^L G$-oper.

\end{itemize}

Finally, to a degenerate $^L G$--oper we can also attach, at least in
some cases, an eigenvector of the Gaudin hamiltonians by generalizing
the Bethe Ansatz procedure, as explained in \secref{last}.

The paper is organized as follows. In \secref{opers1} we introduce
opers and discuss their basic properties. We define opers with
regular singularities and their residues. In \secref{third},
following \cite{FFR}, we define the Gaudin algebra using the
coinvariants of the affine Kac-Moody algebra $\ghat$ of critical
level. We recall the results of \cite{FF:gd,F:wak} on the
isomorphism of the center of the completed universal enveloping
algebra of $\ghat$ at the critical level and the algebra of
functions on the space of $^L G$-opers on the punctured disc.
Using these results and general facts about the spaces of
coinvariants from \cite{FB}, we describe the spectrum of the
Gaudin algebra. In \secref{miura opers} we introduce Miura opers,
Cartan connections and the Miura transformation and describe their
properties, following \cite{F:wak,F:opers}. We use these results
in the next section, \secref{wak and ba}, to describe the Bethe
Ansatz, a construction of eigenvectors of the Gaudin algebra. We
introduce the Wakimoto modules of critical level, following
\cite{FF:si,F:wak}. The Wakimoto modules are naturally
parameterized by the Cartan connections on the punctured disc
introduced in \secref{miura opers}. The action of the center on
the Wakimoto modules is given by the Miura transformation. We
construct the Bethe vectors, following \cite{FFR}, using the
coinvariants of the Wakimoto modules. We show that the Bethe
Ansatz equations which ensure that this vector is an eigenvector
of the Gaudin Hamiltonians coincide with the requirement that the
Miura transformation of the Cartan connection on $\pone$ encoding
the Wakimoto modules has no singularities at the points
$w_1,\ldots,w_m$. Finally, in \secref{findim mod} we consider the
Gaudin model in the case when all modules $M_i$ finite-dimensional
modules. We describe the precise connection between the spectrum
of the Gaudin algebra on the tensor product of finite-dimensional
modules and the set of $^L G$-opers with prescribed singularities
at $z_1,\ldots,z_N,\infty$ and trivial monodromy.

\section{Opers} \label{opers1}

\subsection{Definition of opers} \label{opers}

Let $G$ be a simple algebraic group of adjoint type, $B$ a Borel
subgroup and $N = [B,B]$ its unipotent radical, with the corresponding
Lie algebras $\n \subset \bb\subset \g$. The quotient $H = B/N$ is a
torus. Choose a splitting $H \to B$ of the homomorphism $B \to H$ and
the corresponding splitting $\h \to \bb$ at the level of Lie
algebras. Then we will have a Cartan decomposition $\g = \n_- \oplus
\h \oplus \n$. We will choose generators $\{ e_i \}, i=1,\ldots,\ell$,
of $\n$ and generators $\{ f_i \}, i=1,\ldots,\ell$ of $\n_-$
corresponding to simple roots, and denote by $\crho \in \h$ the sum of
the fundamental coweights of $\g$. Then we will have the following
relations: $[\crho,e_i] = 1, [\crho,f_i] = -1$.

A $G$-oper on a smooth curve $X$ (or a disc $D \simeq \on{Spec}
\C[[t]]$ or a punctured disc $D^\times = \on{Spec} \C((t))$) is by
definition a triple $(\F,\nabla,\F_B)$, where $\F$ is a principal
$G$-bundle $\F$ on $X$, $\nabla$ is a connection on $\F$ and $\F_B$
is a $B$-reduction of $\F$ such that locally on $X$ (with respect to
a local coordinate $t$ and a local trivialization of $\F_B$) the
connection has the form
\begin{equation} \label{form of nabla}
\nabla = \pa_t + \sum_{i=1}^\ell \psi_i(t) f_i + {\mb v}(t),
\end{equation}
where each $\psi_i(t)$ is a nowhere vanishing function, and ${\mb
v}(t)$ is a $\bb$-valued function. The space of $G$-opers on $X$ is
denoted by $\on{Op}_G(X)$.

This definition is due to A. Beilinson and V. Drinfeld \cite{BD} (in
the case when $X$ is the punctured disc opers were first introduced in
\cite{DS}).

In particular, if $U = \on{Spec} R$ is an affine curve with the
ring of functions $R$ and $t$ is a global coordinate on $U$ (for
example, if $U = \on{Spec} \C[[t]]$), then $\on{Op}_G(U)$ is
isomorphic to the quotient of the space of operators of the form
\begin{equation} \label{another form of nabla}
\nabla = \pa_t + \sum_{i=1}^\ell f_i + {\mb v}(t), \qquad {\mb v}(t)
\in \bb(R),
\end{equation}
by the action of the group $N(R)$ (we use the action of $H(R)$ to make
all functions $\psi_i(t)$ equal to $1$). Recall that the gauge
transformation of an operator $\pa_t+A(t)$, where $A(t) \in \g(R)$ by
$g(t) \in G(R)$ is given by the formula
$$
g \cdot (\pa_t + A(t)) = \pa_t + g A(t) g^{-1} - \pa_t g \cdot g^{-1}.
$$

The operator $\on{ad} \crho$ defines the principal gradation on $\bb$,
with respect to which we have a direct sum decomposition $\bb =
\bigoplus_{i\geq 0} \bb_i$. Set
$$
p_{-1} = \sum_{i=1}^\ell f_i.
$$
Let $p_1$ be the unique element of degree 1 in $\n$, such that $\{
p_{-1},2\rv,p_1 \}$ is an $\sw_2$-triple. Let $V_{\can} = \oplus_{i
\in E} V_{\can,i}$ be the space of $\on{ad} p_1$-invariants in
$\n$. Then $p_1$ spans $V_{\on{can},1}$. Choose a linear
generator $p_j$ of $V_{\can,d_j}$ (if the multiplicity of $d_j$ is
greater than one, which happens only in the case $\g=D^{(1)}_{2n},
d_j=2n$, then we choose linearly independent vectors in
$V_{\on{can},d_j}$). The following result is due to Drinfeld and
Sokolov \cite{DS} (the proof is reproduced in Lemma 2.1 of
\cite{F:opers}).

\begin{lem} \label{free}
The gauge action of $N(R)$ on $\wt{\on{Op}}_G(\on{Spec} R)$ is free,
and each gauge equivalence class contains a unique operator of the
form $\nabla = \pa_t + p_{-1} + {\mathbf v}(t)$, where ${\mathbf v}(t)
\in V_{\can}(R)$, so that we can write
\begin{equation} \label{coeff fun}
{\mathbf v}(t) = \sum_{j=1}^\ell v_j(t) \cdot p_j.
\end{equation}
\end{lem}

\subsection{$PGL_2$-opers}

For $\g=\sw_2, G=PGL_2$ we obtain an identification of the space of
$PGL_2$-opers with the space of operators of the form
$$
\pa_t + \left( \begin{array}{ccccc}
0 & v(t) \\
1 & 0
\end{array} \right).
$$
If we make a change of variables $t=\varphi(s)$, then the
corresponding connection operator will become
$$
\pa_s + \left( \begin{array}{ccccc}
0 & \varphi'(s) v(\varphi(s)) \\
\varphi'(s) & 0
\end{array} \right).
$$
Applying the $B$-valued gauge transformation with
$$
\left( \begin{array}{ccccc}
1 & \frac{1}{2} \frac{\varphi''(s)}{\varphi'(s)} \\
0 & 1
\end{array} \right) \left( \begin{array}{ccccc}
(\varphi'(s))^{1/2} & 0 \\
0 & (\varphi'(s))^{1/2}
\end{array} \right),
$$
we obtain the operator
$$
\pa_s + \left( \begin{array}{ccccc}
0 & v(\varphi(s)) \varphi'(s)^2 - \frac{1}{2} \{ \varphi,s \} \\
1 & 0
\end{array} \right),
$$
where
$$
\{ \varphi,s \} = \frac{\varphi'''}{\varphi'} - \frac{3}{2}
\left( \frac{\varphi''}{\varphi'} \right)^2
$$
is the Schwarzian derivative. Thus, under the change of variables
$t = \varphi(s)$ we have
$$
v(t) \mapsto v(\varphi(s)) \varphi'(s)^2 - \frac{1}{2} \{ \varphi,s
\}.
$$
this coincides with the transformation properties of the second order
differential operators $\pa_t^2 - v(t)$ acting from sections of
$\Omega^{-1/2}$ to sections of $\Omega^{3/2}$, where $\Omega$ is the
canonical line bundle on $X$. Such operators are known as {\em
projective connections} on $X$ (see, e.g., \cite{FB}, Sect. 9.2), and
so $PGL_2$-opers are the same as projective connections.

For a general $\g$, the first coefficient function $v_1(t)$ in
\eqref{coeff fun} transforms as a projective connection, and the
coefficient $v_i(t)$ with $i>1$ transforms as a $(d_i+1)$-differential
on $X$. Thus, we obtain an isomorphism
\begin{equation} \label{repr}
\on{Op}_G(X) \simeq {\mc P}roj(X) \times \bigoplus_{i=2}^\el
\Gamma(X,\Omega^{(d_i+1)}).
\end{equation}

\subsection{Opers with regular singularities} \label{reg sing}

Let $x$ be a point of a smooth curve $X$ and $D_x = \on{Spec} \OO_x,
D^\times_x = \on{Spec} \K_x$, where $\OO_x$ is the completion of the
local ring of $x$ and $\K_x$ is the field of fractions of $\OO_x$.
Choose a formal coordinate $t$ at $x$, so that $\OO_x \simeq \C[[t]]$
and $\K_x = \C((t))$. Recall that the space $\on{Op}_G(D_x)$ (resp.,
$\on{Op}_G(D_x^\times)$) of $G$-opers on $D_x$ (resp., $D_x^\times$)
is the quotient of the space of operators of the form \eqref{form of
nabla} where $\psi_i(t)$ and ${\mb v}(t)$ take values in $\OO_x$
(resp., in $\K_x$) by the action of $B(\OO_x)$ (resp., $B(\K_x)$).

A $G$-oper on $D_x$ with regular singularity at $x$ is by definition
(see \cite{BD}, Sect. 3.8.8) a $B(\OO_x)$-conjugacy class of
operators of the form
\begin{equation} \label{oper with RS1}
\nabla = \pa_t + t^{-1} \left( \sum_{i=1}^\ell \psi_i(t) f_i + {\mb
v}(t) \right),
\end{equation}
where $\psi_i(t) \in \OO_x, \psi_i(0) \neq 0$, and ${\mb v}(t) \in
\bb(\OO_x)$. Equivalently, it is an $N(\OO_x)$-equivalence class of
operators
\begin{equation} \label{oper with RS}
\nabla = \pa_t + \frac{1}{t} \left( p_{-1} + {\mb v}(t) \right),
\qquad {\mb v}(t) \in \bb(\OO_x).
\end{equation}
Denote by $\on{Op}_G^{\on{RS}}(D_x)$ the space of opers on $D_x$ with
regular singularity. It is easy to see that the natural map
$\on{Op}_G^{\on{RS}}(D_x) \to \on{Op}_G(D_x^\times)$ is
injective. Therefore an oper with regular singularity may be viewed as
an oper on the punctured disc. But to an oper with regular singularity
one can unambiguously attach a point in
$$\g/G := \on{Spec} \C[\g]^G \simeq \C[\h]^W =: \h/W,$$ its residue,
which in our case is equal to $p_{-1} + {\mb v}(0)$.

Given $\cla \in \h$, we denote by $\on{Op}_G^{\on{RS}}(D_x)_{\cla}$
the subvariety of $\on{Op}_G^{\on{RS}}(D_x)$ which consists of those
opers that have residue $\varpi(- \cla - \crho) \in \h/W$, where
$\varpi$ is the projection $\h \to \h/W$.

In particular, the residue of a regular oper $\pa_t + p_{-1} + {\mb
v}(t)$, where ${\mb v}(t) \in \bb(\OO_x)$, is equal to
$\varpi(-\crho)$ (see \cite{BD}). Indeed, a regular oper may be
brought to the form \eqref{oper with RS} by using the gauge
transformation with $\crho(t) \in B(\K_x)$, after which it takes the
form
$$
\pa_t + \frac{1}{t} \left( p_{-1} - \crho + t \cdot
\crho(t) ({\mb v}(t)) \crho(t)^{-1} \right).
$$
If ${\mb v}(t)$ is regular, then so is $\crho(t) ({\mb v}(t))
\crho(t)^{-1}$. Therefore the residue of this oper in $\h/W$ is
equal to $\varpi(-\crho)$, and so $\on{Op}_G(D_x) =
\on{Op}_G^{\on{RS}}(D_x)_{0}$.

For $G=PGL_2$ we identify $\h$ with $\C$ and so $\cla$ with a complex
number. Then one finds that $\on{Op}_{PGL_2}^{\on{RS}}(D_x)_{\cla}$ is
the space of second order operators of the form
\begin{equation} \label{sl2 reg sing oper}
\pa_t^2 - \frac{\cla(\cla+2)/4}{t^2} - \sum_{n\geq -1} v_n t^n.
\end{equation}

Now suppose that $\cla$ is a dominant integral coweight of
$\g$. Following Drinfeld, introduce the variety
$\on{Op}_G(D_x)_{\cla}$ as the quotient of the space of operators of
the form
\begin{equation} \label{psi la}
\nabla = \pa_t + \sum_{i=1}^\ell \psi_i(t) f_i + {\mb v}(t),
\end{equation}
where $$\psi_i(t) = t^{\langle \al_i,\cla \rangle}(\kappa_i +
t(\ldots)) \in \OO_x, \qquad \kappa_i \neq 0$$ and ${\mb v}(t) \in
\bb(\OO_x)$, by the gauge action of $B(\OO_x)$. Equivalently,
$\on{Op}_G(D_x)_{\cla}$ is the quotient of the space of operators of
the form
\begin{equation} \label{psi la1}
\nabla = \pa_t + \sum_{i=1}^\ell t^{\langle \al_i,\cla \rangle} f_i +
{\mb v}(t),
\end{equation}
where ${\mb v}(t) \in \bb(\OO_x)$, by the gauge action of
$N(\OO_x)$. Considering the $N(\K_x)$-class of such an operator, we
obtain an oper on $D_x^\times$. Thus, we have a map
$\on{Op}_G(D_x)_{\cla} \to \on{Op}_G(D_x^\times)$.

\begin{lem}[\cite{F:opers}, Lemma 2.4] \label{no mon}
The map $\on{Op}_G(D_x)_{\cla} \to \on{Op}_G(D_x^\times)$ is injective
and its image is contained in the subvariety
$\on{Op}_G^{\on{RS}}(D_x)_{\cla}$. Moreover, the points of
$\on{Op}_G(D_x)_{\cla}$ are precisely those $G$-opers with regular
singularity and residue $\cla$ which have no monodromy around $x$.
\end{lem}

\noindent
The space $\on{Op}_{PGL_2}(D_x)_{\cla}$ is the
subspace of codimension
one in $\on{Op}_{PGL_2}(D_x)^{\on{RS}}_{\cla}$. In terms of the
coefficients $v_n, n \geq -1$, appearing in formula \eqref{sl2 reg
sing oper} the corresponding equation has the form $P_\la(v_n) = 0$,
where $P_\la$ is a polynomial of degree $\la+1$, where we set $\deg
v_n = n+2$. For example, $P_0 = v_{-1}, P_2 = 2 v_{-1}^2 - v_0$,
etc. In general, the subspace $\on{Op}_G(D_x)_{\cla} \subset
\on{Op}_G^{\on{RS}}(D_x)_{\cla}$ is defined by $\on{dim} N$ polynomial
equations, where $N$ is the unipotent subgroup of $G$.

\section{The Gaudin model} \label{higher
gaudin} \label{third}

Let $\g$ be a simple Lie algebra. The {\em Langlands dual Lie algebra}
$^L \g$ is by definition the Lie algebra whose Cartan matrix is the
transpose of that of $\g$. We will identify the set of roots of $\g$
with the set of coroots of $^L \g$ and the set of weights of $\g$ with
the set of coweights of $^L \g$. The results on opers from the
previous sections will be applied here to the Lie algebra $^L
\g$. Thus, in particular, $^L G$ will denote the adjoint group of $^L
\g$.

\subsection{The definition of the Gaudin model}

Here we recall the definition of the Gaudin model and the realization
of the Gaudin Hamiltonians in terms of the spaces of conformal blocks
for affine Kac-Moody algebras of critical level. We follow closely the
paper \cite{FFR}.

Choose a non-degenerate invariant inner product $\ka_0$ on $\g$. Let
$\{ J_a \}, a=1,\ldots,d$, be a basis of $\g$ and $\{ J^a \}$ the dual
basis with respect to $\ka_0$. Denote by $\Delta$ the quadratic
Casimir operator from the center of $U(\g)$:
$$\Delta = \frac{1}{2} \sum_{a=1}^d J_a J^a.$$

The {\em Gaudin Hamiltonians} are the elements
\begin{equation} \label{hi}
\Xi_i = \sum_{j\neq i} \sum_{a=1}^d \frac{J_a^{(i)}
J^{a(j)}}{z_i-z_j}, \qquad i=1,\ldots,N,
\end{equation}
of the algebra $U(\g)^{\otimes N}$. Note that they commute with the
diagonal action of $\g$ on $U(\g)^{\otimes N}$ and that $$\sum_{i=1}^N
\Xi_i = 0.$$

\subsection{Gaudin model and coinvariants} \label{coinvariants}
Let $\ghat_{\ka_c}$ be the affine Kac-Moody algebra corresponding
to $\g$. It is the extension of the Lie algebra $\g \otimes \C((t))$
by the one-dimensional center $\C K$. The commutation relations in
$\ghat_{\ka_c}$ read
\begin{equation}
\label{comm rel}
[A \otimes f(t),B \otimes g(t)] = [A,B] \otimes fg - \ka_c(A,B)
\on{Res}_{t=0} f dg \cdot K,
\end{equation}
where $\ka_c$ is the {\em critical} invariant inner product on $\g$
defined by the formula $$\ka_c(A,B) = - \frac{1}{2} \on{Tr}_\g \on{ad}
A \on{ad} B.$$ Note that $\ka_c = - h^\vee \ka_0$, where $\ka_0$ is
the inner product normalized as in \cite{Kac} and $h^\vee$ is the dual
Coxeter number.

Denote by $\G_+$ the Lie subalgebra $\g \otimes \C[[t]] \oplus \C K$
of $\ghat_{\ka_c}$. Let $M$ be a $\g$-module. We extend the action of
$\g$ on $M$ to $\g \otimes \C[[t]]$ by using the evaluation at zero
homomorphism $\g \otimes \C[[t]] \to \g$ and to $\G_+$ by making $K$
act as the identity. Denote by ${\mathbb M}$ the corresponding induced
$\ghat_{\ka_c}$-module
$${\mathbb M} = U(\G_{\ka_c}) \underset{U(\G_+)}\otimes M.$$ By
construction, $K$ acts as the identity on this module. We call such
modules the $\ghat_{\ka_c}$-modules of {\em critical level}. For
example, for $\la \in \h^*$ let $\C_\la$ be the one-dimensional
$\bb$-module on which $\h$ acts by the character $\la: \h \to \C$ and
$\n$ acts by $0$. Let $M_{\la}$ be the Verma module over $\g$ of
highest weight $\la$,
$$
M_\la = U(\g) \underset\otimes{U(\bb)} \C_\la.
$$
The corresponding induced module $\M_{\la}$ is the Verma module over
$\ghat_{\ka_c}$ with highest weight $\la$.

For a dominant integral weight $\la \in \h^*$ denote by $V_\la$ the
irreducible finite-dimensional $\g$-module of highest weight
$\la$. The corresponding induced module $\V_\la$ is called the Weyl
module over $\ghat_{\ka_c}$ with highest weight $\la$.

Consider the projective line $\pone$ with a global coordinate $t$ and
$N$ distinct finite points $z_1,\ldots,z_N \in \pone$. In the
neighborhood of each point $z_i$ we have the local coordinate $t-z_i$
and in the neighborhood of the point $\infty$ we have the local
coordinate $t^{-1}$. Set $\widetilde{\g}(z_i) = \g \otimes
\C((t-z_i))$ and $\wt\g(\infty) = \g \otimes \C((t^{-1}))$. Let $\G_N$
be the extension of the Lie algebra $\bigoplus_{i=1}^N
\widetilde{\g}(z_i) \oplus \wt\g(\infty)$ by a one-dimensional
center $\C K$ whose restriction to each summand $\widetilde{\g}(z_i)$
or $\wt\g(\infty)$ coincides with the above central extension.

Suppose we are given a collection $M_1,\ldots,M_N$ and $M_\infty$ of
$\g$-modules. Then the Lie algebra $\G_N$ naturally acts on the
tensor product $\bigotimes_{i=1}^N {\mathbb M}_i \otimes {\mathbb
M}_\infty$ (in particular, $K$ acts as the identity).

Let $\g_{\zn}=\g_{z_1,\ldots,z_N}$ be the Lie algebra of $\g$-valued
regular functions on $\pone\backslash\{ z_1, \ldots,$ $z_N,\infty \}$
(i.e., rational functions on $\pone$, which may have poles only at the
points $z_1,\ldots,z_N$ and $\infty$). Clearly, such a function can be
expanded into a Laurent power series in the corresponding local
coordinates at each point $z_i$ and at $\infty$. Thus, we obtain an
embedding $$\g_{\zn} \hookrightarrow \bigoplus_{i=1}^N
\widetilde{\g}(z_i) \oplus \wt\g(\infty).$$ It follows from the
residue theorem and formula \eqref{comm rel} that the restriction of
the central extension to the image of this embedding is trivial. Hence
this embedding lifts to the embedding $\g_{\zn} \to \G_N$.

Denote by $H(M_1,\ldots,M_N,M_\infty)$ the space of coinvariants of
$\bigotimes_{i=1}^N {\mathbb M}_i \otimes {\mathbb M}_\infty$ with
respect to the action of the Lie algebra $\g_{\zn}$. By construction,
we have a canonical embedding of a $\g$-module $M$ into the induced
$\ghat_{\ka_c}$-module $\M$:
$$x \in M \arr 1 \otimes x \in \M,$$ which commutes with the action of
$\g$ on both spaces (where $\g$ is embedded into $\ghat_{\ka_c}$ as
the constant subalgebra). Thus we have an embedding
$$
\bigotimes_{i=1}^N M_i \otimes M_\infty \hookrightarrow
\bigotimes_{i=1}^N \M_i \otimes \M_\infty.
$$
The following result is proved in the same way as Lemma 1 in
\cite{FFR}.

\begin{lem} \label{iso}
The composition of this embedding and the projection
$$\bigotimes_{i=1}^N \M_i \otimes \M_\infty \twoheadrightarrow
H(M_1,\ldots,M_N,M_\infty)$$ gives rise to an isomorphism
$$
H(M_1,\ldots,M_N,M_\infty) \simeq (\bigotimes_{i=1}^N M_i
\otimes M_\infty)/\g_{\on{diag}}.
$$
\end{lem}

Let $\V_0$ be the induced $\ghat_{\ka_c}$-module of critical level,
which corresponds to the one-dimensional trivial $\g$-module $V_0$;
it is called the {\em vacuum module}. Denote by $v_0$ the generating
vector of $\V_0$. We assign the vacuum module to a point $u \in \pone$
which is different from $z_1,\ldots,z_N,\infty$. Denote by
$H(M_1,\ldots,M_N,M_\infty,\C)$ the space of $\g_{\zn,u}$-invariant
functionals on $\bigotimes_{i=1}^N \M_i \otimes \M_\infty \otimes
\V_0$ with respect to the Lie algebra $\g_{\zn,u}$. \lemref{iso} tells
us that we have a canonical isomorphism
$$
H(M_1,\ldots,M_N,M_\infty,\C) \simeq H(M_1,\ldots,M_N,M_\infty).
$$

Now observe that by functoriality any endomorphism $X \in
\on{End}_{\ghat_{\ka_c}} \V_0$ gives rise to an endomorphism of the
space of coinvariants $H(M_1,\ldots,M_N,M_\infty,\C)$, and hence of
$H(M_1,\ldots,M_N,M_\infty)$. Thus, we obtain a homomorphism of
algebras
$$
\on{End}_{\ghat_{\ka_c}} \V_0 \to \on{End}_\C
H(M_1,\ldots,M_N,M_\infty).
$$

Let us compute this homomorphism explicitly.

\goodbreak

First of all, we identify the algebra $\on{End}_{\ghat_{\ka_c}} \V_0$
with the space
$$
\zz(\G) = \V_0^{\g[[t]]}
$$
of $\g[[t]]$-invariant vectors in $\V_0$. Indeed, a
$\g[[t]]$-invariant vector $v$ gives rise to an
endomorphism of $\V_0$ commuting with the action of $\ghat_{\ka_c}$
which sends the generating vector $v_0$ to $v$. Conversely,
any $\ghat_{\ka_c}$-endomorphism
of $V_0$ is uniquely determined by the image of $v_0$ which
necessarily belongs to $\zz(\G)$. Thus, we obtain an isomorphism
$\zz(\G) \simeq \on{End}_{\ghat_{\ka_c}}(\V_0)$ which gives $\zz(\G)$
an algebra
structure. The opposite algebra structure on $\zz(\G)$ coincides with
the algebra structure induced by the identification of $\V_0$ with the
algebra $U(\g \otimes t^{-1} \C[t^{-1}])$. But we will see in the next
section that the algebra $\zz(\G)$ is commutative and so the two
algebra structures on it coincide.

Now let $v \in \zz(\ghat) \subset \V_0$. For any $$x \in
(\bigotimes_{i=1}^N M_i \otimes M_\infty)/\g_{\on{diag}} \simeq
H(M_1,\ldots,M_N,M_\infty)$$ take a lifting $\wt{x}$ to
$\bigotimes_{i=1}^N M_i \otimes M_\infty$. By \lemref{iso}, the
projection of the vector $$\wt{x} \otimes v \in \bigotimes_{i=1}^N
\M_i \otimes \M_\infty \otimes \V_0$$ onto
$$H(M_1,\ldots,M_N,M_\infty,\C) \simeq (\bigotimes_{i=1}^N M_i
\otimes M_\infty)/\g_{\on{diag}}$$ is equal to the projection of a
vector of the form $(\Psi_u(v) \cdot \wt{x}) \otimes v_0$, where
$$\Psi_u(v) \cdot \wt{x} \in (\bigotimes_{i=1}^N M_i \otimes
M_\infty)/\g_{\on{diag}}.$$

For $A \in \g$ and $n \in \Z$, denote by $A_n$ the element $A \otimes
t^n \in \ghat_{\ka_c}$. Then $\V_0 \simeq U(\g \otimes t^{-1} \C[t^{-1}])
v_0$ has a basis of lexicographically ordered monomials of the form
$J^{a_1}_{n_1} \ldots J^{a_m}_{n_m} v_0$ with $n_i < 0$. Let us set
$$
{\mathbb J}^a_n(u) = - \sum_{i=1}^N \frac{J^{a(i)}}{(z_i-u)^n}.
$$
Define an anti-homomorphism
$$
\Phi_u: U(\g \otimes t^{-1} \C[t^{-1}]) \to U(\g)^{\otimes N} \otimes
\C[(u-z_i)^{-1}]_{i=1,\ldots,N}
$$
by the formula
\begin{equation} \label{Psi}
\Phi_u(J^{a_1}_{n_1} \ldots J^{a_m}_{n_m} v_0) = {\mathbb
J}^{a_m}_{n_m}(u) {\mathbb J}^{a_{m-1}}_{n_{m-1}}(u)\ldots {\mathbb
J}^{a_1}_{n_1}(u).
\end{equation}
According to the computation presented in the proof of Prop. 1
of \cite{FFR}, we have $$\Psi_u(v) \cdot \wt{x} = \Phi_u(v) \cdot
\wt{x}.$$

In general, $\Phi_u(v)$ does not commute with the diagonal action
of $\g$, and so $\Psi_v(u) \cdot \wt{x}$ depends on the choice of
the lifting $\wt{x}$. But if $v \in \zz(\G) \subset \V_0$, then
$\Phi_u(v)$ commutes with the diagonal action of $\g$ and hence
gives rise to a well-defined endomorphism $(\bigotimes_{i=1}^N M_i
\otimes M_\infty)/\g_{\on{diag}}$.

Thus, we restrict $\Phi_u$ to $\zz(\G)$. This gives us a homomorphism
of algebras
$$
\zz(\ghat) \to \left( U(\g)^{\otimes N} \right)^G \otimes
\C[(u-z_i)^{-1}]_{i=1,\ldots,N},
$$
which we also denote by $\Phi_u$.

For example, consider the Segal-Sugawara vector in $\V_0$:
\begin{equation} \label{sugawara}
S = \frac{1}{2} \sum_{a=1}^d J_{a,-1} J^a_{-1} v_0.
\end{equation}
One shows (see, e.g., \cite{FB}) that this vector belongs to
$\zz(\G)$. Consider the corresponding element $\Phi_u(S)$.

Denote by $\Delta$ the Casimir operator $\dfrac{1}{2} \sum_a J_a J^a$
from $U(\g)$.

\begin{prop}[\cite{FFR}, Prop. 1] \label{coincide}
We have
$$
\Phi_u(S) = \sum_{i=1}^N \frac{\Xi_i}{u-z_i} + \sum_{i=1}^N
\frac{\Delta^{(i)}}{(u-z_i)^2},
$$
where the $\Xi_i$'s are the Gaudin operators \eqref{hi}.
\end{prop}

We wish to study the algebra generated by the image of the map
$\Phi_u$.

\begin{prop}[\cite{FFR}, Prop. 2] \label{commu}
For any $Z_1, Z_2 \in \zz(\G)$ and any points $u_1, u_2$ $\in \pone
\bs \{ z_1,\ldots,z_N,\infty \}$ the linear operators
$\Psi_{Z_1}(u_1)$ and $\Psi_{Z_2}(u_2)$ commute.
\end{prop}

Let $\ZZ(\g)$ be the span in $U(\g)^{\otimes N}$ of the coefficients
in front of the monomials $\prod_{i=1}^N (u-z_i)^{n_i}$ of the series
$\Phi_u(v), v \in \zz(\ghat)$. Since $\Phi_u$ is an algebra
homomorphism, we find that $\ZZ(\g)$ is a subalgebra of $\left(
U(\g)^{\otimes N} \right)^G$, which is commutative by
\propref{commu}. We call it the {\em Gaudin algebra} associated to
$\g$ and the collection $z_1,\ldots,z_N$, and its elements the {\em
generalized Gaudin Hamiltonians}.

\subsection{The center of $\V_0$ and $^L G$-opers}

In order to describe the Gaudin algebra $\ZZ(\g)$ and its spectrum we
need to recall the description of $\zz(\G)$. According to
\cite{FF:gd,F:wak}, $\zz(\G)$ is identified with the algebra $\on{Fun}
\on{Op}_{^L G}(D)$ of (regular) functions on the space $\on{Op}_{^L
G}(D)$ of $^L G$-opers on the disc $D = \on{Spec} \C[[t]]$, where $^L
G$ is the {\em Langlands dual group} to $G$. Since we have assumed
that $G$ is simply-connected, $^L G$ may be defined as the adjoint
group of the Lie algebra $^L \g$ whose Cartan matrix is the transpose
of that of $\g$.

This isomorphism satisfies various properties, one of which we will
now recall. Let $\DerO = \C[[t]] \pa_t$ be the Lie algebra of
continuous derivations of the topological algebra $\OO = \C[[t]]$. The
action of its Lie subalgebra $\on{Der}_0 \OO = t \C[[t]] \pa_t$ on
$\OO$ exponentiates to an action of the group $\AutO$ of formal
changes of variables. Both $\DerO$ and $\AutO$ naturally act on $V_0$
in a compatible way, and these actions preserve $\zz(\G)$. They also
act on the space $\on{Op}_{^L G}(D)$.

Denote by $\on{Fun} \on{Op}_{^L G}(D)$ the algebra of regular
functions on $\on{Op}_{^L G}(D)$. In view of \lemref{free}, it is
isomorphic to the algebra of functions on the space of $\ell$-tuples
$(v_1(t),\ldots,v_\ell(t))$ of formal Taylor series, i.e., the space
$\C[[t]]^\ell$. If we write $v_i(t) = \sum_{n\geq 0} v_{i,n} t^n$,
then we obtain
\begin{equation} \label{descr of opers}
\on{Fun} \on{Op}_{^L G}(D) \simeq \C[v_{i,n}]_{i \in I, n\geq 0}.
\end{equation}
Note that the vector field $-t \pa_t$ acts naturally on $\on{Op}_{^L
G}(D)$ and defines a $\Z$-grading on $\on{Fun} \on{Op}_{^L G}(D)$
such that $\deg v_{i,n} = d_i+n+1$. The vector field $-\pa_t$ acts as
a derivation such that $-\pa_t \cdot v_{i,n} = -(d_i+n+1) v_{i,n+1}$.

\begin{thm}[\cite{FF:gd,F:wak}] \label{center}
There is a canonical isomorphism $$\zz(\G) \simeq \on{Fun} \on{Op}_{^L
G}(D)$$ of algebras which is compatible with the action of $\DerO$ and
$\AutO$.
\end{thm}

We use this result to describe the twist of $\zz(\G)$ by the
$\AutO$-torsor ${\mc A}ut_x$ of formal coordinates at a smooth point
$x$ of an algebraic curve $X$,
$$
\zz(\ghat)_x = {\mc A}ut_x \underset{\AutO}\times \zz(\G)
$$
(see Ch. 6 of \cite{FB} for more details). It follows from the
definition that the corresponding twist of $\on{Fun} \on{Op}_{^L
G}(D)$ by ${\mc A}ut_x$ is nothing but $\on{Fun} \on{Op}_{^L G}(D_x)$,
where $D_x$ is the disc around $x$. Therefore we obtain from
\thmref{center} an isomorphism
\begin{equation} \label{isom x}
\zz(\G)_x \simeq \on{Fun} \on{Op}_{^L G}(D_x).
\end{equation}

The module $\V_0$ has a natural $\Z$-grading defined by the formulas
$\deg v_0 = 0, \deg J^a_n = -n$, and it carries a translation operator
$T$ defined by the formulas $T v_0 = 0, [T,J^a_n] = -n
J^a_{n-1}$. \thmref{center} and the isomorphism \eqref{descr of opers}
imply that there exist non-zero vectors $S_i \in \V_0^{\g[[t]]}$ of
degrees $d_i+1, i \in I$, such that
$$
\zz(\G) = \C[T^n S_i]_{i \in I,n\geq 0} v_0.
$$
Then under the isomorphism of \thmref{center} we have $S_i \mapsto
v_{i,0}$, the $\Z$-gradings on both algebras get identified and the
action of $T$ on $\zz(\G)$ becomes the action of $-\pa_t$ on $\on{Fun}
\on{Op}_{^L G}(D)$. Note that the vector $S_1$ is nothing but the
vector \eqref{sugawara}, up to a non-zero scalar.

Recall from \cite{FB} that $\V_0$ is a vertex algebra, and $\zz(\G)$
is its commutative vertex subalgebra; in fact, it is the center of
$\V_0$. We will also need the center of the completed universal
enveloping algebra of $\ghat$ of critical level. This algebra is
defined as follows.

Let $U_{\ka_c}(\ghat)$ be the quotient of the universal enveloping
algebra $U(\ghat_{\ka_c})$ of $\ghat_{\ka_c}$ by the ideal generated
by $(K-1)$. Define its completion $\wt{U}_{\ka_c}(\ghat)$ as follows:
$$
\wt{U}_{\ka_c}(\ghat) = \underset{\longleftarrow}\lim \;
U_{\ka_c}(\ghat)/U_{\ka_c}(\ghat) \cdot (\g \otimes t^N\C[[t]]).
$$
It is clear that $\wt{U}_{\ka_c}(\ghat)$ is a topological algebra
which acts on all smooth $\ghat_{\ka_c}$-module. By definition, a
smooth $\ghat_{\ka_c}$-module is a $\ghat_{\ka_c}$-module such that
any vector is annihilated by $\g \otimes t^N\C[[t]]$ for sufficiently
large $N$, and $K$ acts as the identity.

Let $Z(\ghat)$ be the center of $\wt{U}_{\ka_c}(\ghat)$.

Denote by $\on{Fun} \on{Op}_{^L G}(D^\times)$ the algebra of regular
functions on the space $\on{Op}_{^L G}(D^\times)$ of $^L G$-opers on
the punctured disc $D^\times = \on{Spec} \C((t))$. In view of
\lemref{free}, it is isomorphic to the algebra of functions on the
space of $\ell$-tuples $(v_1(t),\ldots,v_\ell(t))$ of formal Laurent
series, i.e., the ind-affine space $\C((t))^\ell$. If we write $v_i(t)
= \sum_{n \in \Z} v_{i,n} t^n$, then we obtain that $\on{Fun}
\on{Op}_{^L G}(D)$ is isomorphic to the completion of the polynomial
algebra $\C[v_{i,n}]_{i \in I, n \in \Z}$ with respect to the topology
in which the basis of open neighborhoods of zero is formed by the
ideals generated by $v_{i,n}, i \in I, n \leq N$, for $N\leq 0$.

\begin{thm}[\cite{F:wak}] \label{center1}
There is a canonical isomorphism $$Z(\G) \simeq \on{Fun} \on{Op}_{^L
G}(D^\times)$$ of complete topological algebras which is compatible
with the action of $\DerO$ and $\AutO$.
\end{thm}

If $M$ is a smooth $\ghat_{\ka_c}$-module, then the action of
$Z(\ghat)$ on $M$ gives rise to a homomorphism
$$
Z(\ghat) \to \on{End}_{\ghat_{\ka_c}} M.
$$
For example, if $M=\V_0$, then using Theorems \ref{center} and
\ref{center1} we identify this homomorphism with the surjection
$$
\on{Fun} \on{Op}_{^L G}(D^\times) \twoheadrightarrow \on{Fun}
\on{Op}_{^L G}(D)
$$
induced by the natural embedding
$$
\on{Op}_{^L G}(D) \hookrightarrow \on{Op}_{^L G}(D^\times).
$$

Recall that the Harish-Chandra homomorphism identifies the center
$Z(\g)$ of $U(\g)$ with the algebra $(\on{Fun} \h^*)^W$ of polynomials
on $\h^*$ which are invariant with respect to the action of the Weyl
group $W$. Therefore a character $Z(\g) \to \C$ is the same as a
point in $\on{Spec} (\on{Fun} \h^*)^W$ which is the quotient $\h^*/W$.
For $\la \in \h^*$ we denote by $\varpi(\la)$ its projection onto
$\h^*/W$. In particular, $Z(\g)$ acts on $M_\la$ and $V_\la$ via its
character $\varphi(\la+\rho)$. We also denote by $I_\la$ the maximal
ideal of $Z(\g)$ equal to the kernel of the homomorphism $Z(\g) \to
\C$ corresponding to the character $\varphi(\la+\rho)$.

In what follows we will use the canonical identification between
$\h^*$ and the Cartan subalgebra $^L\h$ of the Langlands dual Lie
algebra $^L\g$. Recall that in \secref{reg sing} we defined the space
$\on{Op}_{^L G}^{\on{RS}}(D)$ of $^L G$-opers on $D^\times$ with
regular singularity and its subspace $\on{Op}_{^L
G}^{\on{RS}}(D)_{\la}$ of opers with residue $\varpi(-\la-\rho)$. We
also defined the subspace
$$
\on{Op}_{^L G}(D)_{\cla} \subset \on{Op}_{^L
G}^{\on{RS}}(D)_{\cla} \subset \on{Op}_{^L G}(D^\times)
$$
of those $^L G$-opers which have trivial monodromy. Here we
identify the coweights of the group $^L G$ with the weights of
$G$.

The following result is obtained by combining Theorem 12.4, Lemma 9.4
and Proposition 12.8 of \cite{F:wak}.

\begin{thm}\label{factor RS}\ %

{\em (1)} Let $\UU$ be the $\ghat_{\ka_c}$-module induced from the
$\g[[t]] \oplus \C K$-module $U(\g)$. Then the homomorphism
$Z(\ghat) \to \on{End}_{\ghat_{\ka_c}} \UU$ factors as
$$
Z(\ghat) \simeq \on{Fun} \on{Op}_{^L G}(D^\times)
\twoheadrightarrow \on{Fun} \on{Op}_{^L G}^{\on{RS}}(D) \to
\on{End}_{\ghat_{\ka_c}} \UU.
$$

{\em (2)}
Let $M$ be a $\g$-module on which the center $Z(\g)$ acts
via its character $\varpi(\la+\rho)$, and let $\M$ be the induced
$\ghat_{\ka_c}$-module. Then the homomorphism $Z(\ghat) \to
\on{End}_{\ghat_{\ka_c}} \M$ factors as follows
$$
Z(\ghat) \simeq \on{Fun} \on{Op}_{^L G}(D^\times)
\twoheadrightarrow \on{Fun} \on{Op}_{^L G}^{\on{RS}}(D)_{\la} \to
\on{End}_{\ghat_{\ka_c}} \M.
$$
Moreover, if $M=M_\la$, then the last map is an isomorphism
$$
\on{End}_{\ghat_{\ka_c}} \M \simeq \on{Fun} \on{Op}_{^L
G}^{\on{RS}}(D)_{\la}.
$$

{\em (3)}
For an integral dominant weight $\la \in \h^*$ the
homomorphism
$$
\on{Fun} \on{Op}_{^L G}(D^\times) \to
\on{End}_{\ghat_{\ka_c}} \V_\la
$$
factors as
$$
\on{Fun} \on{Op}_{^L G}(D^\times) \to \on{Fun} \on{Op}_{^L
G}(D)_\la \to \on{End}_{\G_{\ka_c}} \V_\la,
$$
and the last map is an isomorphism
$$
\on{End}_{\G_{\ka_c}} \V_\la \simeq \on{Fun} \on{Op}_{^L
G}(D)_\la.
$$

\end{thm}

\subsection{Example}

Let us consider the case $\g=\sw_2$ in more detail. Introduce the
Segal-Sugawara operators $S_n, n \in \Z$, by the formula
$$
S(z) = \sum_{n \in \Z} S_n z^{-n-2} = \frac{1}{2} \sum_a \Wick J^a(z)
J_a(z) \Wick \, ,
$$
where the normal ordering is defined as in \cite{FB}. Then the center
$Z(\sw_2)$ is the completion $\C[S_n]^{\sim}_{n \in \Z}$ of the
polynomial algebra $\C[S_n]_{n \in \Z}$ with respect to the topology
in which the basis of open neighborhoods of zero is formed by the
ideals of $S_n, n>N$, for $N \geq 0$.

We have the following diagram of (vertical) isomorphisms and
(horizontal) surjections
$$
\begin{CD}
Z(\su) @>>> \on{End}_{\ghat_{\ka_c}} \UU @>>> \on{End}_{\ghat_{\ka_c}}
\M_\la @>>> \on{End}_{\ghat_{\ka_c}} \V_\la \\ @VVV @VVV @VVV @VVV \\
\C[S_n]^{\sim}_{n \in \Z} @>>> \C[S_n]_{n \leq 0} @>>> \C[S_n]_{n \leq
0}/J_\la @>>> \C[S_n]_{n \leq
0}/J'_\la
\end{CD}
$$
where $J_\la$ is the ideal generated by $(S_0 - \frac{1}{4}
\la(\la+2))$ and $J'_\la$ is the ideal generated by $I_\la$ and
the polynomial $P_\la$ introduced at the end of \secref{reg sing}.

The space of $PGL_2$-opers on $D^\times$ is identified with the space
of projective connections of the form $\pa_t^2 - \sum_{n \in \Z} v_n
t^n$. The isomorphism of \thmref{center1} sends $S_n$ to
$v_{-n-2}$. The relevant spaces of $PGL_2$-opers with regular
singularities were described at the end of \secref{reg sing}, and
these descriptions agree with the above diagram and Theorem
\ref{factor RS}.

\subsection{The Gaudin algebra}

Now we are ready to identify the Gaudin algebra $\ZZ(\g)$ with the
algebra of functions on a certain space of opers on $\pone$.

Let $\on{Op}_{^L G}^{\on{RS}}(\pone)_{(z_i),\infty}$ be the space of
$^L G$-opers on $\pone$ with regular singularities at
$z_1,\ldots,z_N$ and $\infty$. For an arbitrary collection of weights
$\la_1,\ldots,\la_N$ and $\la_\infty$, let $$\on{Op}_{^L
G}^{\on{RS}}(\pone)_{(z_i),\infty;(\la_i),\la_\infty}$$ be its
subspace of those opers whose residue at the point $z_i$ (resp.,
$\infty$) is equal to $\varpi(-\la_i-\rho), i=1,\ldots,N$ (resp.,
$\varpi(-\la_\infty-\rho)$). Finally, if all of the weights
$\la_1,\ldots,\la_N,\la_\infty$ are dominant integral, we introduce a
subset $$\on{Op}_{^L G}(\pone)_{(z_i),\infty;(\la_i),\la_\infty}
\subset \on{Op}_{^L
G}^{\on{RS}}(\pone)_{(z_i),\infty;(\la_i),\la_\infty}$$ which consists
of those $^L G$-opers which have trivial monodromy representation.

On the other hand, for each collection of points $z_1,\ldots,z_N$ on
$\pone \bs \infty$ we have the Gaudin algebra
$$
\ZZ(\g) \subset \left( U(\g)^{\otimes N} \right)^G \simeq \left(
U(\g)^{\otimes (N+1)}/\g_{\on{diag}} \right)^G,
$$
where the second isomorphism is obtained by identifying
$U(\g)^{\otimes (N+1)}/\g_{\on{diag}}$ with $U(\g)^{\otimes N} \otimes 1$.

We have a homomorphism $c_i: Z(\g) \to U(\g) \to U(\g)^{\otimes
(N+1)}$ corresponding to the $i$th factor, for all $i = 1,\ldots,N$,
and a homomorphism $c_\infty: Z(\g) \to U(\g)^{\otimes (N+1)}$
corresponding to the $(N+1)$st factor. It is easy to see that the
images of $c_i, i=1,\ldots,N$, and $c_\infty$ belong to $\ZZ(\g)$. For
a collection of weights $\la_1,\ldots,\la_N$ and $\la_\infty$, let
$I_{(\la_i),\la_\infty}$ be the ideal of $\ZZ(\g)$ generated by
$c_i(I_{\la_i}), i=1,\ldots,N$, and $c_\infty(I_{\la_\infty})$. Let
${\mc Z}_{(z_i),\infty;(\la_i),\la_\infty}$ be the quotient of
$\ZZ(\g)$ by $I_{(\la_i),\la_\infty}$.

The algebra ${\mc Z}_{(z_i),\infty;(\la_i),\la_\infty}(\g)$ acts
on the space of $\g$-coinvariants in $$\bigotimes_{i=1}^N M_i
\otimes M_\infty,$$ where $M_i$ is a $\g$-module with central
character $\varpi(\la_i+\rho), i=1,\ldots,N$, and $M_\infty$ is a
$\g$-module with central character $\varpi(\la_\infty+\rho)$. In
particular, if all the weights $\la_1,\ldots,\la_N,\la_\infty$ are
dominant integral, then we can take as the $M_i$'s the
finite-dimensional irreducible modules $V_{\la_i}$ for
$i=1,\ldots,N$, and as $M_\infty$ the module $V_{\la_\infty}$. The
corresponding space of $\g$-coinvariants is isomorphic to the
space $$\left( \bigotimes_{i=1}^N V_{\la_i} \otimes V_{\la_\infty}
\right)^G$$ of $G$-invariants in $\bigotimes_{i=1}^N V_{\la_i}
\otimes V_{\la_\infty}$. Let $\ol{\mc
Z}_{(z_i),\infty;(\la_i),\la_\infty}(\g)$ be the image of the
algebra ${\mc Z}_{(z_i),\infty;(\la_i),\la_\infty}$ in $\on{End}
\left( \bigotimes_{i=1}^N V_{\la_i} \otimes V_{\la_\infty}
\right)^G$.

We have the following result.

\begin{thm} \label{descr of alg}
\hfill

{\em (1)} The algebra $\ZZ(\g)$ is isomorphic to the algebra of
functions on the space $\on{Op}_{^L
G}^{\on{RS}}(\pone)_{(z_i),\infty}$.

{\em (2)} The algebra ${\mc Z}_{(z_i),\infty;(\la_i),\la_\infty}(\g)$ is
isomorphic to the algebra of functions on the space $\on{Op}_{^L
G}^{\on{RS}}(\pone)_{(z_i),\infty;(\la_i),\la_\infty}$.

{\em (3)} For a collection of dominant integral weights
$\la_1,\ldots,\la_N,\la_\infty$, there is a surjective homomorphism
from the algebra of functions $\on{Op}_{^L
G}(\pone)_{(z_i),\infty;(\la_i),\la_\infty}$ to the algebra $\ol{\mc
Z}_{(z_i),\infty;(\la_i),\la_\infty}(\g)$.

\end{thm}

\begin{proof}
In \cite{FB} we defined, for any quasi-conformal vertex algebra
$V$, a smooth projective curve $X$, a set of points
$x_1,\ldots,x_N \in X$ and a collection of $V$-modules
$M_1,\ldots,M_N$, the space of coinvariants $H_V(X,(x_i),(M_i))$,
which is the quotient of $\bigotimes_{i=1}^N M_i$ by the action of
a certain Lie algebra. This construction (which is recalled in the
proof of Theorem 4.7 in \cite{F:opers}) is functorial in the
following sense. Suppose that we are given a homomorphism $W \to
V$ of vertex algebras (so that each $M_i$ becomes a $V$-module),
a collection $R_1,\ldots,R_N$ of $W$-modules and a collection of
homomorphisms of $W$-modules $M_i \to R_i$ for all
$i=1,\ldots,N$. Then the corresponding map $\bigotimes_{i=1}^N R_i
\to \bigotimes_{i=1}^N M_i$ gives rise to a map of the
corresponding spaces of coinvariants
$$
H_W(X,(x_i),(R_i)) \to H_V(X,(x_i),(M_i)).
$$

Suppose now that $W$ is the center of $V$ (see \cite{FB}). Then the
action of $W$ on any $V$-module $M$ factors through a homomorphism
$\wt{U}(W) \to \on{End}_{\C} M$, where $\wt{U}(W)$ is the enveloping
algebra of $W$ (see \cite{FB}, Sect.~4.3). Let $W(M)$ be the image of
this homomorphism. Then $H_W(X,(x_i),(W(M_i)))$ is an algebra, and we
obtain a natural homomorphism of algebras
\begin{equation} \label{functoriality}
H_W(X,(x_i),(W(M_i))) \to \on{End}_{\C} H_V(X,(x_i),(M_i)).
\end{equation}

If $V=\V_0$, then the center of $V$ is precisely the subspace
$\zz(\G)$ of $\g[[t]]$-invariant vectors in $\V_0$ (see
\cite{FB}). In particular, $\zz(\ghat)$ is a commutative vertex
subalgebra of $\V_0$. A module over the vertex algebra $\zz(\G)$ is
the same as a module over the topological algebra $\wt{U}(\zz(\G))$
which is nothing but the center $Z(\G)$ of $\wt{U}_{\ka_c}(\ghat)$
(see \cite{F:wak}, Sect.~11). The action of $Z(\G)$ on any
$\ghat_{\ka_c}$-module $M$ factors through the homomorphism $Z(\G)
\to \on{End}_{\ghat_{\ka_c}} M$. Let $Z(M)$ denote the image of this
homomorphism. Recall that we have
identified $Z(\G)$ with $\on{Fun} \on{Op}_{^L G}(D^\times)$ in
\thmref{center1}. For each $\ghat_{\ka_c}$-module $M$, the algebra
$Z(M)$ is a quotient of $\on{Fun} \on{Op}_{^L G}(D^\times)$, and hence
$\on{Spec} Z(M)$ is a subscheme in $\on{Op}_{^L G}(D^\times)$ which we
denote by $\on{Op}_{^L G}^M(D^\times)$.

The space of coinvariants $H_{\zz(\G)}(X,(x_i),Z(M_i))$ is computed in
the same way as in Theorem 4.7 of \cite{F:opers}:
\begin{equation} \label{as in 47}
H_{\zz(\G)}(X,(x_i),Z(M_i)) \simeq \on{Fun} \on{Op}_{^L
G}(X,(x_i),(M_i))
\end{equation}
where $\on{Op}_{^L G}(X,(x_i),(M_i))$ is the space of $^L G$-opers on
$X$ which are regular on $X\bs \{ x_1,\ldots,x_N \}$ and such that
their restriction to $D_x^\times$ belongs to $\on{Op}_{^L
G}^{M_i}(D_{x_i}^\times)$ for all $i=1,\ldots,N$.

Let $u$ be an additional point of $X$, different from
$x_1,\ldots,x_N$, and let us insert $\zz(\G) \simeq Z(\V_0)$ at
this point. Then by Theorem 10.3.1 of \cite{FB} we have an isomorphism
$$
H_{\zz(\G)}(X,(x_i),(Z(M_i))) \simeq
H_{\zz(\G)}(X;(z_i),u;(Z(M_i)),Z(\V_0)).
$$
Hence we obtain a homomorphism
\begin{equation} \label{from u}
\zz(\G)_u \simeq Z(\V_0)_u \to H_{\zz(\G)}(X,(x_i),(Z(M_i))).
\end{equation}
The corresponding homomorphism
$$
\on{Fun} \on{Op}_{^L G}(D_u) \to \on{Fun} \on{Op}_{^L
G}(X,(x_i),(M_i))
$$
(see formula \eqref{isom x}) is induced by the embedding
$$
\on{Op}_{^L G}(X,(x_i),(M_i)) \hookrightarrow \on{Op}_{^L G}(D_u)
$$
obtained by restricting an oper to $D_u$.

We apply this construction in the case when the curve $X$ is $\pone$,
the points are $z_1,\ldots,z_N$ and $\infty$, and the modules are
$\G_{\ka_c}$-modules $M_1,\ldots,M_N$ and $M_\infty$. It is proved in
\cite{FB} (see Theorem 9.3.3 and Remark 9.3.10) that the corresponding
space of coinvariants is the space of $\g_{\zn}$-coinvariants of
$\bigotimes_{i=1}^N \M_i \otimes \M_\infty$, which is the space
$H(M_1,\ldots,M_N,M_\infty)$ that we have computed in \lemref{iso}.

The homomorphism \eqref{functoriality} specializes to a homomorphism
\begin{equation} \label{point u}
H_{\zz(\G)}(\pone;(z_i),\infty;(Z(M_i)),Z(M_\infty)) \to \on{End}_\C
H(M_1,\ldots,M_N,M_\infty).
\end{equation}
Observe that by its very definition the homomorphism $$\Phi_u:
\zz(\G)_u \to \on{End}_\C H(M_1,\ldots,M_N,M_\infty)$$ constructed in
\secref{coinvariants} factors through the homomorphisms \eqref{point
u} and \eqref{from u}. Hence the image of $\Phi_u$ is a quotient of
the algebra
$$
H_{\zz(\G)}(\pone;(z_i),\infty;(Z(M_i)),Z(M_\infty)) \simeq \on{Fun}
\on{Op}_{^L G}(\pone;(z_i),\infty;(Z(M_i)),Z(M_\infty)),
$$
according to the isomorphism \eqref{as in 47}.

Let us specialize this result to our setting. First we suppose that
all $M_i$'s and $M_\infty$ are equal to $U(\g)$. By \thmref{factor
RS},(1), we have
$$
\on{Op}_{^L G}^{U(\g)}(D^\times) = \on{Op}_{^L G}^{\on{RS}}(D).
$$
Therefore we find that
\begin{equation} \label{part 1}
H_{\zz(\G)}(\pone;(z_i),\infty;(U(\g)),U(\g)) \simeq \on{Fun}
\on{Op}_{^L G}^{\on{RS}}(\pone)_{(z_i),\infty}.
\end{equation}

Next, by \thmref{factor RS},(2), we have
$$
\on{Op}_{^L G}^{U(\g)/I_\la}(D^\times) = \on{Op}_{^L
G}^{\on{RS}}(D)_\la.
$$
Therefore
\begin{equation} \label{part 2}
H_{\zz(\G)}(\pone;(z_i),\infty;(U(\g)/I_{\la_i}),U(\g)/I_{\la_\infty})
\simeq \on{Fun} \on{Op}_{^L
G}^{\on{RS}}(\pone)_{(z_i),\infty;(\la_i),\la_\infty}.
\end{equation}

Finally, by \thmref{factor RS},(3), we have
$$
\on{Op}_{^L G}^{\V_\la}(D^\times) = \on{Op}_{^L G}(D)_\la.
$$
Therefore
\begin{equation} \label{part 3}
H_{\zz(\G)}(\pone;(z_i),\infty;(\V_{\la_i}),\V_{\la_\infty})
\simeq \on{Fun} \on{Op}_{^L
G}(\pone)_{(z_i),\infty;(\la_i),\la_\infty}.
\end{equation}

Moreover, in all three cases the homomorphism $\Phi_u$ is just the
natural homomorphism from $\on{Fun} \on{Op}_{^L G}(D_u)$ to the above
algebras of functions that is induced by the restrictions of the
corresponding opers to the disc $D_u$.

Now consider the homomorphism \eqref{point u} in the case of the space
of coinvariants given by the left-hand side of formula \eqref{part 1},
$$
\on{Fun} \on{Op}_{^L G}^{\on{RS}}(\pone)_{(z_i),\infty} \to
\on{End}_\C U(\g)^{\otimes (N+1)}/\g_{\on{diag}} \simeq \on{End}_\C
U(\g)^{\otimes N},
$$
where the second isomorphism we use the identification of
$U(\g)^{\otimes (N+1)}/\g_{\on{diag}}$ with $U(\g)^{\otimes N}$
corresponding to the first $N$ factors. It follows from the explicit
computation of this map given above that its image belongs to $$\left(
U(\g)^{\otimes N} \right)^G \subset \on{End}_\C U(\g)^{\otimes N},$$
and so we have a homomorphism
\begin{equation} \label{inj1}
\on{Fun} \on{Op}_{^L G}^{\on{RS}}(\pone)_{(z_i),\infty} \to
\left( U(\g)^{\otimes N} \right)^G.
\end{equation}
By definition, the Gaudin algebra $\ZZ(\g)$ is the image of this
homomorphism.

Likewise, we obtain from formula \eqref{part 2} that the homomorphism
\eqref{point u} gives rise to a homomorphism
\begin{equation} \label{inj2}
\on{Fun}
\on{Op}_{^LG}^{\on{RS}}(\pone)_{(z_i),\infty;(\la_i),\la_\infty}
\to \left( U(\g)^{\otimes (N+1)}/\g_{\on{diag}}
\right)^G/I_{(\la_i),\la_\infty},
\end{equation}
whose image is the algebra
${\mc Z}_{(z_i),\infty;(\la_i),\la_\infty}(\g)$.

Finally, formula \eqref{part 3} gives us a homomorphism
$$
\on{Fun} \on{Op}_{^L G}(\pone)_{(z_i),\infty;(\la_i),\la_\infty}
\to \on{End}_\C \left( \bigotimes_{i=1}^N V_{\la_i} \otimes
V_{\la_\infty} \right),
$$
whose image is the algebra $\ol{\mc
Z}_{(z_i),\infty;(\la_i),\la_\infty}(\g)$. This proves part (3) of the
theorem.

To prove parts (1) and (2), it remains to show that the homomorphisms
\eqref{inj1} and \eqref{inj2} are injective. It is sufficient to prove
that the latter is injective. To see that, we pass to the associate
graded spaces on both sides with respect to natural
filtrations which we now describe.

According to the identification given in formula \eqref{repr}, the
algebra of functions on
$\on{Op}_{^LG}^{\on{RS}}(\pone)_{(z_i),\infty;(\la_i),\la_\infty}$ is
filtered, and the corresponding associated graded algebra is the
algebra of functions on the vector space
$$
C^{\on{RS}}_{(z_i),\infty} = \bigoplus_{i=1}^\ell
\Gamma(\pone,\Omega^{\otimes(d_i+1)}(
-d_i z_1-\ldots-d_i z_N-d_i\infty)),
$$
where $\Omega$ is the canonical line bundle on $\pone$. The algebra
$$
\left( U(\g)^{\otimes (N+1)}/\g_{\on{diag}}
\right)^G/I_{(\la_i),\la_\infty}
$$
carries a PBW filtration, and the associated graded is the algebra of
functions on the space $\mu^{-1}((T^* G/B)^{N+1})/G$, where $\mu: (T^*
G/B)^{N+1} \to \g^*$ is the moment map corresponding to the diagonal
action of $G$ on $(T^* G/B)^{N+1}$. The two filtrations are compatible
according to \cite{F:wak}. The corresponding homomorphism of the
associate graded algebras
\begin{equation} \label{hom of alg of fun}
\on{Fun} C^{\on{RS}}_{(z_i),\infty} \to \on{Fun} \mu^{-1}((T^*
G/B)^{N+1})/G
\end{equation}
is induced by a map
$$
h_{(z_i),\infty}: \mu^{-1}((T^* G/B)^{N+1})/G \to
C^{\on{RS}}_{(z_i),\infty},
$$
that we now describe.

Let us identify the tangent space to a point $gB \subset G/B$ with
$(\g/g \bb g^{-1})^* \simeq g \n g^{-1}$. Then a point in
$\mu^{-1}((T^* G/B)^{N+1})/G$ consists of an $(N+1)$-tuple of
points $g_i B$ of $G/B$ and an $(N+1)$-tuple of vectors
$(\eta_i)$, where $\eta_i \in g_i \n g_i^{-1} \subset \g$ such
that $\sum_{i=1}^{N+1} \eta_i = 0$, considered up to simultaneous
conjugation by $G$.

We attach to it the $\g$-valued one-form
$$
\eta = \sum_{i=1}^N \frac{\eta_i}{t-z_i} dt
$$
on $\pone$ with poles at $z_1,\ldots,z_N,\infty$. Let
$P_1,\ldots,P_\ell$ be generators of the algebra of $G$-invariant
polynomials on $\g$ of degrees $d_i+1$. Then
$$
h_{(z_i),\infty}((g_i),(\eta_i)) = (P_i(\eta))_{i=1}^\ell \in
C^{\on{RS}}_{(z_i),\infty}.
$$

The space $\mu^{-1}((T^* G/B)^{N+1})/G$ is identified with the moduli
space of Higgs fields on the trivial $G$-bundle with parabolic
structures at $z_1,\ldots,z_N,\infty$, and the map $h_{(z_i),\infty}$
is nothing but the Hitchin map (see \cite{ER}). The Hitchin map is
known to be proper, so in particular it is surjective (see, e.g.,
\cite{Markman}). Therefore the corresponding homomorphism \eqref{hom
of alg of fun} of algebras of functions is injective. This implies
that the homomorphism \eqref{inj2} is also injective and completes the
proof of the theorem.
\end{proof}

This theorem has an important application to the question of
simultaneous diagonalization of generalized Gaudin Hamiltonians,
or equivalently, of the commutative algebra $\ZZ(\g)$, on the
tensor product $\bigotimes_{i=1}^N M_i$ of $\g$-modules. Indeed,
the joint eigenvalues of the generalized Gaudin Hamiltonians on
any eigenvector in $\bigotimes_{i=1}^N M_i$ correspond to a point
in the spectrum of the algebra $\ZZ(\g)$, which, according to
\thmref{descr of alg},(1), is a point of the space $\on{Op}_{^L
G}^{\on{RS}}(\pone)_{(z_i),\infty}$.

If we assume in addition that each of the modules $M_i$ admits a
central character $\varpi(\la_i+\rho)$ (for instance, if $M_i$ is the
Verma module $M_{\la_i}$) and we are looking for eigenvectors in the
component of $\bigotimes_{i=1}^N M_i$ corresponding to the central
character $\varpi(-\la_\infty-\rho)$ with respect to the diagonal
action of $\g$, then the joint eigenvalues define a point in the
spectrum of the algebra ${\mc
Z}_{(z_i),\infty;(\la_i),\la_\infty}(\g)$, i.e., a point of
$\on{Op}_{^L G}^{\on{RS}}(\pone)_{(z_i),\infty;(\la_i),\la_\infty}$.

Finally, for a collection of dominant integral weights
$\la_1,\ldots,\la_N,\la_\infty$, the joint eigenvalues of the
generalized Gaudin Hamiltonians on $(\bigotimes_{i=1}^N V_{\la_i} \otimes
V_{\la_\infty})^G$ is a point in the spectrum of the algebra $\ol{\mc
Z}_{(z_i),\infty;(\la_i),\la_\infty}(\g)$, which is a point of $\on{Op}_{^L
G}(\pone)_{(z_i),\infty;(\la_i),\la_\infty}$.

A natural question is whether, conversely, one can attach to a $^L
G$-oper on $\pone$ with regular singularities at
$z_1,\ldots,z_N,\infty$ (and satisfying additional conditions as
above) an eigenvector in $\bigotimes_{i=1}^N M_i$ with such
eigenvalues. It turns out that for general modules this is not true,
but if these modules are finite-dimensional, then we conjecture that
it is true. In order to construct the eigenvectors we use the
procedure called Bethe Ansatz. As shown in \cite{FFR}, this procedure
may be cast in the framework of coinvariants that we have discussed in
this section, using the Wakimoto modules over $\ghat_{\ka_c}$. We will
explain that in \secref{wak and ba} and \secref{last}. But first we
need to introduce Miura opers and Cartan connections.

\section{Miura opers and Cartan connections} \label{miura opers}

By definition (see \cite{F:wak}, Sect. 10.3), a {\em Miura $G$-oper}
on $X$ (which is a smooth curve or a disc) is a quadruple
$(\F,\nabla,\F_B,\F'_B)$, where $(\F,\nabla,\F_B)$ is a $G$-oper on
$X$ and $\F'_B$ is another $B$-reduction of $\F$ which is preserved
by $\nabla$.

We denote the space of Miura $G$-opers on $X$ by $\on{MOp}_G(X)$.

\subsection{Miura opers and flag manifolds}

A $B$-reduction of $\F$ which is preserved by the connection $\nabla$
is uniquely determined by a $B$-reduction of the fiber $\F_x$ of $\F$
at any point $x \in X$ (in the case when $U=D$, $x$ has to be the origin
$0 \in D$). The set of such reductions is the twist
\begin{equation} \label{twist of flag}
(G/B)_{\F_x} = \F_x \us{G}\times G/B = \F'_{B,x} \us{B}\times G/B =
(G/B)_{\F'_{B,x}}
\end{equation}
of the flag manifold $G/B$. If $X$ is a curve or a disc and the oper
connection has a regular singularity and trivial monodromy
representation, then this connection gives us a global (algebraic)
trivialization of the bundle $\F$. Then any $B$-reduction of the
fiber $\F_x$ gives rise to a global (algebraic) $B$-reduction of
$\F$. Thus, we obtain:

\begin{lem} \label{isom with flags}
Suppose that we are given an oper $\tau$ on a curve $X$ (or on the
disc) such that the oper connection has a regular singularity and
trivial monodromy. Then for each $x \in X$ there is a canonical
isomorphism between the space of Miura opers with the underlying oper
$\tau$ and the twist $(G/B)_{\F_{B,x}}$.
\end{lem}

Recall that the $B$-orbits in $G/B$, known as the {\em Schubert
cells}, are parameterized by the Weyl group $W$ of $G$. Let $w_0$ be
the longest element of the Weyl group of $G$. Denote the orbit $B
w_0 w B \subset G/B$ by $S_w$ (so that $S_1$ is the open
orbit). We obtain from the second description of $(G/B)_{\F_x}$ given
in formula \eqref{twist of flag} that $(G/B)_{\F_x}$ decomposes into a
union of locally closed subvarieties $S_{w,\F'_{B,x}}$, which are the
$\F'_{B,x}$-twists of the Schubert cells $S_w$. The $B$-reduction
$\F'_{B,x}$ defines a point in $(G/B)_{\F'_{B,x}}$ which is the one
point orbit $S_{w_0,\F'_{B,x}}$. We will say that the
$B$-reductions $\F_{B,x}$ and $\F'_{B,x}$ are in {\em relative
position} $w$ with if $\F_{B,x}$ belongs to $S_{w,\F'_{B,x}}$. In
particular, if it belongs to the open orbit $S_{1,\F'_{B,x}}$, we will
say that $\F_{B,x}$ and $\F'_{B,x}$ are in generic position.

A Miura $G$-oper is called {\em generic} at the point $x \in X$ if
the $B$-reductions $\F_{B,x}$ and $\F'_{B,x}$ of $\F_x$ are in
generic position. In other words, $\F_{B,x}$ belongs to the stratum
$\on{Op}_G(X) \times S_{1,\F'_{B,x}} \subset \on{MOp}_G(X)$. Being
generic is an open condition. Therefore if a Miura oper is generic at
$x \in X$, then there exists an open neighborhood $U$ of $x$ such that
it is also generic at all other points of $U$. We denote the space of
generic Miura opers on $U$ by $\on{MOp}_G(U)_{\on{gen}}$.

\begin{lem} \label{is generic}
Suppose we are given a Miura oper on the disc $D_x$ around a point $x
\in X$. Then its restriction to the punctured disc $D_x^\times$ is
generic.
\end{lem}

\begin{proof}
Since being generic is an open condition, we obtain that if a Miura
oper is generic at $x$, it is also generic on the entire $D_x$. Hence
we only need to consider the situation where the Miura oper is not
generic at $x$, i.e., the two reductions $\F_{B,x}$ and $\F'_{B,x}$
are in relative position $w \neq 1$. Let us trivialize the $B$-bundle
$\F_B$, and hence the $G$-bundle $\F_G$ over $D_x$. Then $\nabla$
gives us a connection on the trivial $G$-bundle which we can bring to
the canonical form
$$
\nabla = \pa_t + p_{-1} + \sum_{j=1}^\ell v_j(t) \cdot p_j
$$
(see \lemref{free}). It induces a connection on the trivial
$G/B$-bundle. We are given a point $gB$ in the fiber of the latter
bundle which lies in the orbit $S_w = B w_0 w B$, where $w \neq
1$. Consider the horizontal section whose value at $x$ is $gB$, viewed
as a map $D_x \to G/B$. We need to show that the image of this map
lies in the open $B$-orbit $S_1 = B w_0 B$ over $D_x^\times$, i.e.,
it does not lie in the orbit $S_y$ for any $y \neq 1$.

Suppose that this is not so, and the image of the horizontal section
actually lies in the orbit $S_y$ for some $y \neq 1$. Since all
$B$-orbits are $H$-invariant, we obtain that the same would be true
for the horizontal section with respect to the connection $\nabla' = h
\nabla h^{-1}$ for any constant element of $H$. Choosing $h =
\crho(a)$ for $a \in \C^\times$, we can bring the connection to the
form
$$
\pa_t + a^{-1} p_{-1} + \sum_{j=1}^\ell a^{d_j} v_j(t) \cdot p_j.
$$
Changing the variable $t$ to $s = a^{-1} t$, we obtain the connection
$$
\pa_s + p_{-1} + \sum_{j=1}^\ell a^{d_j+1} v_j(t),
$$
so choosing small $a$ we can make the functions $v_j(t)$ arbitrarily
small. Therefore without loss of generality we can consider the case
when our connection operator is $\nabla = \pa_t + p_{-1}$.

In this case our assumption that the horizontal section lies in $S_y,
y \neq 1$, means that the vector field $\xi_{p_{-1}}$ corresponding to
the infinitesimal action of $p_{-1}$ on $G/B$ is tangent to an orbit
$S_y, y \neq 1$, in the neighborhood of some point $gB$ of $S_w
\subset G/B, w \neq 1$. But then, again because of the $H$-invariance
of the $B$-orbits, the vector field $\xi_{h p_{-1} h^{-1}}$ is also
tangent to this orbit for any $h \in H$. For any $i=1\ldots,\ell$,
there exists a one-parameter subgroup $h_{\ep}^{(i)}, \ep \in
\C^\times$ in $H$, such that $\underset{\ep \to 0}\lim \; \ep p_{-1}
\ep^{-1} = f_i$. Hence we obtain that each of the vector fields
$\xi_{f_i}, i=1\ldots,\ell$, is tangent to the orbit $S_y, y \neq 1$,
in the neighborhood of $gB \in S_w, w \neq 1$. But then all
commutators of these vectors fields are also tangent to this
orbit. Hence we obtain that all vector fields of the form $\xi_p, p
\in \n_-$, are tangent to $S_y$ in the neighborhood of $gB \in
S_w$.

Consider any point of $G/B$ that does not belong to the open dense
orbit $S_1$. Then the quotient of the tangent space to this point by
the tangent space to the $B$-orbit passing through this point is
non-zero and the vector fields from the Lie algebra $\n_-$ map
surjectively onto this quotient. Therefore they cannot be tangent to
the orbit $S_y, y \neq 1$, in a neighborhood of $gB$. Therefore our
Miura oper is generic on $D_x^\times$.
\end{proof}

This lemma shows that any Miura oper on any smooth curve $X$ is
generic over an open dense subset.

\subsection{Cartan connections} \label{cartan conn}

Introduce the $H$-bundle $\Omega^{\crho}$ on $X$ which is uniquely
determined by the following property: for any character $\la: H \to
\C^\times$, the line bundle $\Omega^{\crho} \us{H}\times \la$
associated to the corresponding one-dimensional representation of $H$
is $\Omega^{\langle \la,\crho \rangle}$.

Explicitly, connections on $\Omega^{\crho}$ may be described as
follows. If we choose a local coordinate $t$ on $X$, then we
trivialize $\Omega^{\crho}$ and represent the connection as an
operator $\pa_t + {\mb u}(t)$, where ${\mb u}(t)$ is an $\h$-valued
function on $X$. If $s$ is another coordinate such that
$t=\varphi(s)$, then this connection will be represented by the
operator
\begin{equation} \label{trans for conn}
\pa_s + \varphi'(s) {\mb u}(\varphi(s)) - \crho \cdot
\frac{\varphi''(s)}{\varphi'(s)}.
\end{equation}

Let $\on{Conn}(\Omega^{\crho})_X$ be the space of connections on the
$H$-bundle $\Omega^{\crho}$ on $X$. When no confusion can arise, we
will simply write $\on{Conn}_X$. We define a map $${\mb b}_X:
\on{Conn}_X \to \on{MOp}_G(X)_{\on{gen}}.$$

Suppose we are given a connection $\ol\nabla$ on the $H$-bundle
$\Omega^{\crho}$ on $D$. We associate to it a generic Miura oper as
follows. Let us choose a splitting $H \to B$ of the homomorphism $B
\to H$ and set $\F = \Omega^{\crho} \underset{H}\times G, \F_B =
\Omega^{\crho} \underset{H}\times B$, where we consider the adjoint
action of $H$ on $G$ and on $B$ obtained through the above
splitting. The choice of the splitting also gives us the opposite
Borel subgroup $B_-$, which is the unique Borel subgroup in generic
position with $B$ containing $H$. Let again $w_0$ be the longest
element of the Weyl group of $\g$. Then $w_0 B$ is a $B$-torsor
equipped with a left action of $H$, so we define the $B$-subbundle
$\F'_B$ of $\F$ as $\Omega^{\crho} \underset{H}\times w_0 B$.

Observe that the space of connections on $\F$ is isomorphic to the
direct product
$$
\on{Conn}_X \times \bigoplus_{\al \in \De}
\Gamma(X,\Omega^{\al(\crho) + 1}).
$$
Its subspace corresponding to negative simple roots is isomorphic to
the tensor product of $\left( \bigoplus_{i=1}^\ell \g_{-\al_i}
\right)$ and $\on{Fun} X$. Having chosen a basis element $f_i$ of
$\g_{-\al_i}$ for each $i=1,\ldots,\ell$, we now construct an element
$p_{-1} = \sum_{i=1}^\ell f_i \otimes 1$ of this space. Now we set
$\nabla = \ol\nabla + p_{-1}$. By construction, $\nabla$ has the
correct relative position with the $B$-reduction $\F_B$ and preserves
the $B$-reduction $\F'_B$. Therefore the quadruple
$(\F,\nabla,\F_B,\F'_B)$ is a generic Miura oper on $X$. We define the
morphism ${\mb b}_X$ by setting ${\mb b}_X(\ol{\nabla}) =
(\F,\nabla,\F_B,\F'_B)$.

This map is independent of the choice of a splitting $H \to B$ and of
the generators $f_i, i=1,\ldots,\ell$.

\begin{prop}[\cite{F:wak},Prop. 10.4] \label{map beta}
The map ${\mb b}_X$ is an isomorphism of algebraic varieties
$$
\on{Conn}_X \to \on{MOp}_G(X)_{\on{gen}}.
$$
\end{prop}

Thus, generic Miura opers are the same as Cartan connections, which
are much simpler objects than opers.

The composition $\ol{\mb b}_X$ of ${\mb b}_X$ and the forgetful map
$\on{MOp}_G(X)_{\on{gen}} \to \on{Op}_G(X)$ is called the {\em Miura
transformation}.

For example, in the case of $\g=\sw_2$, we have a connection
$\ol{\nabla} = \pa_t - u(t)$ on the line bundle $\Omega^{1/2}$
(equivalently, a connection $\ol{\nabla}^t = \pa_t + u(t)$ on
$\Omega^{-1/2}$), and the Miura transformation assigns to this
connection the $PGL_2$-oper
$$
\pa_t + \left( \begin{array}{ccccc}
- u(t) & 0 \\
1 & u(t)
\end{array} \right).
$$
The oper $B$ reduction $\F_B$ corresponds to the upper triangular
matrices, and the Miura $B$-reduction $\F'_B$ corresponds to the
lower triangular matrices. The corresponding projective connection is
$$
\pa^2_t - v(t) = (\pa_t - u(t))(\pa_t + u(t)),
$$
i.e.,
$$
u(t) \mapsto v(t) = u(t)^2 - u'(t).
$$

\subsection{Singularities of Cartan connections} \label{miura}

Consider the Miura transformation $\ol{\mb b}_{D_x^\times}$ in the
case of the punctured disc $D^\times_x$,
$$
\ol{\mb b}_{D_x^\times}: \Conn_{D_x^\times} \to
\on{Op}_G(D_x^\times).
$$
Let $\Conn^{\on{RS}}_{D_x} \subset \Conn_{D_x^\times}$ be the space of
all connections on the
$H$-bundle $\Omega^{\crho}$ on $D_x$ with regular singularity, i.e.,
those for which the connection operator has the form
$$
\ol{\nabla} = \pa_t + \frac{\cla}{t} + \sum_{n\geq 0} u_n t^n.
$$
We define a map
$$
\on{res}_{\h}: \Conn^{\on{RS}}_{D} \to \h
$$
assigning to such a connection its residue $\cla$.

It follows from the definition of the Miura transformation $\ol{\mb
b}_{D_x^\times}$ that its restriction to $\Conn^{\on{RS}}_{D_x}
\subset \Conn_{D_x^\times}$ takes values in
$\on{Op}_G^{\on{RS}}(D_x)$. Hence we obtain a morphism
$$
\ol{\mb b}_x^{\on{RS}}: \Conn^{\on{RS}}_{D_x} \to
\on{Op}_G^{\on{RS}}(D_x).
$$
Explicitly, after choosing a coordinate $t$ on $D_x$, we can write
$\ol\nabla$ as $\pa_t + t^{-1} {\mb u}(t)$, where ${\mb u}(t) \in
\h[[t]]$. Its residue is ${\mb u}(0)$. Then the corresponding oper
with regular singularity is by definition the $N((t))$-equivalence
class of the operator
$$
\nabla = \pa_t + p_{-1} + t^{-1} {\mb u}(t),
$$
which is the same as the $N[[t]]$-equivalence class of the operator
$$
\crho(t) \nabla \crho(t)^{-1} = \pa_t + t^{-1} (p_{-1} - \crho + {\mb
u}(t)),
$$
so it is indeed an oper with regular singularity.

Therefore it follows from the definition that we have a commutative
diagram
\begin{equation} \label{comm diag}
\begin{CD}
\Conn^{\on{RS}}_{D_x} @>{\ol{\mb b}_x^{\on{RS}}}>>
\on{Op}_G^{\on{RS}}(D_x) \\ @V{\on{res}_{\h}}VV @VV{\on{res}}V \\ \h
@>>> \h/W
\end{CD}
\end{equation}
where the lower horizontal map is the composition of the map $\cla
\mapsto \cla - \crho$ and the projection $\varpi: \h \to \h/W$.

Now let $\on{Conn}^{\on{reg}}_{D_x,\cla}$ be the preimage under
$\ol{\mb b}^{\on{RS}}_x$ of the subspace $\on{Op}_G(D_x)_{\cla}
\subset \on{Op}_G^{\on{RS}}(D_x)$ of $G$-opers on $D_x^\times$
with regular singularity, residue $\varpi(-\cla-\crho)$ and
trivial monodromy. By the commutativity of the above diagram, a
connection in $\on{Conn}^{\on{reg}}_{D_x,\cla}$ necessarily has
residue of the form $-w(\cla+\crho)+\crho$ for some element $w$ of
the Weyl group of $G$, so that $\on{Conn}^{\on{reg}}_{D_x,\cla}$
is the disjoint union of its subsets
$\on{Conn}^{\on{reg}}_{D_x,\cla,w}$ consisting of connections with
residue $-w(\cla+\crho)+\crho$.

The restriction of $\ol{\mb b}_x^{\on{RS}}$ to
$\on{Conn}^{\on{reg}}_{D_x,\cla,w}$ is a map
$$
\ol{\mb b}_{\cla,w}: \on{Conn}^{\on{reg}}_{D_x,\cla,w} \to
\on{Op}_G(D_x)_{\cla}.
$$
Let us recall that by construction of the Miura transformation
$\ol{\mb b}$, each oper on $D_x^\times$ which lies in the image of the
map $\ol{\mb b}$ (hence in particular, in the image of $\ol{\mb
b}_{\cla,w}$) carries a canonical horizontal $B$-reduction $$\F'_B =
\Omega^{\crho} \underset{H}\times w_0 B$$ (i.e., it carries a
canonical structure of Miura oper on $D_x^\times$). But if this oper
is in the image of $\ol{\mb b}_{\cla,w}$, i.e., belongs to
$\on{Op}_G(D_x)_{\cla}$, then the oper $B$-reduction $\F_B$ (and hence
the oper bundle $\F$) has a canonical extension to a $B$-bundle on the
entire disc $D_x$, namely, one for which the oper connection has the
form \eqref{psi la}. Therefore the $B$-reduction $\F'_B =
\Omega^{\crho} \underset{H}\times w_0 B$ may also be extended to
$D_x$.

Therefore we can lift $\ol{\mb b}_{\cla,w}$ to a map
$$
{\mb b}_{\cla,w}: \on{Conn}^{\on{reg}}_{D_x,\cla,w} \to
\on{MOp}_G(D_x)_{\cla}.
$$
Let $\on{MOp}_G(D_x)_{\cla,w} \subset \on{MOp}_G(D_x)_{\cla}$ be the
subvariety of those Miura opers of coweight $\cla$ which have relative
position $w$ at $x$. Then
$$
\on{MOp}_G(D_x)_{\cla,w} \simeq
\on{Op}_G(D_x)_{\cla} \times S_{w,\F'_{B,x}}.
$$

The following result is due to D. Gaitsgory and myself \cite{FG}
(see \cite{F:opers}, Proposition~2.9).

\begin{prop} \label{isom w}
For each $w \in W$ the morphism ${\mb b}_{\cla,w}$ is an isomorphism
between the varieties $\on{Conn}^{\on{reg}}_{D_x,\cla,w}$ and
$\on{MOp}_G(D_x)_{\cla,w}$.
\end{prop}

\begin{proof}
First we observe that at the level of points the map defined by ${\mb
b}_{\cla,w}, w \in W$, from the union of
$\on{Conn}^{\on{reg}}_{D_x,\cla,w}, w \in W$, to
$\on{MOp}_G(D_x)_{\cla}$, is a bijection. Indeed, by \propref{map
beta} we have a map taking a Miura oper from $\on{MOp}_G(D_x)_{\cla}$,
considered as a Miura oper on the punctured disc $D_x^\times$, to a
connection $\ol\nabla$ on the $H$-bundle $\Omega^{\crho}$ over
$D_x^\times$. We have shown above that $\ol\nabla$ has regular
singularity at $x$ and that its residue is of the form
$-w(\cla+\crho)+\crho, w \in W$. Thus, we obtain a map from the set of
points of $\on{MOp}_G(D_x)_{\cla}$ to the union of
$\on{Conn}^{\on{reg}}_{D_x,\cla,w}, w \in W$, and by \propref{map
beta} it is a bijection.

It remains to show that if the Miura oper belongs to
$\on{MOp}_G(D_x)_{\cla,w}$, then the corresponding connection has
residue precisely $-w(\cla+\crho)+\crho$.

Thus, we are given a $G$-oper $(\F,\nabla,\F_B,\F'_B)$ of coweight
$\cla$. Let us choose a trivialization of the $B$-bundle $\F_B$. Then
the connection operator reads
\begin{equation} \label{conn op oper}
\nabla = \pa_t + \sum_{i=1}^\ell t^{\langle \al_i,\cla \rangle} f_i +
{\mb v}(t), \qquad {\mb v}(t) \in \bb[[t]].
\end{equation}
Suppose that the horizontal $B$-reduction $\F'_B$ of our Miura oper
has relative position $w$ with $\F_B$ at $x$ (see \secref{miura opers}
for the definition of relative position). We need to show that the
corresponding connection on $\F'_H \simeq \Omega^{\crho}$ has residue
$-w(\cla+\crho)+\crho$.

This is equivalent to the following statement. Let $\Phi(t)$ be the
$G$-valued solution of the equation
\begin{equation} \label{again equation}
\left( \pa_t + \sum_{i=1}^\ell t^{\langle \al_i,\cla \rangle} f_i +
{\mb v}(t) \right) \Phi(t) = 0,
\end{equation}
such that $\Phi(0) = 1$. Since the connection operator is regular at
$t=0$, this solution exists and is unique. Then $\Phi(t) w^{-1} w_0$
is the unique solution of the equation \eqref{again equation} whose
value at $t=0$ is equal to $w^{-1} w_0$.

By \lemref{is generic}, we have
$$
\Phi(t) w^{-1} w_0 = X_w(t) Y_w(t) Z_w(t) w_0,
$$
where
$$
X_w(t) \in N((t)), \qquad Y_w(t) \in H((t)), \qquad Z_w(t) \in
N_-((t)).
$$
We can write $Y_w(t) = \cmu_w(t) \wt{Y}_w(t)$, where $\cmu_w$ is a coweight
and $\wt{Y}_w(t) \in H[[t]]$.

Since the connection $\nabla$
preserves $$\Phi(t) w_0 \bb_+ w_0 \Phi(t)^{-1} = \Phi(t)
\bb_- \Phi(t)^{-1},$$ the connection $X(t)_w^{-1} \nabla
X_w(t)$ preserves $$Y_w(t) Z_w(t) \bb_- Z_w(t)^{-1} Y_w(t)^{-1} =
\bb_-,$$ and therefore has the form
$$
\pa_t + \sum_{i=1}^\ell t^{\langle \al_i,\cla \rangle} f_i -
\frac{\cmu_w}{t} + {\mb u}(t), \qquad {\mb u}(t) \in \h[[t]].
$$
By conjugating it with $\cla(t)$ we obtain a connection
$$
\pa_t + p_{-1} - \frac{\cla+\cmu_w}{t} + {\mb u}(t), \qquad {\mb
u}(t) \in \h[[t]].
$$
Therefore we need to show that
\begin{equation} \label{desired}
\cmu_w = w(\cla+\crho) - (\cla+\crho).
\end{equation}

To see that, let us apply the identity $\Phi(t) w^{-1} = X_w(t) Y_w(t)
Z_w(t)$ to a non-zero vector $v_{w_0(\nu)}$ of weight $w_0(\nu)$ in a
finite-dimensional irreducible $\g$-module $V_\nu$ of highest weight
$\nu$ (so that $v_{w_0(\nu)}$ is a lowest weight vector and hence is
unique up to scalar). The right-hand side will then be equal to a
$P(t) v_{w_0(\nu)}$ plus the sum of terms of weights greater than
$w_0(\nu)$, where $P(t) = c t^{\langle w_0(\nu),\cmu_w \rangle}, c
\neq 0$, plus the sum of terms of higher degree in $t$. Applying the
left-hand side to $v_{w_0(\nu)}$, we obtain $\Phi(t) v_{w^{-1}
w_0(\nu)}$, where $v_{w^{-1} w_0(\nu)} \in V_\nu$ is a non-zero vector
of weight $w^{-1} w_0(\nu)$ which is also unique up to a scalar.

Thus, we need to show that the coefficient with which $v_{w_0(\nu)}$
enters the expression $\Phi(t) v_{w^{-1} w_0(\nu)}$ is a polynomial in
$t$ whose lowest degree is equal to
$$
\langle w_0(\nu),w(\cla+\crho) - (\cla+\crho) \rangle,
$$
because if this is so for all dominant integral weights $\nu$, then we
obtain the desired equality \eqref{desired}. But this formula is easy
to establish. Indeed, from the form \eqref{conn op oper} of the oper
connection $\nabla$ it follows that we can obtain a vector
proportional to $v_{w_0}$ by applying the operators $\frac{1}{\langle
\al_i,\cla \rangle + 1} t^{\langle \al_i,\cla \rangle + 1} f_i,
i=1,\ldots,\ell$, to $v_{w^{-1} w_0(\nu)}$ in some order. The linear
combination of these monomials appearing in the solution is the term
of the lowest degree in $t$ with which $v_{w_0(\nu)}$ enters $\Phi(t)
v_{w^{-1} w_0(\nu)}$. It follows from \lemref{is generic} that it is
non-zero. The corresponding power of $t$ is nothing but the difference
between the $(\cla+\crho)$-degrees of the vectors $v_{w^{-1}w_0}$ and
$v_{w_0}$, i.e.,
$$
\langle w^{-1} w_0(\nu),\cla+\crho \rangle - \langle
w_0(\nu),\cla+\crho \rangle = \langle w_0(\nu),w(\cla+\crho) -
(\cla+\crho) \rangle,
$$
as desired. This completes the proof.
\end{proof}

Thus, we have identified the space
$\on{Conn}^{\on{reg}}_{D_x,\cla,w}$ of connections on the $H$-bundle
$\Omega^{\crho}$ on $D_x^\times$ of the form
with the space of Miura opers on $D_x^\times$ such that the underlying
oper belongs to $\on{Op}_G(D_x)_{\cla}$ and the corresponding
$B$-reductions $\F_B$ and $\F'_B$ have relative position $w$ at $x$.

The condition that the image under ${\mb b}_{\cla,w}$ of a
connection of the form
\begin{equation*}
\pa_t - \frac{w(\cla+\crho)-\crho}{t} + {\mb u}(t), \qquad {\mb u}(t)
\in \h[[t]],
\end{equation*}
is an oper without monodromy imposes polynomial equations on the
coefficients of the series ${\mb u}(t)$. Consider the simplest case
when $\cla=0$ and $w = s_i$, the $i$th simple reflection. Then
$-w(\cla+\crho)+\crho = \chal_i$, so we write this connection as
\begin{equation} \label{alpha i}
\ol{\nabla} = \pa_t + \frac{\chal_i}{t} + {\mb u}(t), \qquad {\mb
u}(t) \in \h[[t]].
\end{equation}

\begin{lem}[\cite{F:opers},Lemma 2.10] \label{si}
A connection of the form \eqref{alpha i} belongs to
$\on{Conn}^{\on{reg}}_{D_x,s_i}$ (i.e., the
corresponding $G$-oper is regular at $x$) if and only if $\langle
\al_i,{\mb u}(0) \rangle = 0$.
\end{lem}

\section{Wakimoto modules and Bethe Ansatz} \label{wak and ba}

In this section we explain how to construct Bethe eigenvectors of the
generalized Gaudin Hamiltonians. For that we utilize the Wakimoto
modules over $\ghat_{\ka_c}$ which are parameterized by Cartan
connections on the punctured disc. As the result, the eigenvectors
will be parameterized by the Cartan connections on $\pone$ with
regular singularities at $z_1,\ldots,z_N,\infty$ and some additional
points $w_1,\ldots,w_m$ with residues $\la_1,\ldots,\la_N,\la_\infty$
and $-\al_{i_1},\ldots,-\al_{i_m}$ and whose Miura transformation is
an oper that has no singularities at $w_1,\ldots,w_m$.

\subsection{Definition of Wakimoto modules} \label{def of wak}

We recall some of the results of \cite{FF:si,F:wak} on the
construction of the Wakimoto realization.

Let $\cA^{\g}$ be the Weyl algebra with generators $a_{\al,n},
a^*_{\al,n}$,  $\al \in \De_+$, $n \in \Z$, and relations
\begin{equation} \label{commina}
[a_{\al,n},a_{\beta,m}^*] = \delta_{\al,\beta}
\delta_{n,-m}, \hskip.3in
[a_{\al,n},a_{\beta,m}] = [a_{\al,n}^*,a_{\beta,m}^*] = 0.
\end{equation}
Introduce the generating functions
\begin{eqnarray} \label{az}
a_\al(z)&=&\sum_{n\in\Z} a_{\al,n} z^{-n-1}, \\ \label{a*z}
a_\al^*(z)&=&\sum_{n\in\Z} a_{\al,n}^* z^{-n}.
\end{eqnarray}

Let $M_{\g}$ be the Fock representation of $\cA^{\g}$ generated by a
vector $\vac$ such that
$$
a_{\al,n} \vac = 0, \quad n\geq 0; \qquad a^*_{\al,n} \vac = 0, \quad
n>0.
$$
It carries a vertex algebra structure (see \cite{F:wak}).

Let $\pi_0$ be the commutative algebra $\C[b_{i,n}]_{i=1,\ldots,\ell;
n<0}$ with the derivation $T$ given by the formula
$$
T \cdot b_{i_1,n_1} \ldots b_{i_m,n_m} = - \sum_{j=1}^m n_j
b_{i_1,n_1} \ldots b_{i_j,n_j-1} \cdots b_{i_m,n_m}.
$$
Then $\pi_0$ is naturally a commutative vertex algebra (see
\cite{FB}, \S~2.3.9). In particular, we have
$$
Y(b_{i,-1},z) = b_i(z) = \sum_{n<0} b_{i,n} z^{-n-1}.
$$

Recall that ${\mathbb V}_0$ carries the structure of a vertex
algebra. Its vacuum vector will be denoted by $v_{\ka_c}$.

In the next theorem describing the Wakimoto realization we will need
the formulas defining the action of the Lie algebra $\g$ on the space
of functions on the big cell $U$ of the flag manifold $G/B_+$. Here $U$
is the open orbit of the unipotent subgroup $N_+ = [B_+,B_+]$ which is
isomorphic to $N_+$. Since the exponential map $\n_+ \to N_+$ is an
isomorphism, we obtain a system of coordinates $\{ y_\al \}_{\al \in
\De_+}$ on $U$ corresponding to a fixed basis of root vectors $\{
e_\al \}_{\al \in \De_+}$ in $\n_+$.

For each $\chi \in \h^*$ we have a homomorphism $\rho_\chi: \g
\hookrightarrow \D_{\leq 1}(N_+)$ under which the $\D_{\leq
1}(N_+)$-module $\on{Fun} U$ becomes isomorphic to the $\g$-module
$M_\chi^*$ that is contragredient to the Verma module of highest
weight $\chi$.

Explicitly, $\rho_\chi$ looks as follows. Let $e_i, h_i, f_i,
i=1\ldots,\ell$, be the generators of~$\g$. Then
\begin{align} \label{formulas1}
\rho_\chi(e_i) &= \frac{\pa}{\pa y_{\al_i}} + \sum_{\beta \in \Delta_+}
P^i_\beta(y_\al) \frac{\pa}{\pa y_\beta}, \\ \label{formulas2}
\rho_\chi(h_i) &= - \sum_{\beta \in \Delta_+} \beta(h_i) y_\beta
\frac{\pa}{\pa y_\beta} + \chi(h_i), \\ \label{formulas3}
\rho_\chi(f_i) &= \sum_{\beta \in \Delta_+}
Q^i_\beta(y_\al) \frac{\pa}{\pa y_\beta} + \chi(h_i) y_{\al_i},
\end{align}
for some polynomials $P^i_\beta, Q^i_\beta$ in $y_\al, \al \in \De_+$.

In addition, we have a Lie algebra anti-homomorphism $\rho^R: \n_+ \to
\D_{\leq 1}(N_+)$ which corresponds to the {\em right} action of $\n_+$
on $N_+$. The differential operators $\rho^R(x), x \in \n_+$, commute
with the differential operators $\rho_\chi(x'), x' \in \n_+$ (though
their commutation relations with $\rho_\chi(x'), x' \not\in \n_+$, are
complicated in general). We have
$$
\rho^R(e_i) = \frac{\pa}{\pa y_{\al_i}} + \sum_{\beta \in \Delta_+}
P^{R,i}_\beta(y_\al) \frac{\pa}{\pa y_\beta}
$$
for some polynomials $P^{R,i}_\beta, Q^i_\beta$ in $y_\al, \al \in
\De_+$. We let
\begin{equation} \label{e right1}
e^R_i(z) = \sum_{n \in \Z} e^R_{i,n} z^{-n-1} = a_{\al_i}(z) +
\sum_{\beta \in \Delta_+} P^{R,i}_\beta(a^*_\al(z)) a_\beta(z).
\end{equation}

Now we can state the main result, due to \cite{FF:si,F:wak},
concerning the Wakimoto realization at the critical level.

\begin{thm} \label{exist hom}
There exists a homomorphism of vertex algebras $$\ww_{\ka_c}: {\mathbb
V}_0 \to M_\g \otimes \pi_0$$ such that
\begin{align*}
e_i(z) &\mapsto a_{\al_i}(z) + \sum_{\beta \in \Delta_+}
\Wick P^i_\beta(a^*_\al(z)) a_\beta(z) \Wick , \\
h_i(z) &\mapsto - \sum_{\beta \in \Delta_+} \beta(h_i) \Wick a^*_\beta(z)
a_\beta(z) \Wick + b_i(z) , \\
f_i(z) &\mapsto \sum_{\beta \in \Delta_+} \Wick
Q^i_\beta(a^*_\al(z)) a_\beta(z) \Wick + c_i \pa_z a^*_{\al_i}(z) +
b_i(z) a^*_{\al_i}(z),
\end{align*}
where the polynomials $P^i_\beta, Q^i_\beta$ are introduced in
formulas \eqref{formulas1}--\eqref{formulas3}.
\end{thm}

In order to make the homomorphism $w_{\ka_c}$ coordinate-independent,
we need to define the actions of the group $\on{Aut} \OO$ on ${\mathbb
V}_0, M_\g$ and $\pi_0$ which are intertwined by $w_{\ka_c}$. We
already have a natural action of $\AutO$ on ${\mathbb V}_0$ which is
induced by its action on $\ghat_{\ka_c}$ preserving the subalgebra
$\g[[t]]$. Next, we define the $\AutO$-action on $M_\g$ by
stipulating that for each $\al \in \De_+$ the generating functions
$a_\al(z)$ and $a^*_\al(z)$ transform as a one-form and a functions on
$D^\times$, respectively.

Finally, we identify $\pi_0$ with the algebra $\on{Fun} \Con_D$, where
$\Con_D$ is the space of connections on the $^L H$-bundle
$\Omega^{-\rho}$ on $D$. Namely, we write a connection from $\Con_D$
as the first order operator $\pa_t + \chi(t)$, where $\chi(t)$ takes
values in $^L \h \simeq \h^*$. Set $$b_i(t) = \langle \chal_i,\chi(t)
\rangle = \sum_{n<0} b_{i,n} t^{-n-1}.$$ Then we obtain an
identification
$$
\on{Fun} \Co_D \simeq \C[b_{i,n}]_{i=1,\ldots,\ell;n<0} = \pi_0.
$$
Now the natural action of $\AutO$ on $\on{Fun} \Co_D$ gives rise
to an action of $\AutO$ on $\pi_0$.

To write down this action explicitly, suppose that we have a
connection on the $^L H$-bundle $\Omega^{-\rho}$ on a curve $X$ which
with respect to a local coordinate $t$ is represented by the first
order operator $\pa_t + \chi(t)$. If $s$ is another coordinate such that
$t=\varphi(s)$, then this connection will be represented by the
operator
\begin{equation} \label{trans for conn1}
\pa_s + \varphi'(s) \chi(\varphi(s)) + \rho \cdot
\frac{\varphi''(s)}{\varphi'(s)}.
\end{equation}

\begin{prop}[\cite{F:wak}, Cor. 5.4]
The homomorphism $w_{\ka_c}$ commutes with the action of $\AutO$.
\end{prop}

For a connection $\ol\nabla$ on the $^L H$-bundle $\Omega^{-\rho}$,
let $\ol\nabla^*$ be the dual connection on
$\Omega^{\rho}$. Explicitly, if $\ol\nabla$ is given by the operator
$\pa_t + {\mb w}(t)$, then $\ol\nabla^*$ is given by the operator
$\pa_t - {\mb w}(t)$. Sending $\ol\nabla$ to $\ol\nabla^*$, we obtain
an isomorphism $\Con_U \simeq \Co_U$.

Recall the Miura transformation
$$
\ol{\mb b}_{D^\times}: \Co_{D^\times} \to \on{Op}_{^L G}(D^\times)
$$
defined in \secref{miura}. Using the isomorphism $\Con_U \simeq
\Co_U$, we obtain a map
\begin{equation} \label{star}
\ol{\mb b}_{D^\times}^*: \Con_{D^\times} \to \on{Op}_{^L G}(D^\times)
\end{equation}
which by abuse of terminology we will also refer to as Miura
transformation.

Consider the center $\zz(\ghat)$ of the vertex algebra ${\mathbb
V}_0$. Recall that we have
$$
\zz(\ghat) = {\mathbb V}_0^{\g[[t]]} \simeq
\on{End}_{\ghat_{\ka_c}} {\mathbb V}_0.
$$

\begin{thm}[\cite{F:wak}, Theorem 11.3] \label{center to conn}
Under the homomorphism $w_{\ka_c}$ the center $\zz(\ghat) \subset
{\mathbb V}_0$ gets mapped to the vertex subalgebra $\pi_0$ of $M_\g
\otimes \pi_0$. Moreover, we have the following commutative diagram
$$
\begin{CD}
\pi_0 @>{\sim}>> \on{Fun} \Con_D
\\ @AAA @AAA \\ \zz(\ghat) @>{\sim}>> \on{Fun}
\on{Op}_{^L G}(D)
\end{CD}.
$$
where the right vertical map is induced by the Miura transformation
$\ol{\mb b}_{D}^*$ given by~\eqref{star}.
\end{thm}

This theorem implies the following result. Let $N$ be a module over
the vertex algebra $M_\g$ and $R$ a module over the vertex algebra
$\pi_0 \simeq \on{Fun} \Co_D$. Then the homomorphism $w_{\ka_c}$
gives rise to the structure of a ${\mathbb V}_0$-module,
and hence of a smooth $\ghat_{\ka_c}$-module on the tensor product
$N \otimes R$.

It follows from the general theory of \cite{FB}, Ch. 5, that a module
over the vertex algebra $M_\g$ is the same as a smooth module over the
Weyl algebra $\cA^{\g}$, i.e., such that every vector is annihilated
by $a_{\al,n}, a^*_{\al,n}$ for large enough $n$. Likewise, a module
over the commutative vertex algebra $\on{Fun} \Con_D$ is the
same as a smooth module over the commutative topological algebra
$\on{Fun} \Con_{D^\times}$. Note that
$$
\on{Fun} \Con_{D^\times} \simeq \underset{\longleftarrow}\lim \;
\C[b_{i,n}]_{i=1,\ldots,\ell;n\in\Z}/I_N,
$$
where $I_N$ is the ideal generated by $b_{i,n}, i=1,\ldots,\ell; n\geq
N$. A module over this algebra is called smooth if every vector is
annihilated by the ideal $I_N$ for large enough $N$. Thus, for each
choice of $\cA^\g$-module $N$ we obtain a
``semi-infinite induction'' functor from the category of smooth
$\on{Fun} \Con_{D^\times}$-modules to the category of smooth
$\ghat_{\ka_c}$-modules (on which the central element ${\mb 1}$ acts
as the identity), $R \mapsto N \otimes R$.

Now \thmref{center to conn} implies

\begin{cor} \label{through miura}
The action of $Z(\ghat)$ on the $\ghat_{\ka_c}$-module $N \otimes R$
is independent of the choice of $N$ and factors through the
homomorphism $Z(\ghat) \simeq \on{Fun} \on{Op}_{^L G}(D^\times) \to
\on{Fun} \Con_{D^\times}$ induced by the Miura transformation
$\ol{\mb b}^*_{D^\times}$.
\end{cor}

For example, for a connection $\ol{\nabla} = \pa_z + \chi(z)$ in
$\Con_{D^\times}$, let $\C_{\ol{\nabla}} = \C_{\chi(z)}$ be the
corresponding one-dimensional module of $\on{Fun}
\Con_{D^\times}$. Then $M_\g \otimes \C_{\chi(z)}$ is called the {\em
Wakimoto module} of critical level corresponding to $\pa_z+\chi(z)$,
and is denoted by $W_{\chi(z)}$.

We will use below the following result.

\begin{lem}[\cite{FFR}, Lemma 2] \label{sing}
Let $$\chi(z) = - \frac{\al_i}{z} + \sum_{n=0}^\infty \chi_n z^n,
\quad \quad \chi_n \in \h^*.$$ The vector $e^R_{i,-1} \vac \in
W_{\chi(z)}$ is $\g[[t]]$-invariant if and only if $\langle
\chal_i,\chi_0 \rangle = 0$.
\end{lem}

\subsection{Functoriality of coinvariants}    \label{funct of coinv}

We will construct eigenvectors of the generalized Gaudin Hamiltonians
using functionals on the space of coinvariants of Wakimoto modules on
$\pone$. First we explain some general facts about the compatibility
of various spaces of coinvariants.

\looseness=1 Following the general construction of \cite{FB}, Ch. 10,
we define the spaces\break of coinvariants
$\,H_{M_\g}(X,(x_i),(M_i))\,$ (resp., $\,H_{\pi_0}(X,\{ x_i
\},(R_i))\,$) attached to a\break smooth projective curve $X$, points
$x_1,\ldots,x_p$ and $M_\g$-modules $N_1,\ldots,N_p$ (resp.,
$\pi_0$-modules $R_1,\ldots,R_p$). These spaces are defined in
elementary terms in \cite{FFR}, Sect. 5, and in \cite{FB},
Ch. 13--14. We also define spaces of coinvariants for the tensor
product vertex algebra $M_\g \otimes \pi_0$. We have a natural
identification
$$
H_{M_\g \otimes \pi_0}(X,(x_i),(N_i \otimes R_i)) =
H_{M_\g}(X,(x_i),(N_i)) \otimes H_{\pi_0}(X,(x_i),(R_i)).
$$

The vertex algebra homomorphism $\ww_{\ka_c}: {\mathbb V}_0 \to
M_\g \otimes \pi_0$ of \thmref{exist hom} gives rise to a surjective
map of the spaces of coinvariants
$$
H_{{\mathbb V}_0}(X,(x_i),(N_i \otimes R_i))
\twoheadrightarrow H_{M_\g \otimes \pi_0}(X,(x_i),(N_i \otimes
R_i))
$$
and hence a surjective map
\begin{equation} \label{surj map}
H_{{\mathbb V}_0}(X,(x_i),(N_i \otimes R_i))
\twoheadrightarrow H_{M_\g}(X,(x_i),(N_i)) \otimes
H_{\pi_0}(X,(x_i),(R_i)).
\end{equation}

Now recall that in the proof of \thmref{descr of alg} we explained
that since $\zz(\ghat) \simeq \on{Fun} \on{Op}_{^L G}(D)$ is the
center of $\V_0$, we have a homomorphism of algebras
\begin{equation} \label{point u1}
H_{\zz(\G)}(X,(x_i),(Z(M_i))) \to \on{End}_\C H_{\V_0}(X,(x_i),(M_i)).
\end{equation}
In addition, we have an isomorphism
\begin{equation} \label{opm}
H_{\zz(\G)}(X,(x_i),(Z(M_i))) \simeq \on{Fun}
\on{Op}_{^L G}(X,(x_i),(M_i)),
\end{equation}
where $\on{Op}_{^L G}(X,(x_i),(M_i))$ is the space of $^L G$-opers
on $X$ which are regular\break on $X\bs \{ x_1,\ldots,x_N \}$ and such that
their restriction to $D_x^\times$ belongs to the space $\on{Op}_{^L
G}^{M_i}(D_{x_i}^\times) = \on{Spec} Z(M_i)$ for all
$i=1,\ldots,N$.

We can describe the space of coinvariants $H_{\pi_0}(X,(x_i),(R_i))$
for the commutative vertex algebra $\pi_0 \simeq \on{Fun} \Con_D$ in
similar terms. If $R$ is a $\pi_0$-module, or equivalently, a smooth
$\on{Fun} \Con_{D^\times}$-module, let $\pi_0(R)$ be the
image of the corresponding homomorphism
$$
\on{Fun} \Con_{D^\times} \to \on{End}_{\C} R.
$$

Let $\Con(X,(x_i)$, $(R_i))$ be the space of connections on
$\Omega^{-\rho}$ on\break $X \bs \{ x_1,\ldots,x_N \}$ whose restriction to
$D_{x_i}^\times$ belongs to the space $\Con_{D_{x_i}^\times}^{R_i} =
\on{Spec} \pi_0(R_i)$. Then we have an isomorphism
$$
H_{\pi_0}(X,(x_i),(R_i)) \simeq \on{Fun} \Con(X,(x_i),(R_i)).
$$

It follows from \corref{through miura} that if $M_i = N_i \otimes
R_i$, then $Z(N_i \otimes R_i)$ is the image of $Z(\ghat) \simeq
\on{Fun}_{^L G}(D^\times)$ under the homomorphism
$$
\on{Fun}_{^L G}(D^\times) \overset{\ol{\mb
b}^*_{D^\times}}\longrightarrow \on{Fun} \Con_{D^\times} \to
\on{End}_{\C} R_i,
$$
where the first map is induced by the Miura transformation.

Recall that in \secref{cartan conn} we introduced the maps
$$
{\mb b}_U: \Co_U \to \on{MOp}_{^L G}(U)_{\on{gen}}, \qquad \ol{\mb b}_U:
\Co_U \to \on{Op}_{^L G}(U).
$$
Using the identification $\Co_U \simeq \Con_U$, we obtain maps
$$
{\mb b}^*_U: \Con_U \to \on{MOp}_{^L G}(U)_{\on{gen}}, \qquad \ol{\mb
b}^*_U: \Con_U \to \on{Op}_{^L G}(U).
$$
The second map is the Miura transformation. In the case when
$U = X
\bs \{ x_1,\ldots,x_N \}$ the Miura transformation $\ol{\mb
b}^*_U$ restricts to a map
$$
\Con(X,(x_i),(R_i)) \to \on{Op}_{^L G}(X,(x_i),(N_i
\otimes R_i)),
$$
and hence gives rise to a homomorphism
\begin{equation} \label{one more miura}
\on{Fun} \on{Op}_{^L G}(X,(x_i),(N_i
\otimes R_i)) \to \on{Fun} \Con(X,(x_i),(R_i)),
\end{equation}
which does not depend on the $N_i$'s.

Now we have an action of the algebra $\on{Op}_{^L G}(X,(x_i),(N_i
\otimes R_i))$ on both sides of the map \eqref{surj map}. The action
on the left-hand side comes from the homomorphisms \eqref{point u1}
and \eqref{opm}. The action on the right-hand side comes from the
homomorphism \eqref{one more miura}.

The functoriality of the coinvariants implies the following

\begin{lem} \label{compatib}
The map \eqref{surj map} commutes with the action on both sides of the
algebra $\on{Op}_{^L G}(X,(x_i),(N_i \otimes R_i))$.
\end{lem}

Consider the special case where all modules $R_i, i=1,\ldots,N$, are
one-dimen\-sion\-al, i.e., $R_i = \C_{\ol\nabla_i}$, where $\ol\nabla_i$
is a point in $\Con_{D_{x_i}^\times}$, and $R_i$ is the representation
obtained from the homomorphism $\on{Fun} \Con_{D_{x_i}^\times}$
induced by $\ol\nabla_i$. In this case we have the following
description of the space $H_{\pi_0}(X,(x_i),(\C_{\ol\nabla_i}))$ which
follows from the general results on coinvariants of commutative vertex
algebras from Sect. 9.4 of \cite{FB} (see also \cite{FFR}, Prop. 4, in
the case when $X=\pone$).

\begin{prop} \label{onedim}
The space $H_{\pi_0}(X,(x_i),(\C_{\ol\nabla_i}))$ is one-dimensional
if and only if there exists a connection $\ol{\nabla}$ on the $^L
H$-bundle $\Omega^{-\rho}$ on $X \bs \{ x_1,\ldots,x_p,\infty \}$
such that the restriction of $\ol\nabla$ to the punctured disc at each
$x_i$ is equal to $\ol\nabla_i, i=1,\ldots,N$.

Otherwise, $H_{\pi_0}(X,(x_i),(\C_{\ol\nabla_i}))=0$.
\end{prop}

Suppose that such a connection $\ol{\nabla}$ exists and the space
$H_{\pi_0}(X,(x_i),(\C_{\ol\nabla_i}))$ is one-dimensional. Then the
map \eqref{surj map} reads
$$
H_{{\mathbb V}_0}(X,(x_i),(N_i \otimes \C_{\ol\nabla_i}))
\twoheadrightarrow H_{M_\g}(X,(x_i),(N_i)).
$$
Composing this map with any linear functional on
$H_{M_\g}(X,(x_i),(N_i))$, we obtain a linear functional
$$
\phi: H_{{\mathbb V}_0}(X,(x_i),(N_i \otimes \C_{\ol\nabla_i})) \to
\C.
$$

The algebra $\zz(\ghat)_u = \on{Fun} \on{Op}_{^L G}(D_u)$ acts on the
space $H_{{\mathbb V}_0}(X,(x_i),(N_i \otimes \C_{\ol\nabla_i}))$, and
its action factors through the homomorphism \eqref{from u}
$$
\on{Fun} \on{Op}_{^L G}(D_u) \to \on{Fun} \on{Op}_{^L G}(X,(x_i),(N_i
\otimes \C_{\ol\nabla_i})).
$$
According to \lemref{compatib}, this action factors through the
homomorphism \eqref{one more miura}. In our case $R_i =
\C_{\ol\nabla_i}$, and so $\Con(X,(x_i),(R_i))$ is a point
corresponding to the connection $\ol\nabla$. It follows from the
definition that $\on{Op}_{^L G}(X,(x_i),(N_i \otimes
\C_{\ol\nabla_i}))$ is also a point, corresponding to the oper which
is the image of $\ol\nabla$ under the Miura transformation $\ol{\mb
b}^*_{X \bs \{ x_1,\ldots,x_N \}}$. Thus, we obtain:

\begin{prop} \label{eigenv}
The functional $\phi$ is an eigenvector of the algebra $\zz(\G)_u
\simeq\break \on{Op}_{^L G}(D_u)$, and the corresponding homomorphism
$\on{Fun} \on{Op}_{^L G}(D_u) \to \C$ is defined by the point
$\ol{\mb b}^*_{D_u}(\ol\nabla|_{D_u}) \in \on{Op}_{^L G}(D_u)$.
\end{prop}

\subsection{The case of $\pone$}    \label{the case of pone}

We will use \propref{eigenv} in the case when $X = \pone$, equipped
with a global coordinate $t$. Consider a collection of distinct
points $x_1,\ldots,x_p$ and $\infty$ on $\pone$. As the $M_\g$-module
$N_i$ attached to $x_i, i=1,\ldots,p$, we will take $M_\g$. As the
$M_\g$-module $N_\infty$ attached to the point $\infty$ we will take
another module $M'_\g$ generated by a vector $\vac'$ such that
$$
a_{\al,n} \vac' = 0, \quad n > 0; \qquad a^*_{\al,n} \vac' = 0, \quad
n \geq 0.
$$
As the $\pi_0$-module $R_i$ attached to $x_i, i=1,\ldots,p$, we will
take the one-dimensional module $\C_{\ol\nabla_i} = \C_{\chi_i(z)}$,
where $\ol\nabla_i = \pa_z + \chi_i(z)$, $\chi_i(z) \in \h^*((z))$ and
as module $R_\infty$ attached to the point $\infty$ we will take
$\C_{\chi_\infty(z)}$. Thus, $M_i \otimes R_i$ is the Wakimoto module
$W_{\chi_i(z)}$ for all $i=1,\ldots,p$. We denote the module $M'_\g
\otimes \C_{\chi_\infty(z)}$ by $W'_{\chi_\infty(z)}$.

We have the following special case of \propref{onedim}:

\begin{prop}[\cite{FFR}, Prop. 4] \label{one}
The space $H_{M_\g}(\pone;(x_i),\infty;(M_\g),M'_\g)$ is
one-dimensional and the projection of the vector $\vac^{\otimes N}
\otimes \vac'$ on it is non-zero.

The space
$H_{\pi_0}(\pone,(x_i),\infty;(\C_{\chi_i(z)}),\C_{\chi_\infty}(z))$ is
one-dimensional if and only if there exists a connection $\ol{\nabla}$
on the $^L H$-bundle $\Omega^{-\rho}$ on $\pone$ which is regular on
$\pone \bs \{ x_1,\ldots,x_p,\infty \}$ and such that the
restriction of $\ol\nabla$ to the punctured disc at each $x_i$ is
equal to $\pa_t + \chi_i(t-x_i)$, and its restriction to the punctured
disc at $\infty$ is equal to $\pa_{t^{-1}} + \chi_\infty(t^{-1})$.

Otherwise,
$H_{\pi_0}(\pone,(x_i),\infty;(\C_{\chi_i(z)}),\C_{\chi_\infty}(z))=0$.
\end{prop}

Formula \eqref{trans for conn1} shows that if we have a connection on
$\Omega^{-\rho}$ over $\pone$ whose restriction to $\pone \bs \infty$
is represented by the operator $\pa_t + \chi(t)$, then its restriction
to the punctured disc $D_\infty^\times$ at $\infty$ reads, with
respect to the coordinate $u = t^{-1}$, is represented by the operator
$$
\pa_u - u^{-2} \chi(u^{-1}) - 2\rho u^{-1}.
$$

We will choose $\chi_i(z)$ to be of the form $$\chi_i(z) = \chi_i +
\sum_{n\geq 0} \chi_{i,n} z^n,$$ and $\chi_\infty(z)$ to be of the form
$$\chi_\infty(z) = \chi_\infty + \sum_{n\geq 0} \chi_{\infty,n} z^n.$$
The condition of the proposition is then equivalent to saying that the
restriction of $\ol\nabla$ to $\pone \bs \infty$ is represented by the
operator
$$
\pa_t + \sum_{i=1}^p \frac{\chi_i}{t-z_i},
$$
and that $\chi_i(t-z_i)$ is the expansion of $\chi(t)$ at $z_i,
i=1,\ldots,p$, while $\chi_\infty(t^{-1})$ is the expansion of $- t^2
\chi(t) - 2\rho t$ in powers of $t^{-1}$. Note that then we have
$\chi_\infty(u) = \chi_\infty/u + \ldots$, where $\chi_\infty = -2\rho
- \sum_{i=1}^p \chi_i$.

By \propref{one}, the spaces of coinvariants
$H_{M_\g}(\pone;(x_i),\infty;(M_\g),M'_\g)$ and
$H_{\pi_0}(\pone;(x_i),\infty;(\C_{\chi_i(z)}),\C_{\chi_\infty}(z))$
are one-dimensional. Hence the map \eqref{surj map} gives rise to a
non-zero linear functional
\begin{equation} \label{non zero f}
\tau_{(x_i)}:
H_{\V_0}(\pone;(x_i),\infty;(W_{\chi_i(z)}),W'_{\chi_\infty(z)}) \to \C.
\end{equation}
We normalize it so that its value on $\vac^{\otimes N} \otimes \vac'$
is equal to $1$.

The algebra $\zz(\ghat)_u = \on{Fun} \on{Op}_{^L G}(D_u)$ acts on this
space. By \propref{eigenv}, $\tau_{(x_i)}$ is an eigenvector of this
algebra and the corresponding homomorphism $$\on{Fun} \on{Op}_{^L
G}(D_u) \to \C$$ is defined by the point $\ol{\mb
b}^*_{D_u}(\ol\nabla|_{D_u}) \in \on{Op}_{^L G}(D_u)$.

Let us fix highest weights of $\g$, $\la_1,\ldots,\la_N$, and a set of
simple roots of $\g$, $\al_{i_1},\ldots,\al_{i_m}$. Consider a
connection $\ol{\nabla}$ on $\Omega^{-\rho}$ on $\pone$ whose
restriction to $\pone \bs \infty$ is equal to $\pa_t + \la(t)$, where
\begin{equation} \label{sv}
\la(t) = \sum_{i=1}^N \frac{\la_i}{t-z_i} - \sum_{j=1}^m
\frac{\al_{i_j}}{t-w_j}.
\end{equation}
Denote by $\la_i(t-z_i)$ the expansions of $\la(t)$ at the points $z_i,
i=1,\ldots,N$, and by $\mu_j(t-w_j)$ the expansions of $\la(t)$ at the
points $w_j, j=1,\ldots,m$. We have: $$\la_i(z) = \frac{\la_i}{z} + \cdots,
\quad \quad \mu_j(z) = - \frac{\al_{i_j}}{z} + \mu_{j,0} + \cdots,$$ where
\begin{equation}    \label{muj0}
\mu_{j,0} = \sum_{i=1}^N \frac{\la_i}{w_j-z_i} - \sum_{s\neq j}
\frac{\al_{i_s}}{w_j-w_s}.
\end{equation}
Finally, let $\la_\infty(t^{-1})$ be the expansion of $- t^2
\la(t) - 2\rho t$ in powers of $t^{-1}$. Thus we have
\begin{equation} \label{la infty}
\la_\infty(z) = z^{-1} \left( - \sum_{i=1}^N \la_i + \sum_{j=1}^m
\al_{i_j} - 2 \rho \right) + \cdots.
\end{equation}
We have a linear functional $\tau_{(z_i),(w_j)}$ on the corresponding
space of coinvariants
$$
H_{{\mathbb
V}_0}(\pone;(z_i),(w_j),\infty;(W_{\la_i(z)}),(W_{\mu_j(z)}),
W'_{\la_\infty(z)}).
$$
In particular, this space is non-zero. The functional
$\tau_{(z_i),(w_j)}$ is an eigenvector of $\zz(\ghat)_u = \on{Fun}
\on{Op}_{^L G}(D_u)$, and its eigenvalue is the $^L G$-oper on $D_u$
which is the restriction to $D_u$ of the oper in $\on{Op}^{\on{RS}}_{^L
G}(\pone)_{(z_i),(w_j),\infty;(\la_i),(-\al_{i_j}),\la_\infty}$
given by the Miura transformation of the connection $\pa_t + \la(t)$. We
will now use it to construct an eigenvector of the Gaudin algebra.

\subsection{Bethe vectors} \label{bethe vect} \label{miur}

By \lemref{sing}, the vectors $e^R_{i_j,-1} \vac \in W_{\mu_j(z)}$ are
$\g[[t]]$-invariant if and only if the equations $\langle
\chal_{i_j},\mu_{j,0} \rangle = 0$ are satisfied, where $\mu_{j,0}$ is
the constant coefficient in the expansion of $\la(t)$ at $w_j$ given
by formula \eqref{muj0}, i.e., the following system of equations
is satisfied
\begin{equation} \label{bethe}
\sum_{i=1}^N \frac{\langle \chal_{i_j},\la_i \rangle}{w_j-z_i} -
\sum_{s\neq j} \frac{\langle \chal_{i_j},\al_{i_s} \rangle}{w_j-w_s} =
0, \quad j=1,\ldots,m.
\end{equation}
They are called the {\em Bethe Ansatz equations}. This is a system of
equations on the complex numbers $w_j, j=1,\ldots,m$, to each of which
we attach a simple root $\al_{i_j}$.  We have an obvious action of a
product of symmetric groups permuting the points $w_j$ corresponding
to simple roots of the same kind. In what follows, by a {\em solution}
of the Bethe Ansatz equations we will understand a solution defined up
to these permutations.  We will adjoin to the
set of all solutions associated to all possible collections $\{
\al_{i_j} \}$ of simple roots of $\g$, the unique ``empty'' solution,
corresponding to the empty set of simple roots.

Suppose that these equations are satisfied. Then we obtain a
homomorphism of $\G_{\ka_c}$-modules
$$
\bigotimes_{i=1}^N W_{\la_i(z)} \otimes \V_0^{\otimes m} \otimes
W'_{\la_\infty(z)} \to \bigotimes_{i=1}^N
W_{\la_i(z)} \otimes \bigotimes_{j=1}^m W_{\mu_j(z)} \otimes
W'_{\la_\infty(z)},
$$
which sends the vacuum vector in the $j$th copy of $\V_0$ to
$e^R_{i_j,-1} \vac \in W_{\mu_j(z)}$. Hence we obtain the
corresponding map of the spaces of coinvariants. But the insertion of
$\V_0$ does not change the space of coinvariants. Hence the first
space of coinvariants is isomorphic to
\begin{equation*} \label{coinvar}
H_{{\mathbb
V}_0}(\pone;(z_i),\infty;(W_{\la_i(z)}),W'_{\la_\infty(z)}) =
H((W_{\la_i(z)}),W'_{\la_\infty(z)}).
\end{equation*}

Therefore, composing this map with $\tau_{(z_i),(w_j)}$, we obtain a
linear functional
\begin{equation}    \label{linfunct}
\wt\tau_{(z_i),(w_j)}: H((W_{\la_i(z)}),W'_{\la_\infty(z)}) \to \C.
\end{equation}
By construction, $\wt\tau_{(z_i),(w_j)}$ is an eigenvector of
$\zz(\ghat)_u = \on{Fun} \on{Op}_{^L G}(D_u)$, and its eigenvalue is
the $^L G$-oper on $D_u$ which is the restriction of an oper on
$\pone$ with singularities given by the Miura transformation of the
connection $\pa_t + \la(t)$.

Let $U = \pone \bs \{z_1,\ldots,z_N,w_1,\ldots,w_m,\infty \}$ and
$\Con^{\on{gen}}_{(z_i),\infty;(\la_i),\la_\infty}$ be the subset of
$\Con_U$, consisting of all connections on $\Omega^{-\rho}$ on $\pone$
with regular singularities of the form $\pa_t + \la(t)$, where
$\la(t)$ is given by \eqref{sv}, such that their Miura transformation
is a $^L G$-oper in
\begin{equation}    \label{into}
\on{Op}_{^L G}^{\on{RS}}(\pone)_{(z_i),\infty;(\la_i),\la_\infty}
\subset \on{Op}_{^L
G}^{\on{RS}}(\pone)_{(z_i),(w_j),\infty;(\la_i),(-\al_{i_j}),\la_\infty}.
\end{equation}
In other words, $\Con^{\on{gen}}_{(z_i),\infty;(\la_i),\la_\infty}$
consists of those connections $\pa_t + \la(t)$, where $\la(t)$ is
given by \eqref{sv}, whose Miura transformation is a $^L G$-oper on
$U$ whose restriction to $D_{w_j}^\times$ belongs to $\on{Op}_{^L
  G}(D_{w_j}) \subset \on{Op}_{^L G}^{\on{RS}}(D_{w_j}^\times)$.

Thus, we have a map
\begin{equation} \label{miura on pone}
\ol{\mb b}^*_{(z_i),\infty;(\la_i),\la_\infty}:
\Con^{\on{gen}}_{(z_i),\infty;(\la_i),\la_\infty} \to \on{Op}_{^L
G}^{\on{RS}}(\pone)_{(z_i),\infty;(\la_i),\la_\infty}.
\end{equation}

\begin{prop}[\cite{FFR}, Prop. 6] \label{bij1}
There is a bijection between the set of solutions of the Bethe Ansatz
equations \eqref{bethe} and the set
$\Con^{\on{gen}}_{(z_i),\infty;(\la_i),\la_\infty}$.
\end{prop}

\begin{proof}
Applying the Miura transformation to the connection $\pa_t + \la(t)$,
where $\la(t)$ is given by \eqref{sv}, we obtain a $^L G$-oper on
$\pone$ which has regular singularities at the points
$z_1,\ldots,z_N$, and $\infty$ with residues
$\varpi(\la_1+\rho),\ldots,\varpi(\la_N+\rho)$ and
$\varpi(\la_\infty+\rho)$, respectively, and at the points
$w_1,\ldots,w_m$ with residues $\varpi(-\al_{i_j}+\rho) =
\varpi(\rho)$. But by \lemref{si}, this oper is regular at the points
$w_1,\ldots,w_m$ if and only if the equations $\langle
\chal_{i_j},\mu_{j,0} \rangle = 0$ are satisfied, where $\mu_{j,0}$ is
the constant coefficient in the expansion of $\la(t)$ at $w_j$ given
by formula \eqref{muj0}. These equations are precisely the Bethe
Ansatz equations \eqref{bethe}.
\end{proof}

Now let $w_1,\ldots,w_m$ be a solution of the Bethe Ansatz equations
\eqref{bethe}. Consider the corresponding connection $\pa_t + \la(t)$
in $\Con^{\on{gen}}_{(z_i),\infty;(\la_i),\la_\infty}$, where $\la(t)$
is of the form \eqref{sv}. Recall that the linear functional
$\wt\tau_{(z_i),(w_j)}$ from formula \eqref{linfunct} is an
eigenvector of $\zz(\ghat)_u = \on{Fun} \on{Op}_{^L G}(D_u)$, and its
eigenvalue is the $^L G$-oper on $D_u$ obtained by restriction of the
oper in $\on{Op}_{^L
G}^{\on{RS}}(\pone)_{(z_i),\infty;(\la_i),\la_\infty} $ given by the
Miura transformation of the connection $\pa_t + \la(t)$.

Therefore we obtain that the action of $\zz(\ghat)_u$ on the space
$H((W_{\la_i(z)}),W'_{\la_\infty(z)})$ factors through the Gaudin
algebra $\ZZ(\g)$, which acts on it according to the character
corresponding to the $^L G$-oper $\ol{\mb
b}^*_{(z_i),\infty;(\la_i),\la_\infty}(\pa_t+\la(t))$.

We wish to use this fact to construct eigenvectors of the algebra
$\ZZ(\g)$ on the tensor product of Verma modules $\bigotimes_{i=1}^N
M_{\la_i}$.

For $\la(z) \in \h^* \otimes z^{-1}\C[[z]]$ denote by $\la_{-1}$
the residue of $\la(z) dz$ at $z=0$. Let $\wt{W}_{\la(z)}$ the
subspace of a Wakimoto module $W_{\la(z)}$, which is generated
from the vector $\vac$ by the operators $a_{\al,0}^*, \al \in
\De_+$. This space is stable under the action of the constant
subalgebra $\g$ of $\G_{\ka_c}$. It follows from the formulas
defining the action of $\ghat_{\ka_c}$ on $W_{\la(z)}$ that as a
$\g$-module, $\wt{W}_{\la(z)}$ is isomorphic to the module
$M^*_{\la_{-1}}$ contragredient to the Verma module over $\g$ with
highest weight $\la_{-1}$. Likewise, the subspace
$\wt{W}'_{\la(z)}$ generated from the vector $\vac'$ by the
operators $a_{\al,0}, \al \in \De_+$, is stable under the action
of the constant subalgebra $\g$ of $\G_{\ka_c}$ and is isomorphic
as a $\g$-module to the module $M'_{2\rho+\la_{-1}}$ the Verma
module over $\g$ with {\em lowest} weight $2\rho+\la_{-1}$.

Let $\M^*_\chi$ be the $\ghat_{\ka_c}$-module induced from
$M^*_\chi$, and let $\M'_\chi$ be the $\ghat_{\ka_c}$-module induced
from $M'_\chi$ (see \secref{coinvariants}). The embeddings of
$\g$-modules $M^*_{\la_{-1}} \to W_{\la(z)}$ and $M'_{2\rho+\la_{-1}}
\to W'_{\la(z)}$ induce homomorphisms of $\ghat_{\ka_c}$-modules
$\M^*_{\la_{-1}} \to W_{\la(z)}$ and $\M'_{2\rho+\la_{-1}} \to
W'_{\la(z)}$. Thus, we obtain a homomorphism
$$
\bigotimes_{i=1}^N \M^*_{\la_i} \otimes \M'_{2\rho+\la_\infty} \to
\bigotimes_{i=1}^N W_{\la_i(z)} \otimes W'_{\la_\infty(z)},
$$
where
$$
\la_\infty = - 2 \rho - \sum_{i=1}^N \la_i + \sum_{j=1}^m \al_{i_j}.
$$
This homomorphism gives rise to a map of the corresponding spaces of
coinvariants
$$
H((\M^*_{\la_i}),\M'_{2\rho+\la_\infty}) \to
H((W_{\la_i(z)}),W'_{\la_\infty(z)}).
$$

The composition of this map with the functional
$\wt{\tau}_{(z_i),(w_j)}$ defined above is a linear functional
$$
\ol{\tau}_{(z_i),(w_j)}: H((\M^*_{\la_i}),\M'_{2\rho+\la_\infty}) \to
\C.
$$
But we have an isomorphism
$$
H((\M^*_{\la_i}),\M'_{2\rho+\la_\infty}) \simeq (\bigotimes_{i=1}^N
M^*_{\la_i} \otimes M'_{2\rho+\la_\infty})/\g \simeq
(\bigotimes_{i=1}^N M^*_{\la_i})_{-2\rho-\la_\infty}/\n_-,
$$
where the subscript indicates the subspace of
weight $-2\rho-\la_\infty = \sum_{i=1}^N \la_i - \sum_{j=1}^m
\al_{i_j}$. Thus, we obtain an $\n_-$-invariant functional
$$
\psi(w_1^{i_1},\ldots,w_m^{i_m}): ( \bigotimes_{i=1}^N
M^*_{\la_i})_{-2\rho-\la_\infty} \arr \C,
$$
or equivalently an $\n_+$-invariant vector
$$
\phi(w_1^{i_1},\ldots,w_m^{i_m}) \in \otimes_{i=1}^N
M_{\la_i}
$$
of weight $\sum_{i=1}^N \la_i - \sum_{j=1}^m \al_{i_j}$.

Recall that $\wt{\tau}_{(z_i),(w_j)}$ is an eigenvector of the
algebra $\zz(\ghat)_u = \on{Fun} \on{Op}_{^L G}(D_u)$, and its
eigenvalue is the $^L G$-oper on $\pone$ given by the Miura
transformation of the connection $\pa_t + \la(t)$, where $\la(t)$
is given by \eqref{sv}. Hence $\phi(w_1^{i_1},\ldots,w_m^{i_m})$
is also an eigenvector of $\zz(\ghat)_u$ with the same eigenvalue.
But we know from the proof of \thmref{descr of alg} that the
action of $\zz(\ghat)_u$ on the space of coinvariants
$H((\M^*_{\la_i}),\M'_{2\rho+\la_\infty})$ factors through the
Gaudin algebra $\ZZ(\g)$ and further through ${\mc
Z}_{(z_i),\infty;(\la_i),\la_\infty}(\g)$. Therefore
$\phi(w_1^{i_1},\ldots,w_m^{i_m})$ is an eigenvector of the Gaudin
algebra $\ZZ(\g)$, and the eigenvalues of $\ZZ(\g)$ on this vector
are encoded by the $^L G$-oper obtained by the Miura
transformation of the connection $\pa_t + \la(t)$, where $\la(t)$
is given by \eqref{sv}.

Let us summarize our results.

\begin{thm}[\cite{FFR}, Theorem 3] \label{main}
If the Bethe Ansatz equations \eqref{bethe} are satisfied, then the
vector
$$
\phi(w_1^{i_1},\ldots,w_m^{i_m}) \in \otimes_{i=1}^N
M_{\la_i}
$$
is an eigenvector of the algebra
${\mc Z}_{(z_i),\infty;(\la_i),\la_\infty}(\g)$. Its eigenvalues
define a point in $$\on{Spec} {\mc
Z}_{(z_i),\infty;(\la_i),\la_\infty}= \on{Op}_{^L
G}^{\on{RS}}(\pone)_{(z_i),\infty;(\la_i),\la_\infty},$$ which is the
image of the connection $\pa_t + \la(t)$ with $\la(t)$ as in formula
\eqref{sv}, under the Miura transformation $\ol{\mb
b}^*_{(z_i),\infty;(\la_i),\la_\infty}$ given by formula \eqref{miura
on pone}.
\end{thm}

The vector $\phi(w_1^{i_1},\ldots,w_m^{i_m})$ is called the {\em Bethe
vector} corresponding to a solution of the Bethe Ansatz equations
\eqref{bethe}, or equivalently, an element of the set
$\Con^{\on{gen}}_{(z_i),\infty;(\la_i),\la_\infty}$. An explicit formula for this
vector is given in \cite{FFR}, Lemma 3 (see also \cite{ATY}):
\begin{align} \label{genbv}
&\phi(w_1^{i_1},\ldots,w_m^{i_m})\\
&=  (-1)^m \sum_{p=(I^1,\ldots,I^N)}
\prod_{j=1}^N \frac{F_{i^j_1}^{(j)} F_{i^j_2}^{(j)} \ldots
F_{i^j_{a_j}}^{(j)}}{(w_{i^j_1}-w_{i^j_2})(w_{i^j_2}-w_{i^j_3}) \ldots
(w_{i^j_{a_j}}-z_j)} v_{\la_1} \otimes \ldots v_{\la_N}.\nonumber
\end{align}
Here the summation is taken over all {\em ordered} partitions $I^1 \cup
I^2 \cup \ldots \cup I^N$ of the set $\{1,\ldots,m\}$, where $I^j = \{
i^j_1,i^j_2,\ldots,i^j_{a_j} \}$.

The fact that this vector is an eigenvector of the Gaudin
Hamiltonians $\Xi_i$ has also been established by other methods in
\cite{bf,RV}. But \thmref{main} gives us much more: it tells us
that this vector is also an eigenvector of all other generalized
Gaudin Hamiltonians, and identifies the $^L G$-oper encoding
their eigenvalues on this vector as the Miura transformation of
the connection $\pa_t + \la(t)$. We recall that this last result
follows from \corref{through miura} which states that the central
characters on the Wakimoto modules are obtained via the Miura
transformation.

\section{The case of finite-dimensional $\g$-modules}
\label{findim mod}

In this section we consider the eigenvectors of the Gaudin algebra
$\ZZ(\g)$ on the space
$$
H((V_{\la_i}),V_{\la_\infty}) = \left( \bigotimes_{i=1}^N V_{\la_i}
\otimes V_{\la_\infty} \right)^G
$$
(recall that $V_\la$ denotes the irreducible finite-dimensional
$\g$-module with a dominant integral highest weight $\la$). Recall
that the action of $\ZZ(\g)$ on this space factors through the algebra
${\mc Z}_{(z_i),\infty;(\la_i),\la_\infty}(\g)$, so that we have a
homomorphism
$$
{\mc Z}_{(z_i),\infty;(\la_i),\la_\infty} \to \on{End}_\C \;
\bigotimes_{i=1}^N V_{\la_i} \otimes V_{\la_\infty}.
$$
The image of this homomorphism was denoted by $\ol{\mc
Z}_{(z_i),\infty;(\la_i),\la_\infty}(\g)$. Then the set of all common
eigenvalues of the algebra
$\ZZ(\g)$ acting on $\left( \bigotimes_{i=1}^N V_{\la_i}
\otimes V_{\la_\infty} \right)^G$ (without multiplicities)
is precisely the spectrum of the algebra $\ol{\mc
Z}_{(z_i),\infty;(\la_i),\la_\infty}(\g)$.

By \thmref{descr of alg},(3), the algebra $\ol{\mc
Z}_{(z_i),\infty;(\la_i),\la_\infty}(\g)$ is a quotient of the
algebra of functions on the space $\on{Op}_{^L
G}(\pone)_{(z_i),\infty;(\la_i),\la_\infty}$, so that its spectrum
injects into $\on{Op}_{^L
G}(\pone)_{(z_i),\infty;(\la_i),\la_\infty}$. Thus, we have there
following result:

\begin{prop}    \label{inj map1}
There is an injective map
\begin{equation} \label{inj map}
\on{Spec} \ol{\mc Z}_{(z_i),\infty;(\la_i),\la_\infty}(\g) \hookrightarrow
\on{Op}_{^L G}(\pone)_{(z_i),\infty;(\la_i),\la_\infty}.
\end{equation}
In other words, each common eigenvalue of the Gaudin algebra
$\ZZ(\g)$ that occurs in the space $\left( \bigotimes_{i=1}^N
V_{\la_i} \otimes V_{\la_\infty} \right)^G$ is encoded by a point
of $\on{Op}_{^L G}(\pone)_{(z_i),\infty;(\la_i),\la_\infty}$,
i.e., by a $^L G$-oper on $\pone$ with regular singularities at
$z_1,\ldots,z_N,\infty$, with residues
$\la_1,\ldots,\la_N,\la_\infty$ at those points and {\em trivial
monodromy} representation $$\pi_1(\pone \bs \{
z_1,\ldots,z_N,\infty \}) \to {}^L G.$$
\end{prop}

We wish to construct an inverse map from $\on{Op}_{^L
G}(\pone)_{(z_i),\infty;(\la_i),\la_\infty}$ to the spectrum of
$\ZZ(\g)$ on $\left( \bigotimes_{i=1}^N V_{\la_i} \otimes
V_{\la_\infty} \right)^G$.

We start with a more intrinsic geometric description of $\on{Op}_{^L
G}(\pone)_{(z_i),\infty;(\la_i),\la_\infty}$.

\subsection{Geometric description of opers without monodromy}

In order to simplify our notation, in this section we will use the
symbol $G$ instead of $^L G$ to denote our group.

Suppose that we have a $G$-oper $\tau: (\F,\nabla,\F_{B})$ on a
curve $U$. Consider the associated bundle of flag manifolds
$$
({}G/{}B)_{\F} = \F \underset{G}\times {}G/{}B \simeq
\F_B \underset{B}\times {}G/{}B = ({}G/{}B)_{\F_B}.
$$
The reduction $\F_{B}$ gives rise to a section of this bundle.

Now suppose that the monodromy representation $\pi_1(U) \to {}G$
corresponding to this oper is trivial. Picking a point $x \in U$, we
may then use the connection $\nabla$ to identify the fibers of the
oper bundle $\F$ at all points of $U$ with the fiber $\F_x$ at
$x$. Choosing a trivialization $\imath_x: {}G \overset{\sim}\to
\F_x$ of $\F_x$, we then obtain a trivialization of the bundle $\F$,
and hence of the bundle $({}G/{}B)_{\F}$. Let us fix $\imath_x$
and the corresponding trivialization of $\F$. Then the $B$-reduction
$\F_{B}$ gives rise to a map $\phi_\tau: U \to {}G/B$.

The oper condition may be described as follows. Define an open
$B$-orbit ${\bf O}\subset [{}\n,{}\n]^\perp/{}\bb \subset
{}\g/{}\bb$, consisting of vectors which are stabilized by $N\subset
{}B$, and such that all of their negative simple root components, with
respect to the adjoint action of $H = B/N$, are non-zero. This
orbit may also be described as the $B$-orbit of the sum of the
projections of simple root generators $f_i$ of any nilpotent
subalgebra $\n_-$, which is in generic position with $\bb$, onto
$\g/{}\bb$. The torus $H$ acts simply transitively on ${\bf O}$, and
so ${\bf O}$ is an $H$-torsor.

Define an $\ell$-dimensional distribution $T_{\on{Op}} {}G/B$ in
the tangent bundle to $G/{}B$ as follows. The tangent space to
$g{}B$ is identified with the quotient $\g/g{}\bb g^{-1}$. The
fiber of $T_{\on{Op}} {}G/B$ at $g{}B$ is by definition its
subspace $[{}\n,{}\n]^\perp/{}\bb$. It contains an open dense
subset $g {\mb O} g^{-1}$. We denote by
$\overset\circ{T}_{\on{Op}} {}G/B$ the open dense subset of
$T_{\on{Op}} {}G/B$ whose fiber at $g{}B$ is $g {\mb O} g^{-1}$.

A map $\phi: U \to {}G/B$ is said to satisfy the oper condition if
at each point of $U$ the tangent vector to $\phi$ belongs to
$\overset\circ{T}_{\on{Op}} {}G/B$.

For example, if $G = PGL_2$, then the flag manifold is $\pone$, and
a map $\phi: U \to \pone$ satisfies the oper condition if and only if
it has a non-vanishing differential everywhere on $X$.

If $\tau$ is an oper, then the map $\phi_\tau: U \to {}G/B$
satisfies the oper condition (note that
$\overset\circ{T}_{\on{Op}} {}G/B$ is $G$-invariant, and
therefore this condition is independent of the trivialization
$\imath_x$). Conversely, if we are given a map $\phi: U \to {}G/B$
satisfying the oper condition, then the triple
$(\F,\nabla,\F_{B})$, where $\F$ is the trivial $G$-bundle on
$U$, $\nabla$ is the trivial connection, and $\F_{B}$ is the
$B$-reduction of $\F$ defined by the map $\phi$, is a $G$-oper.
Thus, we obtain a bijection between the set of all $G$-opers on
$U$ with trivial monodromy representation and the set of
equivalence classes, with respect to the action of $G$ on $G/B$,
of maps $\phi: U \to {}G/B$ satisfying the oper condition.

Now let $\cla$ be a dominant integral coweight of $G$. A map
$D^\times \to {}H$ is said to vanish at the origin to order $\cla$
if it may be obtained as the composition of an embedding $D^\times
\to C^\times$ at the origin and the cocharacter $\cla: \C^\times
\to {}H$. We will say that a map $\phi: U \to {}G/{}B$ satisfies
the $\cla$-oper condition at $y \in U$ if the restriction of
$\phi$ to the punctured disc $D_y^\times$ satisfies the oper
condition and the corresponding map $\phi_{D_y^\times}: D_y^\times
\to {\mb O} \simeq H$ vanishes to the order $\cla$.

For example, if $G = PGL_2$, then we can identify the dominant
integral weights of $G$ with non-negative integers. Then a map $\phi:
U \to \pone$ satisfies the $\cla$-oper condition at $y \in U$ if its
differential vanishes to the order $\cla$ at $y$.

Recalling the definition of the space $\on{Op}_{G}(D_x)_{\cla}$ given
in \secref{reg sing}, we obtain the following result.

\begin{prop} \label{attach}
There is a bijection between the set
$\on{Op}_{G}(\pone)_{(z_i),\infty;(\cla_i),\cla_\infty}$ and the set
of equivalence classes, with respect to the action of $G$ on $G/B$, of
the maps $\phi: \pone \to {}G/B$ satisfying the oper condition on
$\pone \bs \{ z_1,\ldots,z_N,\infty \})$, the $\cla_i$-oper condition
at $z_i$ and the $\cla_\infty$-oper condition at $\infty$.
\end{prop}

In particular, the points of
$\on{Op}_{PGL_2}(\pone)_{(z_i),\infty;(\cla_i),\cla_\infty}$ are the
same as the equivalence classes of maps $\pone \to \pone$ whose
differential does not vanish anywhere on $\pone \bs \{
z_1,\ldots,z_N,\infty \})$ and vanishes to order $\cla_i$ at $z_i,
i=1,\ldots,N$, and to order $\cla_\infty$ at $\infty$.

\subsection{From opers without monodromy to Bethe vectors}
\label{from}

Let $\tau = (\F,\nabla,\F_{^L B}) \in \on{Op}_{^L
G}(\pone)_{(z_i),\infty;(\la_i),\la_\infty}$. By \propref{attach},
we attach to it a map $\phi_\tau: \pone \to {}^L G/^L B$. A
horizontal $^L B$-reduction on $\tau$ is completely determined by
the choice of $^L B$-reduction at any given point $x \in \pone$,
which is a point of $({}^L G/{}^L B)_{\F_{^L B,x}} \simeq {}^L
G/{}^L B$ (see \lemref{isom with flags}). We will choose the point
$\infty$ as our reference point on $\pone$. Then we obtain

\begin{lem} \label{all miura opers}
For $\tau \in \on{Op}_{^L G}(\pone)_{(z_i),\infty;(\la_i),\la_\infty}$
the set of Miura opers on $\pone$ with the underlying oper $\tau$ is
in bijection with the set of points of $({}^L G/{}^L B)_{\F_{^L
B,\infty}} \simeq {}^L G/{}^L B$.
\end{lem}

Now suppose we are given an oper $\tau \in
\on{Op}_{^L G}(\pone)_{(z_i),\infty;(\la_i),\la_\infty}$. We fix a
trivialization of the fiber $\F_{\infty}$ of $\F$ at $\infty$ and
consider the corresponding map $\phi_\tau: \pone \to {}^L G/{}^L
B$. Then the value of $\phi_\tau$ at $\infty \in \pone$ is the point
$^L B \in {}^L G/{}^L B$. Consider the Miura oper structure on $\tau$
corresponding to $^L B \in {}^L G/{}^L B$, considered as a point in
the fiber $({}^L G/{}^L B)_{\F_{^L B,\infty}}$. This is the unique
Miura oper structure on $\tau$ for which the horizontal Borel
reduction $\F'_{^L B,\infty}$ at the point $\infty$ coincides with the
oper reduction $\F_{^L B,\infty}$ (which is the point $^L B \in {}^L
G/{}^L B$ with respect to our trivialization).

Let us suppose that the oper $\tau$ satisfies the following
conditions:
\begin{itemize}
\item[(1)] $\phi_\tau(z_i)$ is in generic position with $^L B$ for all
$i=1,\ldots,N$;

\item[(2)] the relative position of $\phi_\tau(x)$ and $^L B$ is
either generic or corresponds to a simple reflection $s_i \in W$ for
all $x \in \af = \pone \bs \infty$.
\end{itemize}
Then we will call $\tau$ a {\em non-degenerate} oper.

Let $\tau$ be a non-degenerate oper in $\on{Op}_{^L
G}(\pone)_{(z_i),\infty;(\la_i),\la_\infty}$. Consider the unique
Miura oper $(\F,\nabla,\F_{^L B},\F'_{^L B})$ whose underlying oper is
$\tau$ and such that $\F_{^L B,\infty}$ corresponds to the point $^L B
\in {}^L G/{}^L B \simeq ({}^L G/{}^L B)_{\F_{^L
B,\infty}}$. According to the above conditions, the reductions $\F_{^L
B}$ and $\F'_{^L B}$ are in generic position for all but finitely many
points of $\pone$, which are distinct from
$z_1,\ldots,z_N,\infty$. Let us denote these points by
$w_1,\ldots,w_m$. Then at the point $w_j$ the two reductions have
relative position $s_{i_j}$. According to \propref{map beta}, to this
Miura oper corresponds a connection on the $^L H$-bundle
$\Omega^{\rho}$, and hence a connection $\ol\nabla$ on the dual bundle
$\Omega^{-\rho}$ on $\pone \bs \{ z_1,\ldots,z_N,w_1,\ldots,w_m,\infty
\}$. Moreover, by \propref{isom w}, the connection $\ol\nabla$ extends
to a connection with regular singularity on $\Omega^{-\rho}$ on the
entire $\pone$, with residues $\la_i$ at $z_i, i=1,\ldots,N$,
$w_0(\la_\infty+\rho)-\rho$ at $\infty$, and $-\al_{i_j}$ at $w_j,
j=1,\ldots,m$. (Note that in \propref{isom w} we described the
residues of connections of the bundle $\Omega^\rho$. Since now we
consider connections on $\Omega^{-\rho}$ which is dual to
$\Omega^\rho$, when applying \propref{isom w} we need to change the
signs of all residues.)

Therefore, the restriction of $\ol{\nabla}$ to $\pone\bs \infty$ reads
\begin{equation} \label{main conn}
\ol{\nabla} = \pa_t + \sum_{i=1}^N \frac{\la_i}{t-z_i} - \sum_{j=1}^m
\frac{\al_{i_j}}{t-w_j}.
\end{equation}
Using transformation formula \eqref{trans for conn}, we find that the
residue of $\ol\nabla$ at $\infty$ is equal to
$$
- 2 \rho - \sum_{i=1}^N \la_i + \sum_{j=1}^m \al_{i_j}.
$$
But by our assumption this residue is equal to
$w_0(\la_\infty+\rho)-\rho$. Therefore we find that
\begin{equation} \label{w0la}
-w_0(\la_\infty) = \sum_{i=1}^N \la_i - \sum_{j=1}^m \al_{i_j}.
\end{equation}

In addition, the condition that the oper $\tau$ has trivial
monodromy representation and \lemref{si} imply that if
$$
\pa_t - \al_{i_j}/(t-w_j) + {\mb u}_j(t-w_j), \qquad {\mb u}_j(u)
\in {}^L \h[[u]],
$$
is the expansion of the connection \eqref{main conn} at the point
$w_j$, then $\langle \chal_{i_j},{\mb u}_j(0) \rangle = 0$ for all
$j=1,\ldots,m$. When we write them explicitly, we see immediately
that these equations are precisely the Bethe Ansatz equations
\eqref{bethe}!

Thus, we have attached to a non-degenerate oper $\tau \in \on{Op}_{^L
G}(\pone)_{(z_i),\infty;(\la_i),\la_\infty}$ a connection $\ol\nabla =
\pa_t + \la(t)$ on $\Omega^{-\rho}$ and a solution of the Bethe Ansatz
equations. By construction, the Miura transformation of this
connection is the oper $\tau$, and so we obtain the following

\begin{prop} \label{bije}
Suppose that all opers in $\on{Op}_{^L
G}(\pone)_{(z_i),\infty;(\la_i),\la_\infty}$ are
non-degene\-rate. Then there is a bijection between the set
$\on{Op}_{^L G}(\pone)_{(z_i),\infty;(\la_i),\la_\infty}$ and the set
of solutions of Bethe Ansatz equations \eqref{bethe} such that the
weight $\sum_{i=1}^N \la_i - \sum_{j=1}^m \al_{i_j}$ is dominant.
\end{prop}

Thus, given a non-degenerate oper $\tau \in \on{Op}_{^L
G}(\pone)_{(z_i),\infty;(\la_i),\la_\infty}$, we obtain a solution of
Bethe Ansatz equations \eqref{bethe}. We then follow the procedure of
\secref{bethe vect} to assign to the above solution of the Bethe
Ansatz equations an eigenvector of the Gaudin algebra $\ZZ(\g)$ in
$\bigotimes_{i=1}^N M_{\la_i}$, which is an $\n_+$-invariant vector
$\phi(w_1^{i_1},\ldots,w_m^{i_m})$ of weight $\sum_{i=1}^N \la_i
-\sum_{j=1}^m \al_{i_j}$. Consider its projection onto
$\bigotimes_{i=1}^N V_{\la_i}$. By formula \eqref{w0la}, the weight of
this vector is
\begin{equation} \label{imp cond}
\sum_{i=1}^N \la_i -\sum_{j=1}^m \al_{i_j} = -w_0(\la_\infty),
\end{equation}
and so it may be viewed as an eigenvector of the Gaudin algebra in the
space
$$
V^G_{\laa,\la_\infty} = \left( \bigotimes_{i=1}^N
V_{\la_i} \otimes V_{\la_\infty} \right)^G.
$$

By \thmref{main}, the eigenvalues of the Gaudin algebra $\ZZ(\g)$
on this vector are encoded by the Miura transformation of the
connection $\ol\nabla$, which is precisely the oper $\tau$ with
which we have started. Therefore we obtain the following
result.

\begin{lem}
The two opers encoding the eigenvalues of the Gaudin algebra on
the Bethe vectors corresponding to two different solutions of the
Bethe Ansatz equations are necessarily different.
\end{lem}

Suppose that all opers in $\on{Op}_{^L
G}(\pone)_{(z_i),\infty;(\la_i),\la_\infty}$ are {\em
non-degenerate}\footnote{It follows from the results of Mukhin and
Varchenko in \cite{MV:new} that for some
$\la_1,\ldots,\la_N,\la_\infty$ it is possible that there exist
degenerate opers in $\on{Op}_{^L
G}(\pone)_{(z_i),\infty;(\la_i),\la_\infty}$ even for generic values
of $z_1,\ldots,z_N$.} and that each of the Bethe vectors is {\em
non-zero}. Then we obtain a map
$$
\on{Op}_{^L G}(\pone)_{(z_i),\infty;(\la_i),\la_\infty} \to
\on{Spec} \ol{\mc Z}_{(z_i),\infty;(\la_i),\la_\infty}(\g)
$$
that is inverse to the map \eqref{inj map} discussed at the
beginning of this section.

Let us summarize our results.

\begin{prop} \label{completeness}
Assume that all $^L G$-opers in $\on{Op}_{^L
G}(\pone)_{(z_i),\infty;\laa,\la_\infty}$ are non-degenerate and that
all Bethe vectors obtained from solutions of the Bethe Ansatz
equations \eqref{bethe} satisfying the condition \eqref{imp cond} are
non-zero. Then there is a bijection between the spectrum of the
generalized Gaudin Hamiltonians on $V^G_{\laa,\la_\infty}$ (not
counting multiplicities) and the set $\on{Op}_{^L
G}(\pone)_{(z_i),\infty;\laa,\la_\infty}$ of $^L G$-opers on $\pone$
with regular singularities at $z_1,\ldots,z_N,\infty$ and residues
$\varpi(\la_1+\rho),\ldots,\varpi(\la_N+\rho)$ and
$\varpi(\la_\infty+\rho)$, respectively, which have trivial monodromy.

Moreover, if in addition the Gaudin Hamiltonians are diagonalizable
and have simple spectrum on $V^G_{\laa,\la_\infty}$, then the Bethe
vectors constitute an eigenbasis of $V^G_{\laa,\la_\infty}$.
\end{prop}

The last statement of \propref{completeness} that the Bethe vectors
constitute an eigenbasis of $V^G_{\laa,\la_\infty}$ is referred to as
the completeness of the Bethe Ansatz.

We note that the completeness of the Bethe Ansatz has been previously
proved for $\g=\sw_2$ and generic values of $z_1,\ldots,z_N$ by
Scherbak and Varchenko \cite{SV}.

We will discuss the degenerate opers and the corresponding
eigenvectors of the Gaudin hamiltonians in \secref{last}.

\subsection{The case of $\sw_2$}

In the case when $\g=\sw_2$, we identify the weights with complex
numbers. The Bethe vector has a simple form in this case. We use the
standard basis $\{ e,h,f \}$ of $\sw_2$. Then
\begin{equation}
\phi(w_1,\ldots,w_m) = f(w_1) \ldots f(w_m) v_{\la_1} \otimes \ldots
v_{\la_N},
\end{equation}
where
$$
f(w) = \sum_{i=1}^N \frac{f^{(i)}}{w-z_i}.
$$
The Bethe Ansatz equations read
$$
\sum_{i=1}^N \frac{\la_i}{w_j-z_i} = \sum_{s \neq j}
\frac{2}{w_j-w_s}.
$$
If these equations are satisfied, then this vector is a highest weight
vector of weight $\la_\infty$ (with respect to the diagonal action),
and it is an eigenvector in $\bigotimes_{i=1}^N M_{\la_i}$ of the
Gaudin algebra
$$
{\mc Z}_{(z_i),\infty;(\la_i),\la_\infty} =
\C[\Xi_i]_{i=1,\ldots,N}/\left( \sum_{i=1}^N \Xi_i \right),
$$
where $\Xi_i$ is the $i$th Gaudin Hamiltonian given by formula
\eqref{the gaudin ham}.

Its eigenvalue is represented by the $PGL_2$-oper
$$
\pa_t^2 - \sum_{i=1}^N \frac{\la_i(\la_i+2)/4}{(t-z_i)^2} - \sum_{i=1}^N
\frac{c_i}{t-z_i} = (\pa_t - \la(t))(\pa_t + \la(t)),
$$
where
$$
\la(t) = \sum_{i=1}^N \frac{\la_i/2}{t-z_i} - \sum_{j=1}^m
\frac{1}{t-w_j}.
$$
In other words, the eigenvalue $c_i$ of $\Xi_i$ on
$\phi(w_1,\ldots,w_m)$ is equal to
\begin{equation} \label{ci}
c_i = \la_i \left( \sum_{j \neq i} \frac{\la_j}{z_i-z_j} -
\sum_{j=1}^m \frac{1}{z_i-w_j} \right).
\end{equation}

For $\g=\sw_2$, all opers in $\on{Op}_{^L
G}(\pone)_{(z_i),\infty;(\la_i),\la_\infty}$ are non-degenerate for
generic values of $z_1,\ldots,z_N$ (for example, this follows from the
results of \cite{SV}).

\subsection{Solutions of the Bethe Ansatz equations and flag
manifolds}

In \secref{from} we considered only those solutions of the Bethe
Ansatz equation for which
\begin{equation} \label{weight}
\sum_{i=1}^N \la_i - \sum_{j=1}^m \al_{i_j}
\end{equation}
is a dominant integral weight. The reason for that was that we
wanted to construct eigenvectors of $\ZZ(\g)$ in
$\bigotimes_{i=1}^N V_{\la_i}$ of weight \eqref{weight}. Since the
Bethe vectors are automatically $\n_+$-invariant, we find that
such a vector can be non-zero in $\bigotimes_{i=1}^N V_{\la_i}$
only if this weight is dominant integral. In this case we may
write it in the form $-w_0(\la_\infty)$, where $\la_\infty$ is
another dominant integral weight. Then the corresponding Bethe
vector gives us an eigenvector of $\ZZ(\g)$ in the space
$$V^G_{\laa,\la_\infty} = \left( \bigotimes_{i=1}^N V_{\la_i}
\otimes V_{\la_\infty} \right)^G.$$ We have seen in \secref{from}
that these solutions of the Bethe Ansatz equations are in
one-to-one correspondence with the Miura opers on $\pone$ such
that the corresponding oper in $\on{Op}_{^L
G}(\pone)_{(z_i),\infty;\laa,\la_\infty}$ is non-degenerate, and such
that the horizontal $^L B$-reduction $\F'_{^L B}$ coincides with
the oper $^L B$-reduction $\F_{^L B}$ at $\infty$.

It is natural to consider more general solutions of the Bethe Ansatz
equations, i.e., those for which the weight \eqref{weight} is not
necessarily dominant integral. The first step is to relate these
solutions to Miura opers.

Suppose that we are given an arbitrary solution of the Bethe Ansatz
equations \eqref{bethe}. By \propref{bij1}, we attach to it the
connection $\ol{\nabla}$ on the $^L H$-bundle $\Omega^{-\rho}$ on
$\pone$ whose restriction to $\pone \bs \infty$ reads
\begin{equation} \label{main conn1}
\ol{\nabla} = \pa_t + \sum_{i=1}^N \frac{\la_i}{t-z_i} - \sum_{j=1}^m
\frac{\al_{i_j}}{t-w_j}.
\end{equation}

Recall the Miura transformation $\ol{\mb
b}^*_{(z_i),\infty;(\la_i),\la_\infty}$ introduced in formula
\eqref{miura on pone}. Applying this map to $\ol\nabla$, we obtain a
$^L G$-oper $\tau$. It has regular singularities at $z_1,\ldots,z_N$
and $\infty$, and is regular elsewhere on $\pone$ (recall that the
Bethe Ansatz equations ensure that the oper $\tau$ is regular at the
points $w_1,\ldots,w_m$). Around the point $z_i$ the oper connection
reads
$$\nabla = \pa_t + p_{-1} + \frac{\la_i}{t-z_i} + \on{reg}.$$ and by
applying the gauge transformation with $\la(t)$ we can bring it to the
form \eqref{psi la1}. Therefore, according to \lemref{no mon}, $\tau$
has no monodromy around $z_i, i=1,\ldots,N$.

Let us consider the point $\infty$. Since $\tau$ has no monodromy
anywhere on $\pone \bs \infty$, it cannot have monodromy at $\infty$
either. By \lemref{no mon}, the restriction of $\tau$ to the disc
around $\infty$ belongs to $\on{Op}_G(D_\infty)_{\la_\infty}$ for some
integral dominant weight $\la_\infty$. But then the commutative
diagram \eqref{comm diag} implies that the residue of $\ol{\nabla}$ at
$\infty$ is equal to $w'(\la_\infty+\rho)-\rho$ for some $w' \in
W$. Writing $w'=ww_0$, we find that the residue is equal to
$ww_0(\la_\infty+\rho)-\rho$. On the other hand, using transformation
formula \eqref{trans for conn}, we find that the residue of
$\ol\nabla$ at $\infty$ is equal to
$$
- 2 \rho - \sum_{i=1}^N \la_i + \sum_{j=1}^m \al_{i_j}.
$$
Therefore we must have the equality
\begin{equation} \label{formula}
\sum_{i=1}^N \la_i - \sum_{j=1}^m \al_{i_j} =
w(-w_0(\la_\infty)+\rho)-\rho
\end{equation}
for some dominant integral weight $\la_\infty$ and $w \in W$. Thus, the
oper $\tau$ belongs to $\on{Op}_{^L
G}(\pone)_{(z_i),\infty;(\la_i),\la_\infty}$, and so $\ol{\nabla}$
is an element of $\Con^{\on{gen}}_{(z_i),\infty;(\la_i),\la_\infty}$.

Recall from \secref{funct of coinv} that for any curve $U$ we have
maps
$$
{\mb b}^*_U: \Con_U \to \on{MOp}_{^L G}(U)_{\on{gen}}, \qquad \ol{\mb
b}^*_U: \Con_U \to \on{Op}_{^L G}(U).
$$
Let $U = \pone \bs \{z_1,\ldots,z_N,w_1,\ldots,w_m,\infty \}$. Then
for any weights $\la_1,\ldots,\la_N,\la_\infty$ we
have a subset $\Con^{\on{gen}}_{(z_i),\infty;(\la_i),\la_\infty}$ of
$\Con_U$ which is mapped by the Miura transformation $\ol{\mb
b}^*_{(z_i),\infty;(\la_i),\la_\infty}$ into $\on{Op}_{^L
  G}^{\on{RS}}(\pone)_{(z_i),\infty;(\la_i),\la_\infty}$ (see formula
\eqref{into}). But as we saw above, if the weights
$\la_1,\ldots,\la_N,\la_\infty$ are dominant integral, then the image
is actually contained in even smaller set $\on{Op}_{^L
  G}(\pone)_{(z_i),\infty;(\la_i),\la_\infty}$of those opers which
have trivial monodromy. Therefore we have a map
$$
\Con^{\on{gen}}_{(z_i),\infty;(\la_i),\la_\infty} \to \on{Op}_{^L
  G}(\pone)_{(z_i),\infty;(\la_i),\la_\infty},
$$
which we also denote by $\ol{\mb
  b}^*_{(z_i),\infty;(\la_i),\la_\infty}$. The restriction of the
corresponding map ${\mb b}^*_U$ to
$\Con^{\on{gen}}_{(z_i),\infty;(\la_i),\la_\infty}$ gives us a map
$$
{\mb b}^*_{(z_i),\infty;(\la_i),\la_\infty}:
\Con^{\on{gen}}_{(z_i),\infty;(\la_i),\la_\infty} \to \on{MOp}_{^L
G}(\pone)_{(z_i),\infty;(\la_i),\la_\infty},
$$
where $\on{MOp}_{^L G}(\pone)_{(z_i),\infty;(\la_i),\la_\infty}$ is the
space of Miura opers on $\pone$ whose underlying opers belong to
$\on{Op}_{^L G}(\pone)_{(z_i),\infty;(\la_i),\la_\infty}$. We will
show that the image of this map is open and dense.

Let
$$\on{MOp}_{^L G}(\pone)^{\on{gen}}_{(z_i),\infty;(\la_i),\la_\infty}
\subset \on{MOp}_{^L G}(\pone)_{(z_i),\infty;(\la_i),\la_\infty}$$ be
the subvariety of those Miura opers whose horizontal $^L B$-reduction
satisfies the conditions (1) and (2) of \secref{from}. We call such
Miura opers {\em non-degenerate}.

Recall from \lemref{all miura opers} that the space of all Miura opers
on $\pone$ corresponding to a fixed $^L G$-oper $\tau \in \on{Op}_{^L
G}(\pone)_{(z_i),\infty;(\la_i),\la_\infty}$ is isomorphic to the flag
variety $^L G/{}^L B$. Non-degenerate Miura opers corresponding to
$\tau$ form a subvariety of $^L G/{}^L B$ which we denote by $(^L
G/{}^L B)_\tau$. It is clear that $({}^L G/{}^L B)_\tau$ is an open
and dense subvariety of $^L G/{}^L B$. Indeed, $({}^L G/{}^L B)_\tau$
is contained in the intersection $U_\tau$ of finitely many open and
dense subsets, namely, the sets of points of $^L G/{}^L B$ which are
in generic relative position with $\phi_\tau(z_i)$ (each is isomorphic
to the big Schubert cell). The complement of $({}^L G/{}^L B)_\tau$ in
$U_\tau$ is a subvariety of codimension one. This subvariety consists
of all points in $^L G/{}^L B$ which are in relative position $w$ with
$\phi_\tau(x), x \in \pone \bs \{ z_1,\ldots,z_N,\infty \}$, where $w$
runs over the subset of $W$ of all elements of length $l(w) \geq
2$. The subvariety of these points for fixed $x$ has codimension two,
and therefore their union, as $x$ moves along the curve $\pone \bs \{
z_1,\ldots,z_N,\infty \}$, has codimension (at least) one.

It follows from our construction that the image of the map ${\mb
b}^*_{(z_i),\infty;(\la_i),\la_\infty}$ is contained in
$\on{MOp}^{\on{gen}}_{^L
G}(\pone)_{(z_i),\infty;(\la_i),\la_\infty}$. Conversely, suppose we
are given a Miura oper in $\on{MOp}^{\on{gen}}_{^L
G}(\pone)_{(z_i),\infty;(\la_i),\la_\infty}$. Then we obtain a
connection on $w_0(\Omega^{\rho}) = \F'_{^L B}/{}^L N$, and hence a
connection on $\Omega^{-\rho}$. It follows from \lemref{si} that this
connection belongs to
$\Con^{\on{gen}}_{(z_i),\infty;(\la_i),\la_\infty}$. Moreover, in the
same way as in the proof of \propref{map beta} we find that we obtain
a map
$$
\on{MOp}^{\on{gen}}_{^L
G}(\pone)_{(z_i),\infty;(\la_i),\la_\infty} \to
\Con^{\on{gen}}_{(z_i),\infty;(\la_i),\la_\infty}
$$
that is inverse to ${\mb b}^*_{(z_i),\infty;(\la_i),\la_\infty}$.

Therefore we obtain the following result.

\begin{lem} \label{sets up}
The map ${\mb b}^*_{(z_i),\infty;(\la_i),\la_\infty}$ is a
bijection between $\Con^{\on{gen}}_{(z_i),\infty;(\la_i),\la_\infty}$ and
$\on{MOp}^{\on{gen}}_{^L
G}(\pone)_{(z_i),\infty;(\la_i),\la_\infty}$.
\end{lem}

Combining \lemref{sets up} and \propref{attach}, we
obtain

\begin{thm} \label{another bijection}
The set of solutions of the Bethe Ansatz equations \eqref{bethe} is in
bijection with the set of points of
$\on{MOp}_{^L G}(\pone)^{\on{gen}}_{(z_i),\infty;(\la_i),\la_\infty}$,
where $\la_\infty$ satisfies \eqref{formula} for some $w \in W$.
\end{thm}

By \propref{isom w}, the points of $\on{MOp}_{^L
G}(\pone)^{\on{gen}}_{(z_i),\infty;(\la_i),\la_\infty}$ which satisfy
formula \eqref{formula} with a fixed $w \in W$ are those Miura opers
for which the value of $\phi_\tau$ at $\infty$ belongs to the Schubert
cell
$$
S_{ww_0} = {}^L B w_0 w w_0 {}^L B \subset {}^L G/{}^L B.
$$
Therefore such Miura opers correspond to the points of intersection of
$({}^L G/{}^L B)_\tau$ and $S_{ww_0}$. 

Note that except for the big cell $S_1$, the intersection between the
Schubert cell $S_{ww_0}$ and the open dense subset $({}^L G/{}^L
B)_\tau \subset {}^L G/{}^L B$ could be either an open dense subset of
$S_{ww_0}$ or empty.\footnote{For example, it follows from the results
of Mukhin and Varchenko in \cite{MV:new} that sometimes this open set
may not contain the one point Schubert cell $S_{w_0} \subset G/B$ even
if we allow $z_1,\ldots,z_N$ to be generic.} We make the statement of
\thmref{another bijection} more precise as follows:

\begin{thm} \label{main1}
The set of solutions of the Bethe Ansatz equations \eqref{bethe}
decomposes into a union of disjoint subsets labeled by $\on{Op}_{^L
G}(\pone)_{(z_i),\infty;(\la_i),\la_\infty}$. The set of points
corresponding to $^L G$-oper $\tau \in\on{Op}_{^L
G}(\pone)_{(z_i),\infty;(\la_i),\la_\infty}$ is in bijection with the
set of points of an open and dense subset $({}^L G/{}^L B)_\tau$ of
the flag variety $^L G/{}^L B$.

Further, each of these solution must satisfy the equation
\eqref{formula} for some $w \in W$, and the solutions which satisfy
this equation with a fixed $w \in W$ are in bijection with an open
subset $S_{ww_0} \cap ({}^L G/{}^L B)_\tau$ of the Schubert cell
$S_{ww_0} = {}^L B w_0 w w_0 {}^L B$.
\end{thm}

\begin{remark}
In \cite{MV} (resp., \cite{BM}), it was shown, by a method different
from ours, that in the case when $\g$ is of types $A_n, B_n$ or $C_n$
(resp., $G_2$) and $z_1,\ldots,z_N$ are generic, the set of solutions
of the Bethe Ansatz equations satisfying \eqref{formula} is in
bijection with a disjoint union of open and dense subsets of $S_w$.
The connection between the results of \cite{MV} and our results is
explained in \cite{F:opers}.\qed
\end{remark}

In \cite{F:opers} we obtained a generalization of \thmref{main1} to
the case when $\g$ is an arbitrary Kac-Moody algebra.

\subsection{Degenerate opers}    \label{last}

By \lemref{sets up} we have a bijection between the set
$\Con^{\on{gen}}_{(z_i),\infty;(\la_i),\la_\infty}$ and an open subset
$\on{MOp}^{\on{gen}}_{^L G}(\pone)_{(z_i),\infty;(\la_i),\la_\infty}$
of non-degenerate Miura opers in $$\on{MOp}_{^L
G}(\pone)_{(z_i),\infty;(\la_i),\la_\infty}.$$ Now we wish to extend
this bijection to the entire set $\on{MOp}_{^L
G}(\pone)_{(z_i),\infty;(\la_i),\la_\infty}$. To any point of
$\on{MOp}_{^L G}(\pone)_{(z_i),\infty;(\la_i),\la_\infty}$ we can
still assign, as before, a connection $\ol\nabla$ on $\Omega^{-\rho}$
with regular singularities at $z_1,\ldots,z_N,\infty$ and some
additional points $w_1,\ldots,w_m$. By \propref{isom w}, the residue
of $\ol\nabla$ at $z_i$ (resp., $\infty$, $w_j$) must be equal to
$y_i(\cla_i+\crho)-\crho$ (resp.,
$y_\infty(\cla_\infty+\crho)-\crho, y'_j(\crho)-\crho$) for some
elements $y_i,y_\infty,y'_j \in W$. Hence this connection has the form
\begin{equation}    \label{new conn RS}
\pa_t + \sum_{i=1}^N \frac{y_i(\la_i+\rho)-\rho}{t-z_i}
+ \sum_{j=1}^m \frac{y'_j(\rho)-\rho}{t-w_j}.
\end{equation}
Considering the expansion of this connection at $\infty$, we obtain
the following relation
\begin{equation}    \label{relation}
\sum_{i=1}^N (y_i(\cla_i+\crho)-\crho) + \sum_{j=1}^m
(y'_j(\crho)-\crho) = y_\infty w_0(-w_0(\cla_\infty) +
\crho) - \crho.
\end{equation}

Let $\Con_{(z_i),\infty;(\la_i),\la_\infty}$ be the set of all
connections of the form \eqref{new conn RS} whose Miura transformation
belongs to $\on{Op}_{^L
G}(\pone)_{(z_i),\infty;(\la_i),\la_\infty}$. Using \lemref{all miura
opers}, we then obtain the following generalization of \lemref{sets
up} (see \cite{F:opers}):

\begin{thm}    \label{strongest}
There is a bijection between the sets
$\Con_{(z_i),\infty;(\la_i),\la_\infty}$ and $\on{MOp}_{^L
G}(\pone)_{(z_i),\infty;(\la_i),\la_\infty}$. In particular, the set
of those connections in \linebreak
$\Con_{(z_i),\infty;(\la_i),\la_\infty}$ which correspond to a fixed
$^L G$--oper $\tau \in \on{Op}_{^L
G}(\pone)_{(z_i),\infty;(\cla_i),\cla_\infty}$ is isomorphic to the
set of points of the flag variety $^L G/{}^L B$.

Moreover, the residues of these connections must satisfy the relation
\eqref{relation} for some $y_\infty \in W$. The set of those
connections which satisfy this relation is in bijection with the
Schubert cell $^L B w_0 y_\infty {}^L B$ in $^L G/{}^L B$.
\end{thm}

Consider, in particular, the unique Miura oper corresponding to an
oper $$\tau \in \on{Op}_{^L G}(\pone)_{(z_i),\infty;\laa,\la_\infty}$$
in which the two Borel reductions coincide at $\infty$ (it corresponds
to the one point Schubert cell $^L B$ in $^L G/{}^L B$). Let
$\ol\nabla_\tau$ be the corresponding connection in
$\Con_{(z_i),\infty;(\la_i),\la_\infty}$. It has the form \eqref{new
conn RS} and its residues satisfy the relation \eqref{relation} with
$y_\infty = w_0$:
$$
\sum_{i=1}^N (y_i(\la_i+\rho)-\rho) + \sum_{j=1}^m
(y'_j(\rho)-\rho) = -w_0(\la_\infty).
$$
The oper $\tau \in \on{Op}_{^L
G}(\pone)_{(z_i),\infty;\laa,\la_\infty}$ is degenerate if and only
if either some of the elements $y_i$ are not equal to $1$ or some of
the elements $y'_j$ have lengths greater than $1$ (i.e., are not
simple reflections).

Observe that we can write $y_i(\la_i+\rho)-\rho$ as $\la_i$ minus a
combination of simple roots with non-negative coefficients, and
likewise, $y'_j(\rho)-\rho$ is equal to minus a linear combination of
the simple roots with non-negative coefficients. It is instructive to
think of the connection $\ol\nabla_\tau$ corresponding to degenerate
opers as the limit of a family of connections corresponding to
non-degenerate opers, under which some of the points $w_j$'s either
collide with some of the $z_i$'s or with each other (leading to the
degeneracies of types (1) and (2), respectively). The resulting
residues of the connection become the sums of the residues of the
points that have coalesced together. \propref{isom w} places severe
restrictions on what kinds of residues (and hence collisions) can
occur: namely, at the point $z_i$ the residue must lie in the
$W$--orbit of $\la_i$, whereas at the additional points $w_j$ the
residues must be in the $W$--orbit of $0$ (under the $\rho$--shifted
action of $W$).

We expect that if $z_1,\ldots,z_N$ are generic, then for any $\tau \in
\on{Op}_{^L G}(\pone)_{(z_i),\infty;\laa,\la_\infty}$ we have $y_i =
1$ for all $i=1,\ldots,N$ in the connection $\ol\nabla_\tau$. In other
words, this Miura oper satisfies condition (1) from \secref{from}, but
may not satisfy condition (2), that is at least one of the $y'_j$'s
is not a simple reflection.

Suppose that $\tau$ is such an oper, so in particular we have
\begin{equation}    \label{special condition}
\sum_{i=1}^N \la_i + \sum_{j=1}^m (y'_j(\rho)-\rho) =
-w_0(\la_\infty).
\end{equation}

Then we can still attach to the connection \eqref{new conn RS} an
eigenvector of the generalized Gaudin hamiltonians in
$V_{\laa,\la_\infty}^G$ with eigenvalue $\tau$ by generalizing the
procedure of \secref{bethe vect}.\footnote{the construction outlined
below was found by B. Feigin and myself around the time \cite{FFR} was
written; it was explained in \cite{F:icmp}, Sect. 5.6--5.8 in the case
when $\g=\sw_3$}

More precisely, suppose that the residue of $\la(t)$ at $w_j$ is equal to
$y'_j(\crho)-\crho$ for some $y'_j \in W$. Then the expansion of this
connection at $w_j$ is $\mu_j(t-w_j)$, where
$$
\mu_j(t) = \frac{y'_j(\crho)-\crho}{t} + \sum_{n\geq 0} \mu_{j,n}
t^n, \qquad \mu_{j,n} \in \h^*.
$$
The condition that the oper $\tau$ has no monodromy at $w_j$ (i.e.,
that the Miura transformation of the connection $\pa_t + \mu_j(t)$ is
regular at $t=0$) translates into a system of equations on the
coefficients $\mu_{j,n}$ of this expansion. We denote this system by
$S_j$. Recall from \lemref{si} that when $y'_j = s_{i_j}$ there is
only one equation $\langle \chal_{i_j},\mu_{j,0} \rangle = 0$.

To see what these equations look like, we consider the example when
$\g = \sw_n$. The $PGL_n$--opers may be represented by $n$th order
differential operators of the form
\begin{equation}    \label{first time opers}
\partial_t^n+v_1(t) \partial_t^{n-2}+\ldots+v_{n-1}(t).
\end{equation}
Let us identify the dual Cartan subalgebra of $\sw_n$ with the
hyperplane $\sum_{k=1}^n \epsilon_k = 0$ of the vector space
$\on{span} \{ \epsilon_k \}_{k=1,\ldots,n}$. Then an $^L \h$--valued
connection has the form $\pa_t + \sum_{k=1}^n u_k(t) \epsilon_k$. The
Miura transformation of this connection is given by the formula
\begin{equation}    \label{miura trans for sln}
\partial_t^n+v_1(t) \partial_t^{n-2}+\ldots+v_{n-1}(t) =
(\pa_t - u_1(t)) \ldots (\pa_t - u_n(t)).
\end{equation}
The system of equations $S_j$ is obtained by writing $\mu_j(t)$ as
$\sum_{k=1}^n u_k(t) \epsilon_k$, substituting the functions $u_k(t)$
into formula \eqref{miura trans for sln} and setting to zero all
coefficients in front of the negative powers of $t$ in the resulting
$n$th order differential operator. For example, if $y'_j = s_a s_b,
|a-b|>1$, then we have
$$
u_m(t) = - (\delta_{m,a} - \delta_{m,a+1} + \delta_{m,b} -
\delta_{m,b+1}) t^{-1} + \sum_{r \geq 0} u_{m,r} t^r,
$$
and it is easy to write down the equations on the coefficients
$u_{m,r}$ corresponding to the regularity of the operator \eqref{miura
trans for sln} at $t=0$.

The system of equations $S_j$ also has a nice interpretation in the
theory of the Wakimoto modules. Namely, consider the Wakimoto module
$W_{\mu_j(z)}$. Then it contains a $\g[[t]]$--invariant vector $P_j$
if and only if the system $S_j$ of equations on $\mu_j(z)$ is
satisfied. In the case when $y'_j = s_{i_j}$ and $\mu_j(z) = -
\al_{i_j} z^{-1} + \ldots$, this is the statement of \lemref{sing};
this vector is given by the formula $e^R_{i_j,-1} \vac$ in this
case. In general the formulas for these vectors are more
complicated. They can probably be obtained as the limits of the Bethe
vectors corresponding to small deformations of the highest weights
$\la_i$.

It follows from the definition that for our connection $\ol\nabla_\tau$
the systems $S_j$ are satisfied for all $j=1,\ldots,m$. Therefore each
module $W_{\mu_j(z)}$ contains a $\g[[t]]$--invariant vector
$P_j$. Then, in the same way as in \secref{bethe vect}, we obtain a
homomorphism of $\G_{\ka_c}$-modules
$$
\bigotimes_{i=1}^N W_{\la_i(z)} \otimes \V_0^{\otimes m} \otimes
W'_{\la_\infty(z)} \to \bigotimes_{i=1}^N
W_{\la_i(z)} \otimes \bigotimes_{j=1}^m W_{\mu_j(z)} \otimes
W'_{\la_\infty(z)},
$$
which sends the vacuum vector in the $j$th copy of $\V_0$ to $P_j \in
W_{\mu_j(z)}$. Therefore, composing the corresponding map of the
spaces of coinvariants with the functional $\tau_{(z_i),(w_j)}$
introduced at the end of \secref{the case of pone}, we obtain a linear
functional
$$
\wt\tau_{(z_i),(w_j)}: H((W_{\la_i(z)}),W'_{\la_\infty(z)}) \to \C.
$$
By construction, $\wt\tau_{(z_i),(w_j)}$ is an eigenvector of the
Gaudin algebra $\ZZ(\g)$, and the corresponding eigenvalue is given by
the $^L G$-oper $\ol{\mb
b}^*_{(z_i),\infty;(\la_i),\la_\infty}(\ol\nabla_\tau) = \tau$.

Continuing along the lines of the construction presented in
\secref{bethe vect}, we obtain an $\n_+$--invariant vector in
$\otimes_{i=1}^N M_{\la_i}$ of weight $\sum_{i=1}^N \la_i +
\sum_{j=1}^m (y'_j(\rho)-\rho)$, which, according to formula
\eqref{special condition}, is equal to $-w_0(\la_\infty)$ which is a
dominant integral weight. The projection of this vector onto
$\otimes_{i=1}^N V_{\la_i}$ gives rise to an eigenvector of the Gaudin
algebra $\ZZ(\g)$ in $V_{\laa,\la_\infty}^G$ with the eigenvalue
corresponding to our degenerate oper $\tau$.

Thus, we have assigned to any oper $\tau$ in $\on{Op}_{^L
G}(\pone)_{(z_i),\infty;\laa,\la_\infty}$ such that the corresponding
special connection $\ol\nabla_\tau$ satisfies condition (1) that $y_i
= 1$, for all $i=1,\ldots,N$, in \eqref{new conn RS}, an eigenvector
of the Gaudin algebra in $V_{\laa,\la_\infty}^G$ with the eigenvalue
$\tau$.

The next question is what happens if $\ol\nabla_\tau$ has $y_i \neq 1$
for some $i$. In this case the above construction gives rise to an
eigenvector of the Gaudin algebra $\ZZ(\g)$ in the space
$$(\otimes_{i=1}^N
M_{y_i(\la_i+\rho)-\rho})^{\n_+}_{-w_0(\la_\infty)}$$ of the diagonal
$\n_+$--invariants of $\otimes_{i=1}^N M_{y_i(\la_i+\rho)-\rho}$ of
weight $-w_0(\la_\infty)$.

Now recall that for any dominant integral weight $\la$ and $y \in W$
the Verma module $M_{y(\la+\rho)-\rho}$ is a $\g$--submodule of
$M_\la$. Hence $\otimes_{i=1}^N M_{y_i(\la_i+\rho)-\rho}$ is naturally
a subspace of $\otimes_{i=1}^N M_{\la_i}$. It is also a
$\g_{\on{diag}}$--submodule of $\otimes_{i=1}^N M_{\la_i}$ which is
contained in the maximal proper $\g_{\on{diag}}$--submodule
$I_{(\la_i)}$. The quotient of $\otimes_{i=1}^N M_{\la_i}$ by
$I_{(\la_i)}$ is isomorphic to $\otimes_{i=1}^N V_{\la_i}$. We have
thus associated to $\ol\nabla_\tau$ an eigenvector of the Gaudin
algebra in $\otimes_{i=1}^N M_{\la_i}$ with the eigenvalue $\tau$, but
this eigenvector is inside its maximal proper
$\g_{\on{diag}}$--submodule $I_{(\la_i)}$. Some sample computations
that we have made in the case when $\g=\sw_2$ suggest that in this
case the generalized eigenspace of the Gaudin algebra in
$(\otimes_{i=1}^N M_{\la_i})^{\n_+}_{-w_0(\la_\infty)}$ corresponding
to the eigenvalue $\tau$ has dimension greater than one and that its
projection onto $(\otimes_{i=1}^N V_{\la_i})^{\n_+}_{-w_0(\la_\infty)}
\simeq V_{\laa,\la_\infty}^G$ is non-zero. In other words, we
expect that in this case there also exists an eigenvector of the
Gaudin algebra $\ZZ(\g)$ in $V_{\laa,\la_\infty}^G$ with the same
eigenvalue $\tau$.

If true, this would explain the meaning of the degeneracy of the
connection $\ol\nabla_\tau$ at the point $z_i$. Namely, for some
special values of $z_i$ it may happen that some of the eigenvalues of
the Gaudin operators in $(\otimes_{i=1}^N
V_{\la_i})^{\n_+}_{-w_0(\la_\infty)}$ and in
$(I_{(\la_i)})^{\n_+}_{-w_0(\la_\infty)}$ become equal. Then they may
combine into a Jordan block such that the eigenvector in
$(\otimes_{i=1}^N V_{\la_i})^{\n_+}_{-w_0(\la_\infty)}$ is only a
generalized eigenvector in $(\otimes_{i=1}^N
M_{\la_i})^{\n_+}_{-w_0(\la_\infty)}$. If that happens, we can no
longer obtain this eigenvector by projection from $(\otimes_{i=1}^N
M_{\la_i})^{\n_+}_{-w_0(\la_\infty)}$ (but we could potentially obtain
it by considering the family of eigenvectors corresponding to small
perturbations of the $z_i$'s).

Let $\otimes_{i=1}^N M_{y_i(\la_i+\rho)-\rho}$ be the smallest
$\g^{\otimes N}$--submodule of $\otimes_{i=1}^N M_{\la_i}$ in which the
corresponding ``true'' eigenvector in $(\otimes_{i=1}^N
M_{\la_i})^{\n_+}_{-w_0(\la_\infty)}$ is contained. We believe that
this is precisely the situation when the connection $\ol\nabla_\tau$
develops a singularity with residue $y_i(\la_i+\rho)-\rho$ at the
point $z_i$.

Let us summarize the emerging picture. According to \propref{inj
map1}, we have an injective map from the spectrum of the Gaudin
algebra $\ZZ(\g)$ occurring in $V_{\laa,\la_\infty}^G$ to $\on{Op}_{^L
G}(\pone)_{(z_i),\infty;(\la_i),\la_\infty}$. We wish to construct the
inverse map, in other words, to assign to each oper $\tau \in
\on{Op}_{^L G}(\pone)_{(z_i),\infty;(\la_i),\la_\infty}$ an
eigenvector of $\ZZ(\g)$ with eigenvalue $\tau$ in
$V_{\laa,\la_\infty}^G$. (Here we will not discuss the question as to
whether $\ZZ(\g)$ acts diagonally on $V_{\laa,\la_\infty}^G$ and
whether there are multiplicities. It is known that Jordan blocks occur
even in the case when $\g=\sw_2$, though it is expected that for
generic $z_i$'s the action of $\ZZ(\g)$ is diagonalizable and has
simple spectrum. Here we will only discuss the spectrum as a set
without multiplicities.)

First, we associate to $\tau$ the Miura oper in which the horizontal
Borel reduction coincides with the oper reduction at the point
$\infty$. This Miura oper gives rise to a connection $\ol\nabla_\tau
\in \Con_{(z_i),\infty;(\la_i),\la_\infty}$. It has the form
\eqref{new conn RS}, and the residues satisfy the condition
\eqref{special condition}.

The simplest situation occurs if all elements $y_i$ in \eqref{new conn
RS} are equal to $1$ and all elements $y'_j$ are simple
reflections. This is the situation where Bethe Ansatz is
applicable. Namely, as explained in \secref{bethe vect}, we assign to
$\ol\nabla_\tau$ a certain coinvariant of the tensor product of the
Wakimoto modules, which gives rise to an eigenvector of $\ZZ(\g)$ with
eigenvalue $\tau$ in $V_{\laa,\la_\infty}^G$. This is the Bethe vector
given by an explicit formula \eqref{genbv}. It was believed that this
was in fact a generic situation, i.e., for fixed
$\la_1,\ldots,\la_N,\la_\infty$ and generic $z_1,\ldots,z_N$ the
connection $\ol\nabla_\tau$ satisfies these conditions for all opers
$\tau \in \on{Op}_{^L
G}(\pone)_{(z_i),\infty;(\la_i),\la_\infty}$. But recent results of
Mukhin and Varchenko \cite{MV:new} show that this is not the
case. Hence we analyze the possible degeneracies.

The first type of degeneracy that may occur is the following: while
all elements $y_i$ in \eqref{new conn RS} are still equal to $1$, some
of the elements $y'_j$ are no longer simple reflections.  We expect
this to be the generic situation. In this case we can still construct
an eigenvector of $\ZZ(\g)$ with eigenvalue $\tau$ in
$V_{\laa,\la_\infty}^G$, as explained above. A formula for this vector
will be more complicated than the formula for a Bethe vector, but in
principal such a vector can be obtained by an algorithmic procedure.

The most general and most difficult case is when some of the elements
$y_i$ are not equal to $1$ and some of the elements $y'_j$ are not
simple reflections. In this case, as explained above, we still believe
that there exists an eigenvector of $\ZZ(\g)$ with eigenvalue $\tau$
in $V_{\laa,\la_\infty}^G$. But this eigenvector cannot be constructed
directly, by using the Wakimoto modules, as above. Indeed, our
construction using the Wakimoto modules gives eigenvectors in
$\otimes_{i=1}^N M_{\la_i}$, but we expect that in this most
degenerate case the eigenvector in $V_{\laa,\la_\infty}^G$ cannot be
lifted to an eigenvector in $\otimes_{i=1}^N M_{\la_i}$ (only to a
generalized eigenvector). However, one can probably obtain this
eigenvector by considering families of eigenvectors corresponding to
small perturbations of the $z_i$'s or $\la_i$'s.

Finally, there is a question as to whether all of the eigenvectors
constructed by means of the Wakimoto module construction (including
the Bethe vectors) are non-zero. It was conjectured by S. Chmutov and
I. Scherbak in \cite{CS} that the Bethe vectors are always non-zero
for $\g=\sw_n$.

Now we propose the following conjecture:

\begin{conj}    \label{conjecture}
For any set of integral dominant weights
$\la_1,\ldots,\la_N,\la_\infty$ and an arbitrary collection of
distinct complex numbers $z_1,\ldots,z_N$ there is a bijection between
the set $\on{Op}_{^L G}(\pone)_{(z_i),\infty;(\la_i),\la_\infty}$ and
the spectrum of the Gaudin algebra $\ZZ(\g)$ in
$V_{\laa,\la_\infty}^G$ (not counting multiplicities).
\end{conj}

We do not know whether one can prove this conjecture by following the
strategy outlined above. Rather, one should view the above discussion
only as an explanation of the link between this conjecture and the
traditional Bethe Ansatz conjectures. What we have in mind is a
completely different approach to proving this conjecture that is
suggested by the geometric Langlands correspondence. We conclude this
section by briefly outlining this approach.

As explained in \cite{F:icmp} (following the work \cite{BD} of
A. Beilinson and V. Drinfeld), to each $^L G$--oper $\tau \in
\on{Op}_{^L G}(\pone)_{(z_i),\infty;(\la_i),\la_\infty}$ one assigns a
$D$--module $\wt{\F}_\tau$ on the moduli stack
$\on{Bun}_G(\pone)_{(z_i),\infty}$ of $G$--bundles on $\pone$ with
$B$--reductions at the points $z_1,\ldots,z_N,\infty$. This
$D$--module is a Hecke eigensheaf whose ``eigenvalue'' is the flat $^L
G$--bundle on $\pone \bs \{z_1,\ldots,z_N,\infty \}$ obtained by
forgetting the $^L B$--reduction of $\tau$. This flat $^L G$--bundle
has trivial monodromy and hence is isomorphic to the restriction to
$\pone \bs \{z_1,\ldots,z_N,\infty \}$ of the trivial flat $^L
G$--bundle on $\pone$.

On the other hand, one can show that the $D$--module that we construct
on $\on{Bun}_G(\pone)_{(z_i),\infty}$ ``does not depend'' on the
$B$--reductions at the points $z_i$ and $\infty$, i.e., it is a
pull-back from a $D$--module $\F_\tau$ on the moduli stack
$\on{Bun}_G(\pone)$ of $G$--bundles on $\pone$. This $D$--module is
a Hecke eigensheaf whose eigenvalue is the trivial $^L G$--bundle
with the trivial connection.

The moduli stack $\on{Bun}_G(\pone)$ has an open part
$\on{Bun}^\circ_G(\pone)$ corresponding to the trivial bundle, whose
preimage in $\on{Bun}_G(\pone)_{(z_i),\infty}$ is isomorphic to
$(G/B)^{N+1}/G_{\on{diag}}$.

We expect that the following two properties hold: (a) the Hecke
eigensheaf $\F_\tau$ is non-zero for any $\tau \in \on{Op}_{^L
G}(\pone)_{(z_i),\infty;(\la_i),\la_\infty}$, and (b) the restriction
of any non-zero Hecke eigensheaf on $\on{Bun}_G(\pone)$ to
$\on{Bun}^\circ_G(\pone)$ is necessarily also non-zero.

Suppose that these properties hold (in the case of $\sw_2$ property
(a) can be proved using the Sklyanin separation of variables discussed
in \cite{F:icmp}). Then we obtain that the restriction of the
$D$--module $\wt\F_\tau$ on $\on{Bun}_G(\pone)_{(z_i),\infty}$ to
$(G/B)^{N+1}/G_{\on{diag}}$ is a pull-back of a non-zero $D$--module
on $\on{Bun}^\circ_G(\pone)$, and hence is isomorphic to a direct sum
of copies of its structure sheaf, considered as a $D$--module. This,
in turn, is equivalent to the existence of a non-zero eigenvector of
the Gaudin algebra with eigenvalue $\tau$ in $V_{\laa,\la_\infty}^G$,
as explained in \cite{F:icmp}. Hence we obtain a proof of
\conjref{conjecture}. We hope to return to this question in a future
publication.

\subsection*{Acknowledgements}
This paper reviews the results of our previous works
\cite{FFR,F:icmp,F:opers}, some of which were obtained jointly with
B. Feigin and N. Reshetikhin. I wish to thank I. Scherbak for useful
discussions.

I thank the organizers of the Workshop ``Infinite-dimensional algebras
and quantum integrable systems'' in Faro in July of 2003 for their
invitation to give a talk on this subject and for encouraging me to
write this review.

\end{document}